\input amstex
\documentstyle{amsppt}
%
\catcode`@=11
\redefine\output@{%
  \def\break{\penalty-\@M}\let\par\endgraf
  \ifodd\pageno\global\hoffset=105pt\else\global\hoffset=8pt\fi  
  \shipout\vbox{%
    \ifplain@
      \let\makeheadline\relax \let\makefootline\relax
    \else
      \iffirstpage@ \global\firstpage@false
        \let\rightheadline\frheadline
        \let\leftheadline\flheadline
      \else
        \ifrunheads@ 
        \else \let\makeheadline\relax
        \fi
      \fi
    \fi
    \makeheadline \pagebody \makefootline}%
  \advancepageno \ifnum\outputpenalty>-\@MM\else\dosupereject\fi
}
\def\Beta{\mathchar"0\hexnumber@\rmfam 42}
\catcode`\@=\active
\nopagenumbers
\def\negskp{\hskip -2pt}

\def\lcm{\operatorname{lcm}}
\def\mult{\operatorname{mult}}
\def\compos{\,\raise 1pt\hbox{$\sssize\circ$} \,}
\def\id{\operatorname{id}}
\accentedsymbol\hatgamma{\kern 2pt\hat{\kern -2pt\gamma}}
\accentedsymbol\checkgamma{\kern 2.5pt\check{\kern -2.5pt\gamma}}
\def\blue#1{#1}

\catcode`#=11\def\diez{#}\catcode`#=6
\catcode`&=11\catcode`&=4
\catcode`_=11\def\podcherkivanie{_}\catcode`_=8
\def\mycite#1{\cite{\blue{#1}}\immediate\special{ps:
     ShrHPSdict begin /ShrBORDERthickness 0 def}}
\def\myciterange#1#2#3#4{\cite{\blue{#2#3#4}}\immediate\special{ps:
     ShrHPSdict begin /ShrBORDERthickness 0 def}}
\def\mytag#1{%
    \tag#1}
\def\mythetag#1{\thetag{\blue{#1}}\immediate\special{ps:
     ShrHPSdict begin /ShrBORDERthickness 0 def}}
\def\myrefno#1{\no#1}
\def\myhref#1#2{\blue{#2}\immediate\special{ps:
     ShrHPSdict begin /ShrBORDERthickness 0 def}}
\def\myEarXivlink{\myhref{http://arXiv.org}{http:/\negskp/arXiv.org}}

\def\mytheorem#1{\csname proclaim\endcsname{Theorem #1}}
\def\mytheoremwithtitle#1#2{\csname proclaim\endcsname{Theorem #1#2}}
\def\mythetheorem#1{\blue{#1}\immediate\special{ps:
     ShrHPSdict begin /ShrBORDERthickness 0 def}}
\def\mylemma#1{\csname proclaim\endcsname{Lemma #1}}
\def\mylemmawithtitle#1#2{\csname proclaim\endcsname{Lemma #1#2}}
\def\mythelemma#1{\blue{#1}\immediate\special{ps:
     ShrHPSdict begin /ShrBORDERthickness 0 def}}
\def\mycorollary#1{\csname proclaim\endcsname{Corollary #1}}

\def\myconjecture#1{\csname proclaim\endcsname{Conjecture #1}}
\def\myconjecturewithtitle#1#2{\csname proclaim\endcsname{Conjecture #1#2}}
\def\mytheconjecture#1{\blue{#1}\immediate\special{ps:
     ShrHPSdict begin /ShrBORDERthickness 0 def}}
\def\mytable#1{Table #1}
\def\mythetable#1{\blue{#1}\immediate\special{ps:
     ShrHPSdict begin /ShrBORDERthickness 0 def}}

\pagewidth{360pt}
\pageheight{606pt}
\topmatter
\title
A note on the third
cuboid\\ conjecture. Part~\uppercase\expandafter{\romannumeral 1}.
\endtitle
\author
Ruslan Sharipov
\endauthor
\address Bashkir State University, 32 Zaki Validi street, 450074 Ufa, Russia
\endaddress
\email\myhref{mailto:r-sharipov\@mail.ru}{r-sharipov\@mail.ru}
\endemail
\abstract
    The problem of finding perfect Euler cuboids or proving their non-existence
is an old unsolved problem in mathematics. The third cuboid conjecture is the last
of the three propositions suggested as intermediate stages in proving the 
non-existence of perfect Euler cuboids. It is associated with a certain 
Diophantine equation of the order 12. In this paper a structural theorem for the 
solutions of this Diophantine equation is proved. 
\endabstract
\endtopmatter
\TagsOnRight
\document

\head
1. Introduction.
\endhead
     Let's denote through $P_{abu}(t)$ the following polynomial 
of the order $12$ depending on three integer parameters $a$, $b$, and $u$:
$$
\gathered
P_{abu}(t)=t^{12}+(6\,u^2-2\,a^2-2\,b^2)\,t^{10}
+(a^4+b^4+u^4+4\,a^2\,u^2+\\
+\,4\,b^2\,u^2-12\,b^2\,a^2)\,t^8+(6\,a^4\,u^2+6\,u^2\,b^4-8\,a^2\,b^2\,u^2\,-\\
-\,2\,u^4\,a^2-2\,u^4\,b^2-2\,a^4\,b^2-2\,b^4\,a^2)\,t^6+(4\,u^2\,b^4\,a^2\,+\\
+\,4\,u^2\,a^4\,b^2-12\,u^4\,a^2\,b^2+u^4\,a^4+u^4\,b^4+a^4\,b^4)\,t^4\,+\\
+\,(6\,a^4\,u^2\,b^4-2\,u^4\,a^4\,b^2-2\,u^4\,a^2\,b^4)\,t^2+u^4\,a^4\,b^4.
\endgathered
\quad
\mytag{1.1}
$$
There are some special cases where the polynomial $P_{abu}(t)$ is reducible 
and explicitly splits into lower order factors. Here are these cases:
$$
\xalignat 3
&\hskip -2em
\text{1) \ }a=b; &&\text{3) \ }b\,u=a^2; &&\text{5) \ }a=u;
\qquad\\
\vspace{-1.5ex}
\mytag{1.2}\\
\vspace{-1.5ex}
&\hskip -2em
\text{2) \ }a=b=u; &&\text{4) \ }a\,u=b^2; &&\text{6) \ }b=u.
\qquad
\endxalignat
$$
The special cases \mythetag{1.2} were studied in \mycite{1}, \mycite{2}, and
\mycite{3}. In a general case other than those listed in \mythetag{1.2} the
polynomial \mythetag{1.1} is described by the following conjecture. 
\myconjecturewithtitle{1.1}{ (third cuboid conjecture)} For any three positive 
coprime integer numbers $a$, $b$, and $u$ such that none of the conditions 
\mythetag{1.2} is satisfied the polynomial \mythetag{1.1} is irreducible in the 
ring $\Bbb Z[t]$. 
\endproclaim
     The subcases 2, 5, and 6 in \mythetag{1.2} are trivial. The subcase 1 leads 
to the first cuboid conjecture. The subcases 3 and 4 lead to the second cuboid
conjecture. The first, the second, and the third cuboid conjectures were introduced
\mycite{1}. They are associated with the problem of constructing a perfect Euler 
cuboid (see \mycite{4} and \myciterange{5}{5}{--}{39} for more details). As for the 
polynomial \mythetag{1.1}, it was derived in \mycite{40}.\par
     Let's write the following equation using the polynomial \mythetag{1.1}:
$$
\hskip -2em
P_{abu}(t)=0.
\mytag{1.3}
$$
The equation \mythetag{1.3} can be understood as a Diophantine equation of 
the order 12 with three integer parameters $a$, $b$, and $u$. The third cuboid 
conjecture~\mytheconjecture{1.1} implies the following theorem.
\mytheorem{1.1} For any three positive coprime integer numbers $a$, $b$, and $u$ 
such that none of the conditions \mythetag{1.2} is satisfied the polynomial 
Diophantine equation \mythetag{1.3} has no integer solutions. 
\endproclaim
     A similar theorem associated with the first cuboid conjecture was formulated 
and proved in \mycite{2}. A similar theorem associated with the second cuboid 
conjecture was formulated in \mycite{3}. However, it is not yet proved.\par
     Being a weaker proposition than the conjecture~\mytheconjecture{1.1}, the 
theorem~\mythetheorem{1.1} in our present case is also rather difficult. Probably
it is equally difficult as the third cuboid conjecture itself. Below in section~6 
we formulate and prove a structural theorem for the solutions of the Diophantine 
equation \mythetag{1.3}. This structural theorem is the main result of the present
paper. 
\head
2. The inversion symmetry.
\endhead
     The polynomial $P_{abu}(t)$ in \mythetag{1.1} possesses some special 
properties. They are expressed by the following formulas which can be 
verified by direct calculations:
$$
\xalignat 2
&\hskip -2em
P_{abu}(t)=P_{bau}(t),
&&P_{abu}(-t)=P_{abu}(t).
\mytag{2.1}
\endxalignat
$$
The first equality \mythetag{2.1} means that the polynomial $P_{abu}(t)$ is 
symmetric with respect to the permutation of the parameters $a$ and $b$. The 
second equality is also a symmetry. It is called parity. This symmetry means 
that that the polynomial $P_{abu}(t)$ is an even function of its argument 
$t$. \par
     Apart from the two symmetries \mythetag{2.1}, the polynomial $P_{abu}(t)$
has a third symmetry which is called the inversion symmetry. Having three 
positive integer numbers $a$, $b$, and $u$, we define the transformation
$$
\hskip -2em
\sigma\!:\,(a,b,u)\longmapsto (\tilde a,\tilde b,\tilde u)
\mytag{2.2}
$$
by means of the following three formulas:
$$
\xalignat 3
&\hskip -2em
\tilde a=\frac{\lcm(a,b,u)}{a},
&&\tilde b=\frac{\lcm(a,b,u)}{b},
&&\tilde u=\frac{\lcm(a,b,u)}{u}.
\quad
\mytag{2.3}
\endxalignat
$$
The numerator $\lcm(a,b,u)$ of the fractions \mythetag{2.3} is the least common 
multiple of the integer numbers $a$, $b$, and $u$. The formula
$$
\pagebreak
\hskip -2em
P_{\sigma(abu)}(t)=\frac{P_{abu}(\lcm(a,b,u)/t)\,t^{12}}{a^4\,b^4\,u^4}
\mytag{2.4}
$$
is written in terms of the transformation \mythetag{2.2}. This formula
is easily verified by means of the direct calculations. The formula
\mythetag{2.4} expresses the inversion symmetry of the polynomial 
\mythetag{1.1}.\par
\head
3. Some prerequisites.
\endhead
\mylemma{3.1} For any three positive integer numbers $a$, $b$, and $u$ the
numbers $\tilde a$, $\tilde b$, and $\tilde u$ produced by applying the
transformation \mythetag{2.2} are coprime.
\endproclaim
\demo{Proof} Let $p_1\,\,\ldots,\,p_n$ be the prime factors of the numbers
$a$, $b$, and $u$. Then we can present $a$, $b$, and $u$ in the following 
way:
$$
\xalignat 3
&\hskip -2em
a=\prod^n_{i=1}p_i^{\,\alpha_i},
&&b=\prod^n_{i=1}p_i^{\,\beta_i},
&&u=\prod^n_{i=1}p_i^{\,\omega_i}.
\mytag{3.1}
\endxalignat
$$
The multiplicities $\alpha_i$, $\beta_i$, and $\omega_i$ in \mythetag{3.1}
obey the inequalities
$$
\xalignat 3
&\hskip -2em
\alpha_i\geqslant 0,
&&\beta_i\geqslant 0,
&&\omega_i\geqslant 0.
\mytag{3.2}
\endxalignat
$$
Using the multiplicities \mythetag{3.2}, we define the integer numbers 
$$
\hskip -2em
\theta_i=\max(\alpha_i,\beta_i,\omega_i).
\mytag{3.3}
$$
Then the least common multiple $Z=\lcm(a,b,u)$ in \mythetag{2.3} is expressed
through the above numbers \mythetag{3.3} in the following way:
$$
\hskip -2em
Z=\lcm(a,b,u)=\prod^n_{i=1}p_i^{\,\theta_i}.
\mytag{3.4}
$$
\par
    Let's substitute the formulas \mythetag{3.1} and \mythetag{3.4} into the
formulas \mythetag{2.3}. As a result we derive the following expressions 
for $\tilde a$, $\tilde b$, and $\tilde u$:
$$
\xalignat 3
&\hskip -2em
\tilde a=\prod^n_{i=1}p_i^{\,\theta_i-\alpha_i},
&&\tilde b=\prod^n_{i=1}p_i^{\,\theta_i-\beta_i},
&&\tilde u=\prod^n_{i=1}p_i^{\,\theta_i-\omega_i}.
\mytag{3.5}
\endxalignat
$$
The greatest common divisor of the numbers $\tilde a$, $\tilde b$, and $\tilde u$
in \mythetag{3.5} is calculated by a formula very similar to \mythetag{3.5}. Indeed,
we have
$$
\hskip -2em
\gcd(\tilde a,\tilde b,\tilde u)=\prod^n_{i=1}p_i^{\,\pi_i},
\mytag{3.6}
$$
where the exponents $\pi_i$ are given by the formula
$$
\hskip -2em
\pi_i=\min(\theta_i-\alpha_i,\theta_i-\beta_i,\theta_i-\omega_i). 
\mytag{3.7}
$$
Comparing \mythetag{3.7} with \mythetag{3.3}, we easily see that $\pi_i=0$ for all
$i=1,\,\ldots,\,n$. Substituting $\pi_i=0$ into \mythetag{3.6}, we derive 
$\gcd(\tilde a,\tilde b,\tilde u)=1$. By definition the equality $\gcd(\tilde a,
\tilde b,\tilde u)=1$ means that the numbers $\tilde a$, $\tilde b$, and $\tilde u$
are coprime. Thus, the lemma~\mythelemma{3.1} is proved.\qed\enddemo
\mylemma{3.2} If three positive integer numbers $a$, $b$, $u$ are coprime and if 
the numbers $\tilde a$, $\tilde b$, $\tilde u$ are produced by applying the
transformation \mythetag{2.2} to $a$, $b$, $u$, then, applying the transformation 
\mythetag{2.2} to $\tilde a$, $\tilde b$, $\tilde u$, we get back the numbers 
$a$, $b$, $u$. 
\endproclaim
\demo{Proof} In order to prove the lemma~\mythelemma{3.2} it is convenient to use
the formulas \mythetag{3.1} for the numbers $a$, $b$, $u$ and the formulas 
\mythetag{3.5} for the numbers $\tilde a$, $\tilde b$, and $\tilde u$. The coprimality
condition for $a$, $b$, $u$ is written as 
$$
\hskip -2em
\gcd(a,b,u)=\prod^n_{i=1}p_i^{\,\xi_i}=1.
\mytag{3.8}
$$
For the exponents $\xi_i$ the equality \mythetag{3.8} yields the formula
$$
\hskip -2em
\xi_i=\min(\alpha_i,\beta_i,\omega_i)=0.
\mytag{3.9}
$$\par
     Let's denote through $\hat a$, $\hat b$, $\hat u$ the numbers obtained by 
applying the transformation \mythetag{2.2} to the numbers $\tilde a$, $\tilde b$, 
$\tilde u$. Then we have
$$
\xalignat 3
&\hskip -2em
\hat a=\frac{\lcm(\tilde a,\tilde b,\tilde u)}{\tilde a},
&&\hat b=\frac{\lcm(\tilde a,\tilde b,\tilde u)}{\tilde b},
&&\hat u=\frac{\lcm(\tilde a,\tilde b,\tilde u)}{\tilde u}.
\quad
\mytag{3.10}
\endxalignat
$$
The numerator of the fractions \mythetag{3.10} is calculated according to the
formula
$$
\hskip -2em
\lcm(\tilde a,\tilde b,\tilde u)=\prod^n_{i=1}p_i^{\,\zeta_i},
\mytag{3.11}
$$
where the exponents $\zeta_i$ are given by the formulas 
$$
\hskip -2em
\zeta_i=\max(\theta_i-\alpha_i,\theta_i-\beta_i,\theta_i-\omega_i)
=\theta_i-\min(\alpha_i,\beta_i,\omega_i). 
\mytag{3.12}
$$
The formula \mythetag{3.12} is derived from \mythetag{3.5}, while $\theta_i$ 
are given by the formula \mythetag{3.3}.\par
     Now, applying \mythetag{3.9} to \mythetag{3.12}, we derive $\zeta_i=\theta_i$.
The rest is to substitute $\zeta_i=\theta_i$ into \mythetag{3.11} and then 
substitute \mythetag{3.11} into \mythetag{3.10}. And finally, applying the
formulas \mythetag{3.5} to the transformed formulas \mythetag{3.10}, we derive
$$
\xalignat 3
&\hskip -2em
\hat a=\prod^n_{i=1}p_i^{\,\alpha_i},
&&\hat b=\prod^n_{i=1}p_i^{\,\beta_i},
&&\hat u=\prod^n_{i=1}p_i^{\,\omega_i}.
\mytag{3.13}
\endxalignat
$$
Comparing \mythetag{3.13} with \mythetag{3.1}, we find that the 
lemma~\mythelemma{3.2} is proved. 
\qed\enddemo
     The lemma~\mythelemma{3.2} means that transformation \mythetag{2.2} acts
as an involution upon coprime triples of positive integer numbers $(a,b,u)$, 
i\.\,e\. we have the equality
$$
\hskip -2em
\sigma^2=\sigma\compos\sigma=\id.
\mytag{3.14}
$$
\mylemma{3.3} Let $a$, $b$, $u$ be three positive coprime integer numbers 
and let $\tilde a$, $\tilde b$, $\tilde u$  be the numbers produced from
$a$, $b$, $u$ by applying the transformation \mythetag{2.2}. If one of
the conditions \mythetag{1.2} is fulfilled for $a$, $b$, $u$, \pagebreak 
then the same condition is fulfilled for $\tilde a$, $\tilde b$, $\tilde u$, 
i\.\,e\. 1) $a=b$ implies $\tilde a=\tilde b$, 2) $a=b=u$ implies 
$\tilde a=\tilde b =\tilde u$, 3) $b\,u=a^2$ implies $\tilde b\,\tilde u
=\tilde a^2$, 4) $a\,u=b^2$ implies $\tilde a\,\tilde u=\tilde b^2$, 5) $a=u$ 
implies $\tilde a=\tilde u$, and finally  6) $b=u$ implies $\tilde b=\tilde u$. 
\endproclaim
\mylemma{3.4} Let $a$, $b$, $u$ be three positive coprime integer numbers 
and let $\tilde a$, $\tilde b$, $\tilde u$  be the numbers produced from
$a$, $b$, $u$ by applying the transformation \mythetag{2.2}. If none of
the conditions \mythetag{1.2} is fulfilled for the numbers $a$, $b$, $u$, 
then none of them is fulfilled for the numbers $\tilde a$, $\tilde b$, 
$\tilde u$. 
\endproclaim
     The lemma~\mythelemma{3.3} is proved by means of direct calculations
with the use of the formula \mythetag{2.3}. The lemma~\mythelemma{3.4} is
immediate from the lemma~\mythelemma{3.3}.\par
     Assume that we have an equation \mythetag{1.3} with the parameters 
$a$, $b$, $u$ satisfying the assumptions of the theorem~\mythetheorem{1.1}. 
Then due to the lemma~\mythelemma{3.4} the equation 
$$
\hskip -2em
P_{\sigma(abu)}(t)=0
\mytag{3.15}
$$
is also an equation of the form \mythetag{1.3} whose parameters satisfy
the assumptions of the theorem~\mythetheorem{1.1}. For this reason and due 
to \mythetag{3.14} the equations \mythetag{1.3} and \mythetag{3.15} ate 
called {\it $\sigma$-conjugate cuboid equations}.\par
\head
4. Integer solutions of $\sigma$-conjugate cuboid equations.
\endhead
     Assume that the polynomial $P_{abu}(t)$ has an integer root $t=A_0$. 
Since $a$, $b$, $u$ are nonzero integers, we have $A_0\neq 0$. Then due to the 
inversion symmetry in \mythetag{2.4} the $\sigma$-conjugate polynomial 
$P_{\sigma(abu)}$ has an integer root $t=B_0$, where
$$
\hskip -2em
B_0=\frac{\lcm(a,b,u)}{A_0}. 
\mytag{4.1}
$$
The integer number $B_0$ in \mythetag{4.1} is also nonzero. Applying the parity 
symmetry from \mythetag{2.1}, we conclude that the polynomial $P_{abu}(t)$ has 
the other integer root\linebreak $t=-A_0$, while $P_{\sigma(abu)}(t)$ has the 
other integer root $t=-B_0$. As a result the polynomials $P_{abu}(t)$ and 
$P_{\sigma(abu)}(t)$ split into factors
$$
\xalignat 2
&\hskip -2em
P_{abu}(t)=(t^2-A_0^2)\ C_{10}(t),
&&P_{\sigma(abu)}(t)=(t^2-B_0^2)\ D_{10}(t)\quad
\mytag{4.2}
\endxalignat
$$
with $A_0>0$ and $B_0>0$. Here $C_{10}(t)$ and $D_{10}(t)$ are tenth order 
polynomials complementary to $t^2-A_0^2$ and $t^2-B_0^2$. Applying 
\mythetag{2.1} to \mythetag{4.2} we derive
$$
\xalignat 2
&\hskip -2em
C_{10}(t)=C_{10}(-t),
&&D_{10}(t)=D_{10}(-t).
\quad
\mytag{4.3}
\endxalignat
$$
Due to \mythetag{4.3} the polynomials $C_8(t)$ and $D_8(t)$ are given by the
formulas
$$
\hskip -2em
\aligned
&C_{10}(t)=t^{10}+C_8\,t^8+C_6\,t^6+C_4\,t^4+C_2\,t^2+C_0,\\
&D_{10}(t)=t^{10}+D_8\,t^8+D_6\,t^6+D_4\,t^4+D_2\,t^2+D_0.
\endaligned
\mytag{4.4}
$$
Now let's apply the inversion symmetries from \mythetag{2.4} to
\mythetag{4.2}. As a result we get
$$
C_{10}(t)=-\frac{D_{10}(\lcm(\tilde a,\tilde b,\tilde u)/t)
\,t^{10}}{\tilde a^4\,\tilde b^{\kern 0.2pt 4}\,\tilde u^4/B_0^2},
\quad
D_{10}(t)=-\frac{C_{10}(\lcm(a,b,u)/t)\,t^{10}}
{a^4\,b^{\kern 0.2pt 4}\,u^4/A_0^2}.
\quad
\mytag{4.5}
$$
Applying the symmetries \mythetag{4.5} to \mythetag{4.4}, we derive a series 
of relationships for the coefficients of the polynomials $C_{10}(t)$ and
$D_{10}(t)$:
$$
\xalignat 2
&\hskip -2em
C_0\,Z^2=-a^4\,b^4\,u^4\,B_0^2,
&C_2\,Z^4&=-a^4\,b^4\,u^4\,B_0^2\,D_8,\\
&\hskip -2em
C_4\,Z^6=-a^4\,b^4\,u^4\,B_0^2\,D_6,
&C_6\,Z^8&=-a^4\,b^4\,u^4\,B_0^2\,D_4,
\mytag{4.6}\\
&\hskip -2em
C_8\,Z^{10}=-a^4\,b^4\,u^4\,B_0^2\,D_2,
&Z^{12}&=-a^4\,b^4\,u^4\,B_0^2\,D_0,\\
\vspace{2ex}
&\hskip -2em
D_0\,Z^2=-\tilde a^4\,\tilde b^4\,\tilde u^4\,A_0^2,
&D_2\,Z^4&=-\tilde a^4\,\tilde b^4\,\tilde u^4\,A_0^2\,C_8,\\
&\hskip -2em
D_4\,Z^6=-\tilde a^4\,\tilde b^4\,\tilde u^4\,A_0^2\,C_6,
&D_6\,Z^8&=-\tilde a^4\,\tilde b^4\,\tilde u^4\,A_0^2\,C_4,
\mytag{4.7}\\
&\hskip -2em
D_8\,Z^{10}=-\tilde a^4\,\tilde b^4\,\tilde u^4\,A_0^2\,C_2,
&Z^{12}&=-\tilde a^4\,\tilde b^4\,\tilde u^4\,A_0^2\,C_0.
\endxalignat
$$
Here we use the notations \mythetag{3.4}, i\.\,e\. the relationships
$$
\hskip -2em
Z=\lcm(a,b,u)=\lcm(\tilde a,\tilde b,\tilde u)
\mytag{4.8}
$$
are fulfilled for the parameter $Z$ in \mythetag{4.6} and \mythetag{4.7}.
The equations \mythetag{4.6} and \mythetag{4.7} are excessive. Due to 
\mythetag{2.3} and \mythetag{4.8} some of them are equivalent to some others. 
For this reason we can eliminate excessive variables:
$$
\xalignat 2
&\hskip -2em
C_0=-\frac{a^4\,b^4\,u^4\,B_0^2}{Z^2\strut},
&&C_2=-\frac{a^4\,b^4\,u^4\,B_0^2\,D_8}{Z^4\strut},\\
&\hskip -2em
C_4=-\frac{a^4\,b^4\,u^4\,B_0^2\,D_6}{Z^6\strut},
&&D_0=-\frac{\tilde a^4\,\tilde b^4\,\tilde u^4\,A_0^2}{Z^2\strut},
\mytag{4.9}\\
&\hskip -2em
D_2=-\frac{\tilde a^4\,\tilde b^4\,\tilde u^4\,A_0^2\,C_8}{Z^4\strut},
&&D_4=-\frac{\tilde a^4\,\tilde b^4\,\tilde u^4\,A_0^2\,C_6}{Z^6\strut}.
\endxalignat
$$
Substituting \mythetag{4.9} into the formulas \mythetag{4.4} for $C_{10}(t)$
and $D_{10}(t)$, we get 
$$
\hskip -2em
\aligned
&\gathered
C_{10}(t)=t^{10}+C_8\,t^8+C_6\,t^6-a^4\,b^4\,u^4\,B_0^2\,D_6\,Z^{-6}\,t^4\,-\\
-\,a^4\,b^4\,u^4\,B_0^2\,D_8\,Z^{-4}\,t^2-a^4\,b^4\,u^4\,B_0^2\,Z^{-2},
\endgathered\\
\vspace{2ex}
&\gathered
D_{10}(t)=t^{10}+D_8\,t^8+D_6\,t^6-\tilde a^4\,\tilde b^4
\,\tilde u^4\,A_0^2\,C_6\,Z^{-6}\,t^4\,-\\
-\,\tilde a^4\,\tilde b^4\,\tilde u^4\,A_0^2\,C_8\,Z^{-4}\,t^2
-\tilde a^4\,\tilde b^4\,\tilde u^4\,A_0^2\,Z^{-2}.
\endgathered
\endaligned
\mytag{4.10}
$$\par
     Having derived the formulas \mythetag{4.10}, we substitute them back into
the relationships \mythetag{4.2}. As a result we derive the following formulas:
$$
\allowdisplaybreaks
\gather
\hskip -2em
\gathered
P_{abu}(t)=t^{12}+(C_8-A_0^2)\,t^{10}+(C_6-A_0^2\,C_8)\,t^8\,-\\
\vspace{2ex}
-\,\frac{A_0^2\,C_6\,Z^6+B_0^2\,a^4\,b^4\,u^4\,D_6}{Z^6}\,t^6
-\,\frac{B_0^2\,a^4\,b^4\,u^4\,(D_8\,Z^2-A_0^2\,D_6)}{Z^6}\,t^4\,-\\
\vspace{1ex}
-\,\frac{B_0^2\,a^4\,b^4\,u^4\,(Z^2-A_0^2\,D_8)}{Z^4}\,t^2
+\frac{A_0^2\,B_0^2\,a^4\,b^4\,u^4}{Z^2},
\endgathered\quad
\mytag{4.11}\\
\vspace{1ex}
\hskip -2em
\gathered
P_{\sigma(abu)}(t)=t^{12}+(D_8-B_0^2)\,t^{10}+(D_6-B_0^2\,D_8)\,t^8\,-\\
\vspace{2ex}
-\,\frac{B_0^2\,D_6\,Z^6+A_0^2\,\tilde a^4\,\tilde b^4\,\tilde u^4\,C_6}
{Z^6}\,t^6-\frac{A_0^2\,\tilde a^4\,\tilde b^4\,\tilde u^4
\,(C_8\,Z^2-B_0^2\,C_6)}{Z^6}\,t^4\,-\\
\vspace{1ex}
-\,\frac{A_0^2\,\tilde a^4\,\tilde b^4\,\tilde u^4\,(Z^2-B_0^2\,C_8)}{Z^4}
\,t^2+\frac{B_0^2\,A_0^2\,\tilde a^4\,\tilde b^4\,\tilde u^4}{Z^2}.
\endgathered\quad
\mytag{4.12}
\endgather
$$
The polynomial $P_{abu}(t)$ in \mythetag{4.11} is initially given by the 
formula \mythetag{1.1}. As for the polynomial $P_{\sigma(abu)}(t)$ it is 
produced from the polynomial \mythetag{1.1} by substituting $\tilde a$, 
$\tilde b$, and $\tilde u$ for the parameters $a$, $b$, and $u$ respectively:
$$
\gathered
P_{\sigma(abu)}(t)=t^{12}+(6\,\tilde u^2-2\,\tilde a^2-2\,\tilde b^2)\,t^{10}
+(\tilde a^4+\tilde b^4+\tilde u^4+4\,\tilde a^2\,\tilde u^2\,+\\
+\,4\,\tilde b^2\,\tilde u^2-12\,\tilde b^2\,\tilde a^2)\,t^8+(6\,\tilde a^4
\,\tilde u^2+6\,\tilde u^2\,\tilde b^4-8\,\tilde a^2\,\tilde b^2\,\tilde u^2\,-\\
-\,2\,\tilde u^4\,\tilde a^2-2\,\tilde u^4\,\tilde b^2-2\,\tilde a^4\,\tilde b^2
-2\,\tilde b^4\,\tilde a^2)\,t^6+(4\,\tilde u^2\,\tilde b^4\,\tilde a^2\,+\\
+\,4\,\tilde u^2\,\tilde a^4\,\tilde b^2-12\,\tilde u^4\,\tilde a^2\,\tilde b^2
+\tilde u^4\,\tilde a^4+\tilde u^4\,\tilde b^4+\tilde a^4\,\tilde b^4)\,t^4\,+\\
+\,(6\,\tilde a^4\,\tilde u^2\,\tilde b^4-2\,\tilde u^4\,\tilde a^4\,\tilde b^2
-2\,\tilde u^4\,\tilde b^4\,\tilde a^2)\,t^2+\tilde u^4\,\tilde a^4\,\tilde b^4.
\endgathered
\quad
\mytag{4.13}
$$
Comparing the formula \mythetag{4.11} with \mythetag{1.1} and comparing 
the formula \mythetag{4.12} with \mythetag{4.13}, we derive twelve equations
for the coefficients of the polynomials \mythetag{4.10}. Two of them are
equivalent to the equation \mythetag{4.1} written as
$$
\hskip -2em
A_0\,B_0=Z. 
\mytag{4.14}
$$ 
The other ten of these equations are written as follows:
$$
\gather
\hskip -2em
\gathered
C_8-A_0^2=6\,u^2-2\,a^2-2\,b^2,\\
D_8-B_0^2=6\,\tilde u^2-2\,\tilde a^2-2\,\tilde b^2,
\endgathered\quad
\mytag{4.15}\\
\vspace{1ex}
\hskip -2em
\gathered
C_6-A_0^2\,C_8=a^4+b^4+u^4+4\,a^2\,u^2+4\,b^2\,u^2-12\,b^2\,a^2,\\
D_6-B_0^2\,D_8=\tilde a^4+\tilde b^4+\tilde u^4
+4\,\tilde a^2\,\tilde u^2+4\,\tilde b^2\,\tilde u^2
-12\,\tilde b^2\,\tilde a^2,
\endgathered\quad
\mytag{4.16}\\
\vspace{1ex}
\hskip -2em
\gathered
-(A_0^2\,C_6\,Z^6+B_0^2\,a^4\,b^4\,u^4\,D_6)\,Z^{-6}=6\,a^4\,u^2+6\,u^2\,b^4\,-\\
-\,8\,a^2\,b^2\,u^2-2\,u^4\,a^2-2\,u^4\,b^2-2\,a^4\,b^2-2\,b^4\,a^2,\\
-(B_0^2\,D_6\,Z^6+A_0^2\,\tilde a^4\,\tilde b^4\,\tilde u^4\,C_6)\,Z^{-6}
=6\,\tilde a^4\,\tilde u^2+6\,\tilde u^2\,\tilde b^4\,-\\
-\,8\,\tilde a^2\,\tilde b^2\,\tilde u^2-2\,\tilde u^4\,\tilde a^2
-2\,\tilde u^4\,\tilde b^2-2\,\tilde a^4\,\tilde b^2-2\,\tilde b^4\,\tilde a^2,
\endgathered\quad
\mytag{4.17}\\
\vspace{1ex}
\hskip -2em
\gathered
B_0^2\,a^4\,b^4\,u^4\,(A_0^2\,D_6-D_8\,Z^2)\,Z^{-6}=4\,u^2\,b^4\,a^2\,+\\
+\,4\,u^2\,a^4\,b^2-12\,u^4\,a^2\,b^2+u^4\,a^4+u^4\,b^4+a^4\,b^4,\\
A_0^2\,\tilde a^4\,\tilde b^4\,\tilde u^4\,(B_0^2\,C_6-C_8\,Z^2)\,Z^{-6}
=4\,\tilde u^2\,\tilde b^4\,\tilde a^2\,+\\
+\,4\,\tilde u^2\,\tilde a^4\,\tilde b^2-12\,\tilde u^4\,\tilde a^2\,\tilde b^2
+\tilde u^4\,\tilde a^4+\tilde u^4\,\tilde b^4+\tilde a^4\,\tilde b^4,
\endgathered\quad
\mytag{4.18}\\
\vspace{1ex}
\hskip -2em
\gathered
B_0^2\,a^4\,b^4\,u^4\,(A_0^2\,D_8-Z^2)\,Z^{-4}=6\,a^4\,u^2\,b^4
-2\,u^4\,a^4\,b^2-2\,u^4\,b^4\,a^2,\\
A_0^2\,\tilde a^4\,\tilde b^4\,\tilde u^4\,(B_0^2\,C_8-Z^2)\,Z^{-4}
=6\,\tilde a^4\,\tilde u^2\,\tilde b^4-2\,\tilde u^4\,\tilde a^4\,\tilde b^2
-2\,\tilde u^4\,\tilde b^4\,\tilde a^2.
\endgathered\quad
\mytag{4.19}
\endgather
$$\par
     Note that the parameters $a$, $b$, $u$ and $\tilde a$, $\tilde b$, 
$\tilde u$ are related with each other by means of the formulas \mythetag{2.3}. 
Now we write these formulas as follows: 
$$
\xalignat 3
&\hskip -2em
a\,\tilde a=Z,
&&a\,\tilde b=Z,
&&u\,\tilde u=Z.
\mytag{4.20}
\endxalignat
$$
The equations \mythetag{4.15}, \mythetag{4.16}, \mythetag{4.17}, 
\mythetag{4.18}, \mythetag{4.19} are excessive. Indeed, the equations
\mythetag{4.18} can be derived from the equations \mythetag{4.16} by applying 
the equalities \mythetag{4.14} and \mythetag{4.20}. Similarly, the equations 
\mythetag{4.19} can be derived from the equations \mythetag{4.15} by applying 
the equalities \mythetag{4.14} and \mythetag{4.20}.
As for the equations \mythetag{4.15}, \mythetag{4.16}, \mythetag{4.17}, when
complemented with the equation \mythetag{4.14}, they constitute a system of 
Diophantine equations with respect to the integer variables $C_8$, $D_8$, 
$C_6$, $D_6$, $A_0$ and $B_0$. The results of the above calculations are 
summarized as a lemma.
\mylemma{4.1} For any three positive coprime integer numbers $a$, $b$, and $u$ 
the polynomial $P_{abu}(t)$ has integer roots if and only if the system of 
Diophantine equations \mythetag{4.14}, \mythetag{4.15}, \mythetag{4.16}, and
\mythetag{4.17} is solvable with respect to the integer variables $C_8$, $D_8$, 
$C_6$, $D_6$, $A_0>0$, and $B_0>0$.
\endproclaim 
     Now let's consider the equations \mythetag{4.17}. They are not independent.
The second equation \mythetag{4.17} can be derived from the first one. Indeed,
it is sufficient to multiply the first equation \mythetag{4.17} by 
$tu^4\,ta^4\,tb^4\,Z^{-6}$ and then apply the relationships \mythetag{4.20}.
Due to this observation we can omit the second equation \mythetag{4.17}
preserving the first equation \mythetag{4.17} only. We write this equation
as follows:
$$
\hskip -2em
\gathered
\tilde a^2\,\tilde b^2\,\tilde u^2\,A_0^2\,C_6+a^2\,b^2\,u^2\,B_0^2\,D_6
=Z^4\,(8\,Z^2-6\,\tilde b^2\,a^2\,-\\
-\,6\,b^2\,\tilde a^2+2\,\tilde a^2\,u^2+2\,\tilde b^2\,u^2
+2\,\tilde u^2\,a^2+2\,\tilde u^2\,b^2).
\endgathered
\mytag{4.21}
$$
The equation \mythetag{4.21} is produced from the first equation \mythetag{4.17}
by multiplying it by $\tilde a^2\,\tilde b^2\,\tilde u^2$ and then applying the
relationships \mythetag{4.20}. In terms of the equation \mythetag{4.21} the above
lemma~\mythelemma{4.1} is reformulated as follows. 
\mylemma{4.2} For any three positive coprime integer numbers $a$, $b$, and $u$ 
the polynomial $P_{abu}(t)$ has integer roots if and only if the system of 
Diophantine equations \mythetag{4.14}, \mythetag{4.15}, \mythetag{4.16}, and
\mythetag{4.21} is solvable with respect to the integer variables $C_8$, $D_8$, 
$C_6$, $D_6$, $A_0>0$, and $B_0>0$.
\endproclaim
     Note that the equations \mythetag{4.15} can be explicitly resolved with 
respect to the variables $C_8$ and $D_8$. As a result we get
$$
\hskip -2em
\gathered
C_8=A_0^2+6\,u^2-2\,a^2-2\,b^2,\\
D_8=B_0^2+6\,\tilde u^2-2\,\tilde a^2-2\,\tilde b^2.
\endgathered
\mytag{4.22}
$$
Similarly, the equations \mythetag{4.16} can be explicitly resolved with respect 
to the variables $C_6$ and $D_6$. Resolving them, we get
$$
\hskip -2em
\gathered
C_6=A_0^2\,C_8+a^4+b^4+u^4+4\,a^2\,u^2+4\,b^2\,u^2-12\,b^2\,a^2,\\
D_6=B_0^2\,D_8+\tilde a^4+\tilde b^4+\tilde u^4
+4\,\tilde a^2\,\tilde u^2+4\,\tilde b^2\,\tilde u^2
-12\,\tilde b^2\,\tilde a^2.
\endgathered\quad
\mytag{4.23}
$$
Then we can substitute the expressions \mythetag{4.22} for $C_8$ and $D_8$ into
the equations \mythetag{4.23}. \pagebreak As a result we get the expressions for 
$C_6$ and $D_6$ directly through $A_0$ and $B_0$. We write these expressions in 
the following way:
$$
\gather
\hskip -2em
\aligned
C_6&=A_0^2\,(A_0^2+6\,u^2-2\,a^2-2\,b^2)
+a^4+\\
&+\,b^4+u^4+4\,a^2\,u^2+4\,b^2\,u^2-12\,b^2\,a^2,
\endaligned
\mytag{4.24}\\
\hskip -2em
\aligned
D_6&=B_0^2\,(B_0^2+6\,\tilde u^2-2\,\tilde a^2-2\,\tilde b^2)
+\tilde a^4+\\
&+\,\tilde b^4+\tilde u^4
+4\,\tilde a^2\,\tilde u^2+4\,\tilde b^2\,\tilde u^2
-12\,\tilde b^2\,\tilde a^2.
\endaligned
\mytag{4.25}
\endgather
$$
The next step is to substitute \mythetag{4.24} and \mythetag{4.25} into the
equation \mythetag{4.21}. As a result we obtain the following equation for
the variables $A_0$ and $B_0$:
$$
\gathered
\tilde a^2\,\tilde b^2\,\tilde u^2\,A_0^2\,\bigl(A_0^2\,(A_0^2+6\,u^2
-2\,a^2-2\,b^2)+a^4+b^4+u^4+4\,a^2\,u^2+\\
+\,4\,b^2\,u^2-12\,b^2\,a^2\bigr)+a^2\,b^2\,u^2\,B_0^2\,
\bigl(B_0^2\,(B_0^2+6\,\tilde u^2-2\,\tilde a^2-2\,\tilde b^2)\,+\\
+\,\tilde a^4+\tilde b^4+\tilde u^4+4\,\tilde a^2\,\tilde u^2
+4\,\tilde b^2\,\tilde u^2-12\,\tilde b^2\,\tilde a^2\bigr)
-Z^4\,(8\,Z^2\,-\\
-\,6\,\tilde b^2\,a^2-6\,b^2\,\tilde a^2+2\,\tilde a^2\,u^2
+2\,\tilde b^2\,u^2+2\,\tilde u^2\,a^2+2\,\tilde u^2\,b^2)=0.
\endgathered
\mytag{4.26}
$$
With the use of the equation \mythetag{4.26} now the lemma~\mythelemma{4.2}
is reformulated as follows.
\mylemma{4.3} For any three positive coprime integer numbers $a$, $b$, and $u$ 
the polynomial $P_{abu}(t)$ has integer roots if and only if the system of 
Diophantine equations \mythetag{4.14} and \mythetag{4.26} is solvable with 
respect to the integer variables $A_0>0$, and $B_0>0$.
\endproclaim
\head
5. The prime factors structure.
\endhead
     Below we continue studying the equations \mythetag{4.14} and \mythetag{4.26}
implicitly assuming $a$, $b$, and $u$ to be three positive coprime integer 
numbers. Assuming $p_1,\,\ldots,\,p_n$ to be the prime factors of $a$, $b$, and 
$u$, we apply the formulas \mythetag{3.1} with the multiplicities $\alpha_i$, 
$\beta_i$, and $\omega_i$ obeying the inequalities \mythetag{3.2}. For the 
least common multiple $Z$ of the numbers $a$, $b$, and $u$ in \mythetag{4.8}
we use the formula \mythetag{3.4}, where the exponents $\theta_i$ are given by
the formula \mythetag{3.3}.\par
     The equation \mythetag{4.14} combined with the formula \mythetag{3.4} means
that the numbers $A_0$ and $B_0$ cannot have prime factors other than $p_1,\,
\ldots,\,p_n$. Therefore we write
$$
\xalignat 2
&\hskip -2em
A_0=p_1^{\kern 0.7pt\mu_1}\cdot\ldots\cdot p_n^{\kern 0.7pt\mu_m},
&&B_0=p_1^{\kern 0.7pt\eta_1}\cdot\ldots\cdot p_n^{\kern 0.7pt\eta_m}.
\mytag{5.1}
\endxalignat
$$
In terms of \mythetag{5.1} and \mythetag{3.3} the equation \mythetag{4.14} yields 
the equalities
$$
\hskip -2em
\mu_i+\eta_i=\theta_i
\mytag{5.2}
$$
for each particular value of the index $i=1,\,\ldots,\,n$. The coprimality condition 
for the numbers $a$, $b$, $u$ is \ $\gcd(a,b,u)=1$. It leads to \mythetag{3.8} and 
\mythetag{3.9}. Due to \mythetag{3.9} at least one of the three options is fulfilled 
for each particular $i=1,\,\ldots,\,n$:
$$
\xalignat 5
&\alpha_i=0,
&&\text{or}
&&\beta_i=0,
&&\text{or}
&&\omega_i=0.
\qquad
\mytag{5.3}
\endxalignat
$$
Note that the multiplicities $\alpha_i$, $\beta_i$, and $\omega_i$ in \mythetag{5.3} 
cannot vanish simultaneously. For this reason, applying the formula \mythetag{3.3}, 
we derive 
$$
\hskip -2em
\theta_i=\max(\alpha_i,\beta_i,\omega_i)>0.
\mytag{5.4}
$$
The inequality \mythetag{5.4} means that the multiplicities $\mu_i$ and $\eta_i$ in
\mythetag{5.2} cannot vanish simultaneously either. 

\par
     In order to investigate the equation \mythetag{4.26} we introduce the notation 
$\mult_p(N)$ for the multiplicity of the prime number $p$ in the prime factors 
expansion of $N$: 
$$
\hskip -2em
\mult_p(N)=k\text{\ \ means \ }N=N'\cdot p^k\text{, \ where \ }
N'\not\equiv 0\ (\kern -0.9em\mod p).
\mytag{5.5}
$$
Assume that an integer number $N$ is a sum of some other integer numbers:
$$
\hskip -2em
N=N_1+\ldots+N_m.
\mytag{5.6}
$$
Let's denote through $k_i=\mult_p(N_i)$ the multiplicities of the summands in
\mythetag{5.6} and denote through $k_{\ssize\text{min}}$ the minimum of these
multiplicities:
$$
\hskip -2em
k_{\ssize\min}=\min(k_1,\ldots,k_m).
\mytag{5.7}
$$
In terms of the notations \mythetag{5.5}, \mythetag{5.6}, and \mythetag{5.7} 
we can formulate the following three simple lemmas. Their proofs are obvious. 
\mylemma{5.1} If exactly one term $N_s$ in the sum \mythetag{5.6} has the minimal
multiplicity $k_s=k_{\ssize\min}$, then the multiplicity of the sum in whole
is equal to this minimal multiplicity, i\.\,e\. $\mult_p(N)=k_s=k_{\ssize\min}$.
\endproclaim
\mylemma{5.2} If more than one term in the sum \mythetag{5.6} has the minimal
multiplicity $k_{\ssize\min}$, then $\mult_p(N)\geqslant k_{\ssize\min}$.
\endproclaim
\mylemma{5.3} If exactly one term $N_s$ in the sum \mythetag{5.6} has the minimal
multiplicity $k_s=k_{\ssize\min}$, then the sum $N$ cannot vanish, i\.\,e\.
$N\neq 0$.
\endproclaim
\noindent Using the notation \mythetag{5.5} and the above three lemmas~\mythelemma{5.1},
\mythelemma{5.2}, and \mythelemma{5.3}, below we study several options derived from 
\mythetag{5.3}.\par
     {\bf The case $\alpha_i>\beta_i>\omega_i=0$ and $p_i\neq 2$}. In this case the
multiplicities of the parameters $a$, $b$, and $u$ obey the following equalities and
inequalities:
$$
\hskip -2em
\mult_{p_i}(a)=\alpha_i>\mult_{p_i}(b)=\beta_i>\mult_{p_i}(u)=\omega_i=0. 
\mytag{5.8}
$$
Applying the formulas \mythetag{3.3}, \mythetag{3.4}, and  \mythetag{4.20} to 
\mythetag{5.8}, we derive 
$$
\xalignat 2
\hskip -2em
\theta_i=&\mult_{p_i}(Z)=\alpha_i,
&\tilde\alpha_i=&\mult_{p_i}(\tilde a)=0,\\
\vspace{-1.5ex}
\mytag{5.9}\\
\vspace{-1.5ex}
\hskip -2em
\tilde\beta_i=&\mult_{p_i}(\tilde b)=\alpha_i-\beta_i,
&\tilde\omega_i=&\mult_{p_i}(\tilde u)=\alpha_i.\quad
\endxalignat
$$ 
Combining \mythetag{5.9} with the inequalities \mythetag{5.8}, we get
$$
\hskip -2em
\mult_{p_i}(\tilde u)=\tilde\omega_i>\mult_{p_i}(\tilde b)=\tilde\beta_i
>\mult_{p_i}(\tilde a)=\tilde\alpha_i=0. 
\mytag{5.10}
$$\par 
     In order to continue studying the equation \mythetag{4.26}, we write this
equation as
$$
\hskip -2em
\tilde a^2\,\tilde b^2\,\tilde u^2\,A_0^2\,C_6
+a^2\,b^2\,u^2\,B_0^2\,D_6-Z^4\,E=0,
\mytag{5.11}
$$
where $C_6$ and $D_6$ are given by the formulas \mythetag{4.23}, while $E$
is a new parameter:
$$
E=8\,Z^2-6\,\tilde b^2\,a^2-6\,b^2\,\tilde a^2+2\,\tilde a^2\,u^2
+2\,\tilde b^2\,u^2+2\,\tilde u^2\,a^2+2\,\tilde u^2\,b^2.\quad
\mytag{5.12}
$$
The right hand side of the formula \mythetag{5.12} is a sum of seven terms.
Applying the formulas \mythetag{5.8}, \mythetag{5.9}, and \mythetag{5.10}
and taking into account that $p_i\neq 2$, we derive
$$
\xalignat 2
&\hskip -2em
\mult_{p_i}(8\,Z^4)=2\,\alpha_i,
&&\vphantom{a}\kern -6ex
\mult_{p_i}(-6\,\tilde b^2\,a^2)\geqslant 4\,\alpha_i-2\,\beta_i,\\
&\hskip -2em
\mult_{p_i}(-6\,b^2\,\tilde a^2)\geqslant 2\,\beta_i,
&&\vphantom{a}\kern -6ex
\mult_{p_i}(2\,\tilde b^2\,u^2)=2\,\alpha_i-2\,\beta_i,
\mytag{5.13}\\
&\hskip -2em
\mult_{p_i}(2\,\tilde u^2\,a^2)=4\,\alpha_i,
&&\vphantom{a}\kern -6ex
\mult_{p_i}(2\,\tilde u^2\,b^2)=2\,\alpha_i+2\,\beta_i,\\
\endxalignat
$$ 
The only term with the minimal multiplicity in the sum \mythetag{5.12} is
the term $2\,\tilde a^2\,u^2$:
$$
\hskip -2em
\mult_{p_i}(2\,\tilde a^2\,u^2)=0.
\mytag{5.14}
$$
Applying the lemma~\mythelemma{5.1}, from \mythetag{5.13} and \mythetag{5.14} 
we derive the multiplicity of $E$:
$$
\hskip -2em
\mult_{p_i}(E)=0.
\mytag{5.15}
$$
Using \mythetag{5.15}, we can calculate the multiplicity of the last term 
in \mythetag{5.11}:
$$
\hskip -2em
\mult_{p_i}(-Z^4\,E)=4\,\alpha_i.
\mytag{5.16}
$$\par
    Assume that both multiplicities $\mu_i$ and $\eta_i$ in \mythetag{5.2} are 
nonzero. Under this assumption for the terms in the right hand side of the 
formulas \mythetag{4.23} we have 
$$
\xalignat 2
&\hskip -2em
\mult_{p_i}(A_0^2\,C_8)\geqslant 2\,\mu_i,
&&\mult_{p_i}(a^4)=4\,\alpha_i,\\
&\hskip -2em
\mult_{p_i}(b^4)=4\,\beta_i,
&&\mult_{p_i}(u^4)=0,\\
&\hskip -2em
\mult_{p_i}(4\,a^2\,u^2)=2\,\alpha_i,
&&\mult_{p_i}(4\,b^2\,u^2)=2\,\beta_i,\\
&\hskip -2em
\mult_{p_i}(-12\,b^2\,a^2)=2\,\alpha_i+2\,\beta_i,\kern -4ex
&&\mult_{p_i}(B_0^2\,D_8)\geqslant 2\,\eta_i,
\mytag{5.17}\\
&\hskip -2em
\mult_{p_i}(\tilde a^4)=0,
&&\mult_{p_i}(\tilde b^4)=4\,\alpha_i-4\,\beta_i,\\
&\hskip -2em
\mult_{p_i}(\tilde u^4)=4\,\alpha_i,
&&\mult_{p_i}(4\,\tilde a^2\,\tilde u^2)=2\,\alpha_i,\\
&\hskip -2em
\mult_{p_i}(4\,\tilde b^2\,\tilde u^2)=4\,\alpha_i-2\,\beta_i,
&&\mult_{p_i}(-12\,\tilde b^2\,\tilde a^2)=2\,\alpha_i-2\,\beta_i.
\kern -4ex
\endxalignat
$$ 
Applying the formulas \mythetag{5.17} and the lemma~\mythelemma{5.1} to 
\mythetag{4.22}, we find
$$
\xalignat 2
&\hskip -2em
\mult_{p_i}(C_6)=0,
&&\mult_{p_i}(D_6)=0.
\mytag{5.18}
\endxalignat
$$
Now we apply \mythetag{5.18} to the equation \mythetag{5.11}. As a result 
we get
$$
\hskip -2em
\aligned
&\mult_{p_i}(\tilde a^2\,\tilde b^2\,\tilde u^2\,A_0^2\,C_6)=4\,\alpha_i
-2\,\beta_i+2\,\mu_i,\\
&\mult_{p_i}(a^2\,b^2\,u^2\,B_0^2\,D_6)=2\,\alpha_i+2\,\beta_i+2\,\eta_i.
\endaligned
\mytag{5.19}
$$
Due to \mythetag{5.19} and \mythetag{5.16} the lemma~\mythelemma{5.3} applied 
to the equation \mythetag{5.11} means that at least one of the following three 
conditions should be fulfilled:
$$
\gather
\hskip -2em
4\,\alpha_i-2\,\beta_i+2\,\mu_i=4\,\alpha_i\leqslant 
2\,\alpha_i+2\,\beta_i+2\,\eta_i,
\mytag{5.20}\\
\hskip -2em
2\,\alpha_i+2\,\beta_i+2\,\eta_i=4\,\alpha_i\leqslant 
4\,\alpha_i-2\,\beta_i+2\,\mu_i,
\mytag{5.21}\\
\hskip -2em
2\,\alpha_i+2\,\beta_i+2\,\eta_i=4\,\alpha_i-2\,\beta_i+2\,\mu_i
\leqslant 4\,\alpha_i.
\mytag{5.22}
\endgather
$$
\par
     The equality in \mythetag{5.20} is easily resolvable. Resolving this 
equality, we obtain
$$
\xalignat 2
&\hskip -2em
\mu_i=\beta_i,
&&\eta_i=\alpha_i-\beta_i.
\mytag{5.23}
\endxalignat
$$
Substituting \mythetag{5.23} back into \mythetag{5.20}, we find that the 
inequality \mythetag{5.20} turns to the equality and the condition 
\mythetag{5.20} in whole appears to be fulfilled.\par 
     The equality in \mythetag{5.21} is also easily resolvable. The solution 
of this equality coincides with \mythetag{5.23}. Substituting \mythetag{5.23} 
into \mythetag{5.20}, we again find that the inequality \mythetag{5.21} turns 
to the equality and the condition \mythetag{5.21} in whole appears to be 
fulfilled. Similarly, the solution of the equality \mythetag{5.22} coincides
with \mythetag{5.23} and the condition \mythetag{5.22} in whole appears to
be fulfilled upon substituting \mythetag{5.23} into it.\par 
     {\bf The subcase $\mu_i=0$} is slightly different. The formula \mythetag{5.16}
remains unchanged, while the formulas \mythetag{5.19} in this subcase are replaced by 
the following ones: 
$$
\hskip -2em
\aligned
&\mult_{p_i}(\tilde a^2\,\tilde b^2\,\tilde u^2\,A_0^2\,C_6)=\zeta_i\geqslant 
4\,\alpha_i-2\,\beta_i,\\
&\mult_{p_i}(a^2\,b^2\,u^2\,B_0^2\,D_6)=4\,\alpha_i+2\,\beta_i.
\endaligned
\mytag{5.24}
$$
Due to \mythetag{5.24} and \mythetag{5.16} the lemma~\mythelemma{5.3} applied 
to the equation \mythetag{5.11} means that at least one of the following three 
conditions should be fulfilled:
$$
\gather
\hskip -2em
\zeta_i=4\,\alpha_i\leqslant 
4\,\alpha_i+2\,\beta_i,
\mytag{5.25}\\
\hskip -2em
4\,\alpha_i+2\,\beta_i=4\,\alpha_i\leqslant\zeta_i,
\mytag{5.26}\\
\hskip -2em
4\,\alpha_i+2\,\beta_i=\zeta_i\leqslant 4\,\alpha_i.
\mytag{5.27}
\endgather
$$
The conditions \mythetag{5.26} and \mythetag{5.27} are inconsistent since 
$\beta_i>0$. However, the subcase $\mu_i=0$ in whole is consistent because
\mythetag{5.25} is consistent. In this subcase $\eta_i=\alpha_i$ due to the
relationships \mythetag{5.2} and \mythetag{5.9}.\par
     {\bf The subcase $\eta_i=0$} is another option. In this subcase the formula 
\mythetag{5.16} remains unchanged, while the formulas \mythetag{5.19} are 
replaced by the following ones: 
$$
\hskip -2em
\aligned
&\mult_{p_i}(\tilde a^2\,\tilde b^2\,\tilde u^2\,A_0^2\,C_6)=6\,\alpha_i
-2\,\beta_i,\\
&\mult_{p_i}(a^2\,b^2\,u^2\,B_0^2\,D_6)=\xi_i\geqslant 2\,\alpha_i+2\,\beta_i.
\endaligned
\mytag{5.28}
$$
Due to \mythetag{5.28} and \mythetag{5.16} the lemma~\mythelemma{5.3} applied 
to the equation \mythetag{5.11} means that at least one of the following three 
conditions should be fulfilled:
$$
\gather
\hskip -2em
6\,\alpha_i-2\,\beta_i=4\,\alpha_i\leqslant 
\xi_i,
\mytag{5.29}\\
\hskip -2em
\xi_i=4\,\alpha_i\leqslant 
6\,\alpha_i-2\,\beta_i,
\mytag{5.30}\\
\hskip -2em
\xi_i=6\,\alpha_i-2\,\beta_i
\leqslant 4\,\alpha_i.
\mytag{5.31}
\endgather
$$
The conditions \mythetag{5.29} and \mythetag{5.31} are inconsistent since 
$\alpha_i>\beta_i$. However, the subcase $\eta_i=0$ in whole is consistent 
because \mythetag{5.30} is consistent. In this subcase $\mu_i=\alpha_i$ due 
to the relationships \mythetag{5.2} and \mythetag{5.9}.\par
     The cases and subcases are too numerous. In order to describe them we
use tables. For this purpose let's introduce the following notations:
$$
\xalignat 2
&\hskip -2em
m_C=\mult_{p_i}(\tilde a^2\,\tilde b^2\,\tilde u^2\,A_0^2\,C_6),
&&m_D=\mult_{p_i}(a^2\,b^2\,u^2\,B_0^2\,D_6).\qquad
\mytag{5.32}
\endxalignat
$$
Besides \mythetag{5.32}, let's denote through $m_E$ the multiplicity of the 
last term in \mythetag{5.11}:
$$
\hskip -2em
m_E=\mult_{p_i}(-Z^4\,E).
\mytag{5.33}
$$
In terms of \mythetag{5.32} and \mythetag{5.33} we build the table for the 
first case considered above.\par
\vskip 1ex
\noindent\mytable{5.1}\par
\newdimen\tablesize
\tablesize=\hsize
\advance\tablesize -3pt
\vbox{
\hrule width\hsize
\noindent
\vrule height 13pt depth 5pt
\vtop{\hsize=\tablesize
\hbox to\hsize{\ \ $\alpha_i>\beta_i>\omega_i=0$, $p_i\neq 2$\ \ 
\vrule height 13pt depth 5pt\ \ $\mu_i>0$ and $\eta_i>0\ \Rightarrow
\ \mu_i=\beta_i$ and $\eta_i=\alpha_i-\beta_i$
\hss}
}
\vrule height 13pt depth 5pt
\hrule width\hsize
\noindent
\vrule height 13pt depth 5pt
\vtop{\hsize=\tablesize
\hbox to\hsize{\ \ \ \ \ $m_C=4\,\alpha_i-2\,\beta_i+2\,\mu_i$, \ \ \ \ \ \ \ \ 
$m_D=2\,\alpha_i+2\,\beta_i+2\,\eta_i$, \ \ \ \ \ \ \ \ \ 
$m_E=4\,\alpha_i$
\hss}
}
\vrule height 13pt depth 5pt
\hrule width\hsize
\noindent
\vrule height 13pt depth 5pt
\vtop{\hsize=\tablesize
\hbox to\hsize{\kern 8em 
$4\,\alpha_i-2\,\beta_i+2\,\mu_i=4\,\alpha_i\leqslant 
2\,\alpha_i+2\,\beta_i+2\,\eta_i$
\hss\vrule height 13pt depth 5pt\ \,$\surd$\ }
}
\vrule height 13pt depth 5pt\vskip -1pt
\noindent
\vrule height 13pt depth 5pt
\vtop{\hsize=\tablesize
\hbox to\hsize{\kern 8em 
$2\,\alpha_i+2\,\beta_i+2\,\eta_i=4\,\alpha_i\leqslant 
4\,\alpha_i-2\,\beta_i+2\,\mu_i$
\hss\vrule height 13pt depth 5pt\ \,$\surd$\ }
}
\vrule height 13pt depth 5pt\vskip -1pt
\noindent
\vrule height 13pt depth 5pt
\vtop{\hsize=\tablesize
\hbox to\hsize{\kern 8em 
$2\,\alpha_i+2\,\beta_i+2\,\eta_i=4\,\alpha_i-2\,\beta_i+2\,\mu_i
\leqslant 4\,\alpha_i$
\hss\vrule height 13pt depth 5pt\ \,$\surd$\ }
}
\vrule height 13pt depth 5pt
\hrule width\hsize
}
\vskip 1ex
\noindent 
The first row of the table~\mythetable{5.1} is a general information. The second 
row of this table reflects the formulas \mythetag{5.19} and \mythetag{5.16}. The 
rest of this table is the conditions \mythetag{5.20}, \mythetag{5.21}, and 
\mythetag{5.22}. Check marks say that all of these conditions are consistent.
\par
\vskip 1ex
\noindent\mytable{5.2}\par
\vbox{
\hrule width\hsize
\noindent
\vrule height 13pt depth 5pt
\vtop{\hsize=\tablesize
\hbox to\hsize{\ \ $\alpha_i>\beta_i>\omega_i=0$, $p_i\neq 2$\ \ 
\vrule height 13pt depth 5pt\ \ $\mu_i=0$ and $\eta_i=\alpha_i$
\hss}
}
\vrule height 13pt depth 5pt
\hrule width\hsize
\noindent
\vrule height 13pt depth 5pt
\vtop{\hsize=\tablesize
\hbox to\hsize{\ \ \ \ \ $m_C=\zeta_i\geqslant 4\,\alpha_i-2\,\beta_i$, 
\kern 5em
$m_D=4\,\alpha_i+2\,\beta_i$,\kern 5em
$m_E=4\,\alpha_i$
\hss}
}
\vrule height 13pt depth 5pt
\hrule width\hsize
\noindent
\vrule height 13pt depth 5pt
\vtop{\hsize=\tablesize
\hbox to\hsize{\kern 12em 
$\zeta_i=4\,\alpha_i\leqslant 4\,\alpha_i+2\,\beta_i$
\hss\vrule height 13pt depth 5pt\ \,$\surd$\ }
}
\vrule height 13pt depth 5pt\vskip -1pt
\noindent
\vrule height 13pt depth 5pt
\vtop{\hsize=\tablesize
\hbox to\hsize{\kern 12em 
$4\,\alpha_i+2\,\beta_i=4\,\alpha_i\leqslant\zeta_i$
\hss\vrule height 13pt depth 5pt\ \,$\phantom{\surd}$\ }
}
\vrule height 13pt depth 5pt\vskip -1pt
\noindent
\vrule height 13pt depth 5pt
\vtop{\hsize=\tablesize
\hbox to\hsize{\kern 12em 
$4\,\alpha_i+2\,\beta_i=\zeta_i\leqslant 4\,\alpha_i$
\hss\vrule height 13pt depth 5pt\ \,$\phantom{\surd}$\ }
}
\vrule height 13pt depth 5pt
\hrule width\hsize
}
\vskip 1ex
\noindent
The last two rows in the table~\mythetable{5.2} are not check marked. 
This means that the corresponding conditions \mythetag{5.26} and
\mythetag{5.27} are not consistent.\par
\vskip 1ex
\noindent\mytable{5.3}\par
\vbox{
\hrule width\hsize
\noindent
\vrule height 13pt depth 5pt
\vtop{\hsize=\tablesize
\hbox to\hsize{\ \ $\alpha_i>\beta_i>\omega_i=0$, $p_i\neq 2$\ \ 
\vrule height 13pt depth 5pt\ \ $\mu_i=\alpha_i$ and $\eta_i=0$
\hss}
}
\vrule height 13pt depth 5pt
\hrule width\hsize
\noindent
\vrule height 13pt depth 5pt
\vtop{\hsize=\tablesize
\hbox to\hsize{\ \ \ \ \ $m_C=6\,\alpha_i-2\,\beta_i$, 
\kern 5em
$m_D=\xi_i\geqslant 2\,\alpha_i+2\,\beta_i$,\kern 5em
$m_E=4\,\alpha_i$
\hss}
}
\vrule height 13pt depth 5pt
\hrule width\hsize
\noindent
\vrule height 13pt depth 5pt
\vtop{\hsize=\tablesize
\hbox to\hsize{\kern 12em 
$6\,\alpha_i-2\,\beta_i=4\,\alpha_i\leqslant\xi_i$
\hss\vrule height 13pt depth 5pt\ \,$\phantom{\surd}$\ }
}
\vrule height 13pt depth 5pt\vskip -1pt
\noindent
\vrule height 13pt depth 5pt
\vtop{\hsize=\tablesize
\hbox to\hsize{\kern 12em 
$\xi_i=4\,\alpha_i\leqslant 6\,\alpha_i-2\,\beta_i$
\hss\vrule height 13pt depth 5pt\ \,$\surd$\ }
}
\vrule height 13pt depth 5pt\vskip -1pt
\noindent
\vrule height 13pt depth 5pt
\vtop{\hsize=\tablesize
\hbox to\hsize{\kern 12em 
$\xi_i=6\,\alpha_i-2\,\beta_i\leqslant 4\,\alpha_i$
\hss\vrule height 13pt depth 5pt\ \,$\phantom{\surd}$\ }
}
\vrule height 13pt depth 5pt
\hrule width\hsize
}
\vskip 1ex
     Below we list other cases and subcases in a tabular form without any 
detailed comments for them. 
\vskip 1ex
\noindent\mytable{5.4}\par
\vbox{
\hrule width\hsize
\noindent
\vrule height 13pt depth 5pt
\vtop{\hsize=\tablesize
\hbox to\hsize{\ \ $\alpha_i>\beta_i>\omega_i=0$, $p_i=2$\ \ 
\vrule height 13pt depth 5pt\ \ $\mu_i>0$ and $\eta_i>0\ \Rightarrow
\ \mu_i=\beta_i$ and $\eta_i=\alpha_i-\beta_i$
\hss}
}
\vrule height 13pt depth 5pt
\hrule width\hsize
\noindent
\vrule height 13pt depth 5pt
\vtop{\hsize=\tablesize
\hbox to\hsize{\ \ \ \ \ $m_C=4\,\alpha_i-2\,\beta_i+2\,\mu_i$, \ \ \ \ \ \ \ \ 
$m_D=2\,\alpha_i+2\,\beta_i+2\,\eta_i$, \ \ \ \ \ \ \ \ \ 
$m_E=4\,\alpha_i+1$
\hss}
}
\vrule height 13pt depth 5pt
\hrule width\hsize
\noindent
\vrule height 13pt depth 5pt
\vtop{\hsize=\tablesize
\hbox to\hsize{\kern 7em 
$4\,\alpha_i-2\,\beta_i+2\,\mu_i=4\,\alpha_i+1\leqslant 
2\,\alpha_i+2\,\beta_i+2\,\eta_i$
\hss\vrule height 13pt depth 5pt\ \,$\phantom{\surd}$\ }
}
\vrule height 13pt depth 5pt\vskip -1pt
\noindent
\vrule height 13pt depth 5pt
\vtop{\hsize=\tablesize
\hbox to\hsize{\kern 7em 
$2\,\alpha_i+2\,\beta_i+2\,\eta_i=4\,\alpha_i+1\leqslant 
4\,\alpha_i-2\,\beta_i+2\,\mu_i$
\hss\vrule height 13pt depth 5pt\ \,$\phantom{\surd}$\ }
}
\vrule height 13pt depth 5pt\vskip -1pt
\noindent
\vrule height 13pt depth 5pt
\vtop{\hsize=\tablesize
\hbox to\hsize{\kern 7em 
$2\,\alpha_i+2\,\beta_i+2\,\eta_i=4\,\alpha_i-2\,\beta_i+2\,\mu_i
\leqslant 4\,\alpha_i+1$
\hss\vrule height 13pt depth 5pt\ \,$\surd$\ }
}
\vrule height 13pt depth 5pt
\hrule width\hsize
}
\noindent\mytable{5.5}\par
\vbox{
\hrule width\hsize
\noindent
\vrule height 13pt depth 5pt
\vtop{\hsize=\tablesize
\hbox to\hsize{\ \ $\alpha_i>\beta_i>\omega_i=0$, $p_i=2$\ \ 
\vrule height 13pt depth 5pt\ \ $\mu_i=0$ and $\eta_i=\alpha_i$
\hss}
}
\vrule height 13pt depth 5pt
\hrule width\hsize
\noindent
\vrule height 13pt depth 5pt
\vtop{\hsize=\tablesize
\hbox to\hsize{\ \ \ \ \ $m_C=\zeta_i\geqslant 4\,\alpha_i-2\,\beta_i$, 
\kern 4.5em
$m_D=4\,\alpha_i+2\,\beta_i$,\kern 4.5em
$m_E=4\,\alpha_i+1$
\hss}
}
\vrule height 13pt depth 5pt
\hrule width\hsize
\noindent
\vrule height 13pt depth 5pt
\vtop{\hsize=\tablesize
\hbox to\hsize{\kern 12em 
$\zeta_i=4\,\alpha_i+1\leqslant 4\,\alpha_i+2\,\beta_i$
\hss\vrule height 13pt depth 5pt\ \,$\surd$\ }
}
\vrule height 13pt depth 5pt\vskip -1pt
\noindent
\vrule height 13pt depth 5pt
\vtop{\hsize=\tablesize
\hbox to\hsize{\kern 12em 
$4\,\alpha_i+2\,\beta_i=4\,\alpha_i+1\leqslant\zeta_i$
\hss\vrule height 13pt depth 5pt\ \,$\phantom{\surd}$\ }
}
\vrule height 13pt depth 5pt\vskip -1pt
\noindent
\vrule height 13pt depth 5pt
\vtop{\hsize=\tablesize
\hbox to\hsize{\kern 12em 
$4\,\alpha_i+2\,\beta_i=\zeta_i\leqslant 4\,\alpha_i+1$
\hss\vrule height 13pt depth 5pt\ \,$\phantom{\surd}$\ }
}
\vrule height 13pt depth 5pt
\hrule width\hsize
}
\vskip 1ex plus 3pt
\noindent\mytable{5.6}\par
\vbox{
\hrule width\hsize
\noindent
\vrule height 13pt depth 5pt
\vtop{\hsize=\tablesize
\hbox to\hsize{\ \ $\alpha_i>\beta_i>\omega_i=0$, $p_i=2$\ \ 
\vrule height 13pt depth 5pt\ \ $\mu_i=\alpha_i$ and $\eta_i=0$
\hss}
}
\vrule height 13pt depth 5pt
\hrule width\hsize
\noindent
\vrule height 13pt depth 5pt
\vtop{\hsize=\tablesize
\hbox to\hsize{\ \ \ \ \ $m_C=6\,\alpha_i-2\,\beta_i$, 
\kern 4.5em
$m_D=\xi_i\geqslant 2\,\alpha_i+2\,\beta_i$,\kern 4.5em
$m_E=4\,\alpha_i+1$
\hss}
}
\vrule height 13pt depth 5pt
\hrule width\hsize
\noindent
\vrule height 13pt depth 5pt
\vtop{\hsize=\tablesize
\hbox to\hsize{\kern 12em 
$6\,\alpha_i-2\,\beta_i=4\,\alpha_i+1\leqslant\xi_i$
\hss\vrule height 13pt depth 5pt\ \,$\phantom{\surd}$\ }
}
\vrule height 13pt depth 5pt\vskip -1pt
\noindent
\vrule height 13pt depth 5pt
\vtop{\hsize=\tablesize
\hbox to\hsize{\kern 12em 
$\xi_i=4\,\alpha_i+1\leqslant 6\,\alpha_i-2\,\beta_i$
\hss\vrule height 13pt depth 5pt\ \,$\surd$\ }
}
\vrule height 13pt depth 5pt\vskip -1pt
\noindent
\vrule height 13pt depth 5pt
\vtop{\hsize=\tablesize
\hbox to\hsize{\kern 12em 
$\xi_i=6\,\alpha_i-2\,\beta_i\leqslant 4\,\alpha_i+1$
\hss\vrule height 13pt depth 5pt\ \,$\phantom{\surd}$\ }
}
\vrule height 13pt depth 5pt
\hrule width\hsize
}
\vskip 1ex plus 3pt
\noindent\mytable{5.7}\par
\vbox{
\hrule width\hsize
\noindent
\vrule height 13pt depth 5pt
\vtop{\hsize=\tablesize
\hbox to\hsize{\ \ $\alpha_i>\beta_i=\omega_i=0$, $p_i\neq 2,3$\ \ 
\vrule height 13pt depth 5pt\ \ $\eta_i>0 \Rightarrow \eta_i=\alpha_i$
and $\mu_i=0$
\hss}
}
\vrule height 13pt depth 5pt
\hrule width\hsize
\noindent
\vrule height 13pt depth 5pt
\vtop{\hsize=\tablesize
\hbox to\hsize{\ \ \ \ \ $m_C=\zeta_i\geqslant 4\,\alpha_i+2\,\mu_i$, 
\kern 4em
$m_D=2\,\alpha_i+2\,\eta_i$,\kern 4em
$m_E=\varkappa_i\geqslant 4\,\alpha_i$
\hss}
}
\vrule height 13pt depth 5pt
\hrule width\hsize
\noindent
\vrule height 13pt depth 5pt
\vtop{\hsize=\tablesize
\hbox to\hsize{\kern 12.5em 
$\varkappa_i=\zeta_i\leqslant 2\,\alpha_i+2\,\eta_i$
\hss\vrule height 13pt depth 5pt\ \,$\surd$\ }
}
\vrule height 13pt depth 5pt\vskip -1pt
\noindent
\vrule height 13pt depth 5pt
\vtop{\hsize=\tablesize
\hbox to\hsize{\kern 12.5em 
$\varkappa_i=2\,\alpha_i+2\,\eta_i\leqslant\zeta_i$
\hss\vrule height 13pt depth 5pt\ \,$\surd$\ }
}
\vrule height 13pt depth 5pt\vskip -1pt
\noindent
\vrule height 13pt depth 5pt
\vtop{\hsize=\tablesize
\hbox to\hsize{\kern 12.5em 
$\zeta_i=2\,\alpha_i+2\,\eta_i\leqslant\varkappa_i$
\hss\vrule height 13pt depth 5pt\ \,$\surd$\ }
}
\vrule height 13pt depth 5pt
\hrule width\hsize
}
\vskip 1ex plus 3pt
\noindent\mytable{5.8}\par
\vbox{
\hrule width\hsize
\noindent
\vrule height 13pt depth 5pt
\vtop{\hsize=\tablesize
\hbox to\hsize{\ \ $\alpha_i>\beta_i=\omega_i=0$, $p_i\neq 2,3$\ \ 
\vrule height 13pt depth 5pt\ \ $\eta_i=0$ and $\mu_i=\alpha_i$
\hss}
}
\vrule height 13pt depth 5pt
\hrule width\hsize
\noindent
\vrule height 13pt depth 5pt
\vtop{\hsize=\tablesize
\hbox to\hsize{\ \ \ \ \ $m_C=\zeta_i\geqslant 6\,\alpha_i$, 
\kern 5.7em
$m_D=\xi_i\geqslant 2\,\alpha_i$,\kern 5.7em
$m_E=\varkappa_i\geqslant 4\,\alpha_i$
\hss}
}
\vrule height 13pt depth 5pt
\hrule width\hsize
\noindent
\vrule height 13pt depth 5pt
\vtop{\hsize=\tablesize
\hbox to\hsize{\kern 15em 
$\varkappa_i=\zeta_i\leqslant\xi_i$
\hss\vrule height 13pt depth 5pt\ \,$\surd$\ }
}
\vrule height 13pt depth 5pt\vskip -1pt
\noindent
\vrule height 13pt depth 5pt
\vtop{\hsize=\tablesize
\hbox to\hsize{\kern 15em 
$\varkappa_i=\xi_i\leqslant\zeta_i$
\hss\vrule height 13pt depth 5pt\ \,$\surd$\ }
}
\vrule height 13pt depth 5pt\vskip -1pt
\noindent
\vrule height 13pt depth 5pt
\vtop{\hsize=\tablesize
\hbox to\hsize{\kern 15em 
$\zeta_i=\xi_i\leqslant\varkappa_i$
\hss\vrule height 13pt depth 5pt\ \,$\surd$\ }
}
\vrule height 13pt depth 5pt
\hrule width\hsize
}
\vskip 1ex
     The following tables correspond to the special values of the prime
factor $p_i$, i\.\,e\. to $p_i=2$ and to $p_i=3$.\par 
\vskip 1ex plus 3pt
\noindent\mytable{5.9}\par
\vbox{
\hrule width\hsize
\noindent
\vrule height 13pt depth 5pt
\vtop{\hsize=\tablesize
\hbox to\hsize{\ \ $\alpha_i>\beta_i=\omega_i=0$, $p_i=2$\ \ 
\vrule height 13pt depth 5pt\ \ $\eta_i>0 \Rightarrow \eta_i=\alpha_i$
and $\mu_i=0$
\hss}
}
\vrule height 13pt depth 5pt
\hrule width\hsize
\noindent
\vrule height 13pt depth 5pt
\vtop{\hsize=\tablesize
\hbox to\hsize{\ \ \ \ \ $m_C=\zeta_i\geqslant 4\,\alpha_i+2\,\mu_i$, 
\kern 3em
$m_D=2\,\alpha_i+2\,\eta_i$,\kern 3em
$m_E=\varkappa_i\geqslant 4\,\alpha_i+1$
\hss}
}
\vrule height 13pt depth 5pt
\hrule width\hsize
\noindent
\vrule height 13pt depth 5pt
\vtop{\hsize=\tablesize
\hbox to\hsize{\kern 13em 
$\varkappa_i=\zeta_i\leqslant 2\,\alpha_i+2\,\eta_i$
\hss\vrule height 13pt depth 5pt\ \,$\phantom{\surd}$\ }
}
\vrule height 13pt depth 5pt\vskip -1pt
\noindent
\vrule height 13pt depth 5pt
\vtop{\hsize=\tablesize
\hbox to\hsize{\kern 13em 
$\varkappa_i=2\,\alpha_i+2\,\eta_i\leqslant\zeta_i$
\hss\vrule height 13pt depth 5pt\ \,$\phantom{\surd}$\ }
}
\vrule height 13pt depth 5pt\vskip -1pt
\noindent
\vrule height 13pt depth 5pt
\vtop{\hsize=\tablesize
\hbox to\hsize{\kern 13em 
$\zeta_i=2\,\alpha_i+2\,\eta_i\leqslant\varkappa_i$
\hss\vrule height 13pt depth 5pt\ \,$\surd$\ }
}
\vrule height 13pt depth 5pt
\hrule width\hsize
}\pagebreak
\noindent\mytable{5.10}\par
\vbox{
\hrule width\hsize
\noindent
\vrule height 13pt depth 5pt
\vtop{\hsize=\tablesize
\hbox to\hsize{\ \ $\alpha_i>\beta_i=\omega_i=0$, $p_i=2$\ \ 
\vrule height 13pt depth 5pt\ \ $\eta_i=0$ and $\mu_i=\alpha_i$
\hss}
}
\vrule height 13pt depth 5pt
\hrule width\hsize
\noindent
\vrule height 13pt depth 5pt
\vtop{\hsize=\tablesize
\hbox to\hsize{\ \ \ \ \ $m_C=\zeta_i\geqslant 6\,\alpha_i$, 
\kern 5em
$m_D=\xi_i\geqslant 2\,\alpha_i$,\kern 5em
$m_E=\varkappa_i\geqslant 4\,\alpha_i+1$
\hss}
}
\vrule height 13pt depth 5pt
\hrule width\hsize
\noindent
\vrule height 13pt depth 5pt
\vtop{\hsize=\tablesize
\hbox to\hsize{\kern 15em 
$\varkappa_i=\zeta_i\leqslant\xi_i$
\hss\vrule height 13pt depth 5pt\ \,$\surd$\ }
}
\vrule height 13pt depth 5pt\vskip -1pt
\noindent
\vrule height 13pt depth 5pt
\vtop{\hsize=\tablesize
\hbox to\hsize{\kern 15em 
$\varkappa_i=\xi_i\leqslant\zeta_i$
\hss\vrule height 13pt depth 5pt\ \,$\surd$\ }
}
\vrule height 13pt depth 5pt\vskip -1pt
\noindent
\vrule height 13pt depth 5pt
\vtop{\hsize=\tablesize
\hbox to\hsize{\kern 15em 
$\zeta_i=\xi_i\leqslant\varkappa_i$
\hss\vrule height 13pt depth 5pt\ \,$\surd$\ }
}
\vrule height 13pt depth 5pt
\hrule width\hsize
}
\vskip 1ex plus 3pt
\noindent\mytable{5.11}\par
\vbox{
\hrule width\hsize
\noindent
\vrule height 13pt depth 5pt
\vtop{\hsize=\tablesize
\hbox to\hsize{\ \ $\alpha_i>\beta_i=\omega_i=0$, $p_i=3$\ \ 
\vrule height 13pt depth 5pt\ \ $\eta_i>0 \Rightarrow \eta_i=\alpha_i$
and $\mu_i=0$
\hss}
}
\vrule height 13pt depth 5pt
\hrule width\hsize
\noindent
\vrule height 13pt depth 5pt
\vtop{\hsize=\tablesize
\hbox to\hsize{\ \ \ \ \ $m_C=\zeta_i\geqslant 4\,\alpha_i+2\,\mu_i$, 
\kern 5.3em
$m_D=2\,\alpha_i+2\,\eta_i$,\kern 5.3em
$m_E=4\,\alpha_i$
\hss}
}
\vrule height 13pt depth 5pt
\hrule width\hsize
\noindent
\vrule height 13pt depth 5pt
\vtop{\hsize=\tablesize
\hbox to\hsize{\kern 13em 
$4\,\alpha_i=\zeta_i\leqslant 2\,\alpha_i+2\,\eta_i$
\hss\vrule height 13pt depth 5pt\ \,$\surd$\ }
}
\vrule height 13pt depth 5pt\vskip -1pt
\noindent
\vrule height 13pt depth 5pt
\vtop{\hsize=\tablesize
\hbox to\hsize{\kern 13em 
$4\,\alpha_i=2\,\alpha_i+2\,\eta_i\leqslant\zeta_i$
\hss\vrule height 13pt depth 5pt\ \,$\surd$\ }
}
\vrule height 13pt depth 5pt\vskip -1pt
\noindent
\vrule height 13pt depth 5pt
\vtop{\hsize=\tablesize
\hbox to\hsize{\kern 13em 
$\zeta_i=2\,\alpha_i+2\,\eta_i\leqslant4\,\alpha_i$
\hss\vrule height 13pt depth 5pt\ \,$\surd$\ }
}
\vrule height 13pt depth 5pt
\hrule width\hsize
}
\vskip 1ex plus 3pt
\noindent\mytable{5.12}\par
\vbox{
\hrule width\hsize
\noindent
\vrule height 13pt depth 5pt
\vtop{\hsize=\tablesize
\hbox to\hsize{\ \ $\alpha_i>\beta_i=\omega_i=0$, $p_i=3$\ \ 
\vrule height 13pt depth 5pt\ \ $\eta_i=0$ and $\mu_i=\alpha_i$
\hss}
}
\vrule height 13pt depth 5pt
\hrule width\hsize
\noindent
\vrule height 13pt depth 5pt
\vtop{\hsize=\tablesize
\hbox to\hsize{\ \ \ \ \ $m_C=\zeta_i\geqslant 6\,\alpha_i$, 
\kern 7em
$m_D=\xi_i\geqslant 2\,\alpha_i$,\kern 7em
$m_E=4\,\alpha_i$
\hss}
}
\vrule height 13pt depth 5pt
\hrule width\hsize
\noindent
\vrule height 13pt depth 5pt
\vtop{\hsize=\tablesize
\hbox to\hsize{\kern 15em 
$4\,\alpha_i=\zeta_i\leqslant\xi_i$
\hss\vrule height 13pt depth 5pt\ \,$\phantom{\surd}$\ }
}
\vrule height 13pt depth 5pt\vskip -1pt
\noindent
\vrule height 13pt depth 5pt
\vtop{\hsize=\tablesize
\hbox to\hsize{\kern 15em 
$4\,\alpha_i=\xi_i\leqslant\zeta_i$
\hss\vrule height 13pt depth 5pt\ \,$\surd$\ }
}
\vrule height 13pt depth 5pt\vskip -1pt
\noindent
\vrule height 13pt depth 5pt
\vtop{\hsize=\tablesize
\hbox to\hsize{\kern 15em 
$\zeta_i=\xi_i\leqslant4\,\alpha_i$
\hss\vrule height 13pt depth 5pt\ \,$\phantom{\surd}$\ }
}
\vrule height 13pt depth 5pt
\hrule width\hsize
}
\vskip 1ex plus 3pt
    In the following cases the multiplicity $\beta_i$ is not zero. It is 
equal to the multiplicity $\alpha_i$ and, according to \mythetag{3.3}, we
have $\theta_i=\max(\alpha_i,\beta_i,\omega_i)=\alpha_i$. 
\vskip 1ex plus 3pt
\noindent\mytable{5.13}\par
\vbox{
\hrule width\hsize
\noindent
\vrule height 13pt depth 5pt
\vtop{\hsize=\tablesize
\hbox to\hsize{\ \ $\alpha_i=\beta_i>\omega_i=0$, $p_i\neq 2$\ \ 
\vrule height 13pt depth 5pt\ \ $\mu_i>0 \ \Rightarrow
\ \mu_i=\alpha_i$ and $\eta_i=0$
\hss}
}
\vrule height 13pt depth 5pt
\hrule width\hsize
\noindent
\vrule height 13pt depth 5pt
\vtop{\hsize=\tablesize
\hbox to\hsize{\ \ \ \ \ $m_C=2\,\alpha_i+2\,\mu_i$,\kern 4.2em
$m_D=\xi_i\geqslant 4\,\alpha_i+2\,\eta_i$,\kern 4.2em 
$m_E=\varkappa_i\geqslant 4\,\alpha_i$
\hss}
}
\vrule height 13pt depth 5pt
\hrule width\hsize
\noindent
\vrule height 13pt depth 5pt
\vtop{\hsize=\tablesize
\hbox to\hsize{\kern 13.5em 
$2\,\alpha_i+2\,\mu_i=\varkappa_i\leqslant\xi_i$
\hss\vrule height 13pt depth 5pt\ \,$\surd$\ }
}
\vrule height 13pt depth 5pt\vskip -1pt
\noindent
\vrule height 13pt depth 5pt
\vtop{\hsize=\tablesize
\hbox to\hsize{\kern 13.5em 
$\xi_i=\varkappa_i\leqslant 2\,\alpha_i+2\,\mu_i$
\hss\vrule height 13pt depth 5pt\ \,$\surd$\ }
}
\vrule height 13pt depth 5pt\vskip -1pt
\noindent
\vrule height 13pt depth 5pt
\vtop{\hsize=\tablesize
\hbox to\hsize{\kern 13.5em 
$\xi_i=2\,\alpha_i+2\,\mu_i\leqslant\varkappa_i$
\hss\vrule height 13pt depth 5pt\ \,$\surd$\ }
}
\vrule height 13pt depth 5pt
\hrule width\hsize
}
\vskip 1ex plus 3pt
\noindent\mytable{5.14}\par
\vbox{
\hrule width\hsize
\noindent
\vrule height 13pt depth 5pt
\vtop{\hsize=\tablesize
\hbox to\hsize{\ \ $\alpha_i=\beta_i>\omega_i=0$, $p_i\neq 2$\ \ 
\vrule height 13pt depth 5pt\ \ $\mu_i=0$ and $\eta_i=\alpha_i$
\hss}
}
\vrule height 13pt depth 5pt
\hrule width\hsize
\noindent
\vrule height 13pt depth 5pt
\vtop{\hsize=\tablesize
\hbox to\hsize{\ \ \ \ \ $m_C=\zeta_i\geqslant 2\,\alpha_i$,\kern 6em
$m_D=\xi_i\geqslant 6\,\alpha_i$,\kern 6em 
$m_E=\varkappa_i\geqslant 4\,\alpha_i$
\hss}
}
\vrule height 13pt depth 5pt
\hrule width\hsize
\noindent
\vrule height 13pt depth 5pt
\vtop{\hsize=\tablesize
\hbox to\hsize{\kern 15em 
$\zeta_i=\varkappa_i\leqslant\xi_i$
\hss\vrule height 13pt depth 5pt\ \,$\surd$\ }
}
\vrule height 13pt depth 5pt\vskip -1pt
\noindent
\vrule height 13pt depth 5pt
\vtop{\hsize=\tablesize
\hbox to\hsize{\kern 15em 
$\xi_i=\varkappa_i\leqslant\zeta_i$
\hss\vrule height 13pt depth 5pt\ \,$\surd$\ }
}
\vrule height 13pt depth 5pt\vskip -1pt
\noindent
\vrule height 13pt depth 5pt
\vtop{\hsize=\tablesize
\hbox to\hsize{\kern 15em 
$\xi_i=\zeta_i\leqslant\varkappa_i$
\hss\vrule height 13pt depth 5pt\ \,$\surd$\ }
}
\vrule height 13pt depth 5pt
\hrule width\hsize
}\pagebreak
     The case $p_i=3$ for $\alpha_i=\beta_i>\omega_i=0$ leads to the same 
equalities and inequalities as the other cases in the tables \mythetable{5.13} 
and \mythetable{5.14}.
\vskip 1ex plus 3pt
\noindent\mytable{5.15}\par
\vbox{
\hrule width\hsize
\noindent
\vrule height 13pt depth 5pt
\vtop{\hsize=\tablesize
\hbox to\hsize{\ \ $\alpha_i=\beta_i>\omega_i=0$, $p_i=2$\ \ 
\vrule height 13pt depth 5pt\ \ $\mu_i>0 \ \Rightarrow
\ \mu_i=\alpha_i$ and $\eta_i=0$
\hss}
}
\vrule height 13pt depth 5pt
\hrule width\hsize
\noindent
\vrule height 13pt depth 5pt
\vtop{\hsize=\tablesize
\hbox to\hsize{\ \ \ \ \ $m_C=2\,\alpha_i+2\,\mu_i$,\kern 3.5em
$m_D=\xi_i\geqslant 4\,\alpha_i+2\,\eta_i$,\kern 3.5em 
$m_E=\varkappa_i\geqslant 4\,\alpha_i+1$
\hss}
}
\vrule height 13pt depth 5pt
\hrule width\hsize
\noindent
\vrule height 13pt depth 5pt
\vtop{\hsize=\tablesize
\hbox to\hsize{\kern 13.5em 
$2\,\alpha_i+2\,\mu_i=\varkappa_i\leqslant\xi_i$
\hss\vrule height 13pt depth 5pt\ \,$\phantom{\surd}$\ }
}
\vrule height 13pt depth 5pt\vskip -1pt
\noindent
\vrule height 13pt depth 5pt
\vtop{\hsize=\tablesize
\hbox to\hsize{\kern 13.5em 
$\xi_i=\varkappa_i\leqslant 2\,\alpha_i+2\,\mu_i$
\hss\vrule height 13pt depth 5pt\ \,$\phantom{\surd}$\ }
}
\vrule height 13pt depth 5pt\vskip -1pt
\noindent
\vrule height 13pt depth 5pt
\vtop{\hsize=\tablesize
\hbox to\hsize{\kern 13.5em 
$\xi_i=2\,\alpha_i+2\,\mu_i\leqslant\varkappa_i$
\hss\vrule height 13pt depth 5pt\ \,$\surd$\ }
}
\vrule height 13pt depth 5pt
\hrule width\hsize
}
\vskip 1ex plus 3pt
\noindent\mytable{5.16}\par
\vbox{
\hrule width\hsize
\noindent
\vrule height 13pt depth 5pt
\vtop{\hsize=\tablesize
\hbox to\hsize{\ \ $\alpha_i=\beta_i>\omega_i=0$, $p_i=2$\ \ 
\vrule height 13pt depth 5pt\ \ $\mu_i=0$ and $\eta_i=\alpha_i$
\hss}
}
\vrule height 13pt depth 5pt
\hrule width\hsize
\noindent
\vrule height 13pt depth 5pt
\vtop{\hsize=\tablesize
\hbox to\hsize{\ \ \ \ \ $m_C=\zeta_i\geqslant 2\,\alpha_i$,\kern 5.2em
$m_D=\xi_i\geqslant 6\,\alpha_i$,\kern 5.2em
$m_E=\varkappa_i\geqslant 4\,\alpha_i+1$
\hss}
}
\vrule height 13pt depth 5pt
\hrule width\hsize
\noindent
\vrule height 13pt depth 5pt
\vtop{\hsize=\tablesize
\hbox to\hsize{\kern 15em 
$\zeta_i=\varkappa_i\leqslant\xi_i$
\hss\vrule height 13pt depth 5pt\ \,$\surd$\ }
}
\vrule height 13pt depth 5pt\vskip -1pt
\noindent
\vrule height 13pt depth 5pt
\vtop{\hsize=\tablesize
\hbox to\hsize{\kern 15em 
$\xi_i=\varkappa_i\leqslant\zeta_i$
\hss\vrule height 13pt depth 5pt\ \,$\surd$\ }
}
\vrule height 13pt depth 5pt\vskip -1pt
\noindent
\vrule height 13pt depth 5pt
\vtop{\hsize=\tablesize
\hbox to\hsize{\kern 15em 
$\xi_i=\zeta_i\leqslant\varkappa_i$
\hss\vrule height 13pt depth 5pt\ \,$\surd$\ }
}
\vrule height 13pt depth 5pt
\hrule width\hsize
}
\vskip 1ex plus 3pt
    In the following cases the multiplicity $\beta_i$ is zero, while 
$\omega_i$ is nonzero. 
\vskip 1ex plus 3pt
\noindent\mytable{5.17}\par
\vbox{
\hrule width\hsize
\noindent
\vrule height 13pt depth 5pt
\vtop{\hsize=\tablesize
\hbox to\hsize{\ \ $\alpha_i>\omega_i>\beta_i=0$, $p_i\neq 2,3$\ \ 
\vrule height 13pt depth 5pt\ \ $\mu_i>0$ and $\eta_i>0\ \Rightarrow
\ \mu_i=\omega_i$ and $\eta_i=\alpha_i-\omega_i$
\hss}
}
\vrule height 13pt depth 5pt
\hrule width\hsize
\noindent
\vrule height 13pt depth 5pt
\vtop{\hsize=\tablesize
\hbox to\hsize{\ \ \ \ \ $m_C=4\,\alpha_i-2\,\omega_i+2\,\mu_i$,\kern 3.6em
$m_D=2\,\alpha_i+2\,\omega_i+2\,\eta_i$,\kern 3.6em 
$m_E=4\,\alpha_i$
\hss}
}
\vrule height 13pt depth 5pt
\hrule width\hsize
\noindent
\vrule height 13pt depth 5pt
\vtop{\hsize=\tablesize
\hbox to\hsize{\kern 8em 
$4\,\alpha_i-2\,\omega_i+2\,\mu_i=4\,\alpha_i\leqslant 2\,\alpha_i
+2\,\omega_i+2\,\eta_i$
\hss\vrule height 13pt depth 5pt\ \,$\surd$\ }
}
\vrule height 13pt depth 5pt\vskip -1pt
\noindent
\vrule height 13pt depth 5pt
\vtop{\hsize=\tablesize
\hbox to\hsize{\kern 8em 
$2\,\alpha_i+2\,\omega_i+2\,\eta_i=4\,\alpha_i\leqslant 4\,\alpha_i
-2\,\omega_i+2\,\mu_i$
\hss\vrule height 13pt depth 5pt\ \,$\surd$\ }
}
\vrule height 13pt depth 5pt\vskip -1pt
\noindent
\vrule height 13pt depth 5pt
\vtop{\hsize=\tablesize
\hbox to\hsize{\kern 8em 
$2\,\alpha_i+2\,\omega_i+2\,\eta_i=4\,\alpha_i-2\,\omega_i+2\,\mu_i
\leqslant 4\,\alpha_i$
\hss\vrule height 13pt depth 5pt\ \,$\surd$\ }
}
\vrule height 13pt depth 5pt
\hrule width\hsize
}
\vskip 1ex plus 3pt
\noindent\mytable{5.18}\par
\vbox{
\hrule width\hsize
\noindent
\vrule height 13pt depth 5pt
\vtop{\hsize=\tablesize
\hbox to\hsize{\ \ $\alpha_i>\omega_i>\beta_i=0$, $p_i\neq 2,3$\ \ 
\vrule height 13pt depth 5pt\ \ $\mu_i=0$ and $\eta_i=\alpha_i$
\hss}
}
\vrule height 13pt depth 5pt
\hrule width\hsize
\noindent
\vrule height 13pt depth 5pt
\vtop{\hsize=\tablesize
\hbox to\hsize{\ \ \ \ \ $m_C=\zeta_i\geqslant 4\,\alpha_i-2\,\omega_i$,
\kern 5.25em
$m_D=4\,\alpha_i+2\,\omega_i$,\kern 5.25em
$m_E=4\,\alpha_i$
\hss}
}
\vrule height 13pt depth 5pt
\hrule width\hsize
\noindent
\vrule height 13pt depth 5pt
\vtop{\hsize=\tablesize
\hbox to\hsize{\kern 13em 
$\zeta_i=4\,\alpha_i\leqslant 4\,\alpha_i+2\,\omega_i$
\hss\vrule height 13pt depth 5pt\ \,$\surd$\ }
}
\vrule height 13pt depth 5pt\vskip -1pt
\noindent
\vrule height 13pt depth 5pt
\vtop{\hsize=\tablesize
\hbox to\hsize{\kern 13em 
$4\,\alpha_i+2\,\omega_i=4\,\alpha_i\leqslant\zeta_i$
\hss\vrule height 13pt depth 5pt\ \,$\phantom{\surd}$\ }
}
\vrule height 13pt depth 5pt\vskip -1pt
\noindent
\vrule height 13pt depth 5pt
\vtop{\hsize=\tablesize
\hbox to\hsize{\kern 13em 
$4\,\alpha_i+2\,\omega_i=\zeta_i\leqslant 4\,\alpha_i$
\hss\vrule height 13pt depth 5pt\ \,$\phantom{\surd}$\ }
}
\vrule height 13pt depth 5pt
\hrule width\hsize
}
\vskip 1ex plus 3pt
\noindent\mytable{5.19}\par
\vbox{
\hrule width\hsize
\noindent
\vrule height 13pt depth 5pt
\vtop{\hsize=\tablesize
\hbox to\hsize{\ \ $\alpha_i>\omega_i>\beta_i=0$, $p_i\neq 2,3$\ \ 
\vrule height 13pt depth 5pt\ \ $\mu_i=\alpha_i$ and $\eta_i=0$
\hss}
}
\vrule height 13pt depth 5pt
\hrule width\hsize
\noindent
\vrule height 13pt depth 5pt
\vtop{\hsize=\tablesize
\hbox to\hsize{\ \ \ \ \ $m_C=6\,\alpha_i-2\,\omega_i$,
\kern 5.25em
$m_D=\xi_i\geqslant 2\,\alpha_i+2\,\omega_i$,\kern 5.25em
$m_E=4\,\alpha_i$
\hss}
}
\vrule height 13pt depth 5pt
\hrule width\hsize
\noindent
\vrule height 13pt depth 5pt
\vtop{\hsize=\tablesize
\hbox to\hsize{\kern 13em 
$6\,\alpha_i-2\,\omega_i=4\,\alpha_i\leqslant\xi_i$
\hss\vrule height 13pt depth 5pt\ \,$\phantom{\surd}$\ }
}
\vrule height 13pt depth 5pt\vskip -1pt
\noindent
\vrule height 13pt depth 5pt
\vtop{\hsize=\tablesize
\hbox to\hsize{\kern 13em 
$\xi_i=4\,\alpha_i\leqslant 6\,\alpha_i-2\,\omega_i$
\hss\vrule height 13pt depth 5pt\ \,$\surd$\ }
}
\vrule height 13pt depth 5pt\vskip -1pt
\noindent
\vrule height 13pt depth 5pt
\vtop{\hsize=\tablesize
\hbox to\hsize{\kern 13em 
$\xi_i=6\,\alpha_i-2\,\omega_i\leqslant 4\,\alpha_i$
\hss\vrule height 13pt depth 5pt\ \,$\phantom{\surd}$\ }
}
\vrule height 13pt depth 5pt
\hrule width\hsize
}
\noindent\mytable{5.20}\par
\vbox{
\hrule width\hsize
\noindent
\vrule height 13pt depth 5pt
\vtop{\hsize=\tablesize
\hbox to\hsize{\ \ $\alpha_i>\omega_i>\beta_i=0$, $p_i=2,3$\ \ 
\vrule height 13pt depth 5pt\ \ $\mu_i>0$ and $\eta_i>0\ \Rightarrow
\ \mu_i=\omega_i$ and $\eta_i=\alpha_i-\omega_i$
\hss}
}
\vrule height 13pt depth 5pt
\hrule width\hsize
\noindent
\vrule height 13pt depth 5pt
\vtop{\hsize=\tablesize
\hbox to\hsize{\ \ \ \ \ $m_C=4\,\alpha_i-2\,\omega_i+2\,\mu_i$,\kern 2.9em
$m_D=2\,\alpha_i+2\,\omega_i+2\,\eta_i$,\kern 2.9em 
$m_E=4\,\alpha_i+1$
\hss}
}
\vrule height 13pt depth 5pt
\hrule width\hsize
\noindent
\vrule height 13pt depth 5pt
\vtop{\hsize=\tablesize
\hbox to\hsize{\kern 8em 
$4\,\alpha_i-2\,\omega_i+2\,\mu_i=4\,\alpha_i+1\leqslant 2\,\alpha_i
+2\,\omega_i+2\,\eta_i$
\hss\vrule height 13pt depth 5pt\ \,$\phantom{\surd}$\ }
}
\vrule height 13pt depth 5pt\vskip -1pt
\noindent
\vrule height 13pt depth 5pt
\vtop{\hsize=\tablesize
\hbox to\hsize{\kern 8em 
$2\,\alpha_i+2\,\omega_i+2\,\eta_i=4\,\alpha_i+1\leqslant 4\,\alpha_i
-2\,\omega_i+2\,\mu_i$
\hss\vrule height 13pt depth 5pt\ \,$\phantom{\surd}$\ }
}
\vrule height 13pt depth 5pt\vskip -1pt
\noindent
\vrule height 13pt depth 5pt
\vtop{\hsize=\tablesize
\hbox to\hsize{\kern 8em 
$2\,\alpha_i+2\,\omega_i+2\,\eta_i=4\,\alpha_i-2\,\omega_i+2\,\mu_i
\leqslant 4\,\alpha_i+1$
\hss\vrule height 13pt depth 5pt\ \,$\surd$\ }
}
\vrule height 13pt depth 5pt
\hrule width\hsize
}
\vskip 1ex plus 3pt
\noindent\mytable{5.21}\par
\vbox{
\hrule width\hsize
\noindent
\vrule height 13pt depth 5pt
\vtop{\hsize=\tablesize
\hbox to\hsize{\ \ $\alpha_i>\omega_i>\beta_i=0$, $p_i=2,3$\ \ 
\vrule height 13pt depth 5pt\ \ $\mu_i=0$ and $\eta_i=\alpha_i$
\hss}
}
\vrule height 13pt depth 5pt
\hrule width\hsize
\noindent
\vrule height 13pt depth 5pt
\vtop{\hsize=\tablesize
\hbox to\hsize{\ \ \ \ \ $m_C=\zeta_i\geqslant 4\,\alpha_i-2\,\omega_i$,
\kern 4.4em
$m_D=4\,\alpha_i+2\,\omega_i$,\kern 4.4em
$m_E=4\,\alpha_i+1$
\hss}
}
\vrule height 13pt depth 5pt
\hrule width\hsize
\noindent
\vrule height 13pt depth 5pt
\vtop{\hsize=\tablesize
\hbox to\hsize{\kern 12em 
$\zeta_i=4\,\alpha_i+1\leqslant 4\,\alpha_i+2\,\omega_i$
\hss\vrule height 13pt depth 5pt\ \,$\surd$\ }
}
\vrule height 13pt depth 5pt\vskip -1pt
\noindent
\vrule height 13pt depth 5pt
\vtop{\hsize=\tablesize
\hbox to\hsize{\kern 12em 
$4\,\alpha_i+2\,\omega_i=4\,\alpha_i+1\leqslant\zeta_i$
\hss\vrule height 13pt depth 5pt\ \,$\phantom{\surd}$\ }
}
\vrule height 13pt depth 5pt\vskip -1pt
\noindent
\vrule height 13pt depth 5pt
\vtop{\hsize=\tablesize
\hbox to\hsize{\kern 12em 
$4\,\alpha_i+2\,\omega_i=\zeta_i\leqslant 4\,\alpha_i+1$
\hss\vrule height 13pt depth 5pt\ \,$\phantom{\surd}$\ }
}
\vrule height 13pt depth 5pt
\hrule width\hsize
}
\vskip 1ex plus 3pt
\noindent\mytable{5.22}\par
\vbox{
\hrule width\hsize
\noindent
\vrule height 13pt depth 5pt
\vtop{\hsize=\tablesize
\hbox to\hsize{\ \ $\alpha_i>\omega_i>\beta_i=0$, $p_i=2,3$\ \ 
\vrule height 13pt depth 5pt\ \ $\mu_i=\alpha_i$ and $\eta_i=0$
\hss}
}
\vrule height 13pt depth 5pt
\hrule width\hsize
\noindent
\vrule height 13pt depth 5pt
\vtop{\hsize=\tablesize
\hbox to\hsize{\ \ \ \ \ $m_C=6\,\alpha_i-2\,\omega_i$,
\kern 4.5em
$m_D=\xi_i\geqslant 2\,\alpha_i+2\,\omega_i$,\kern 4.5em
$m_E=4\,\alpha_i+1$
\hss}
}
\vrule height 13pt depth 5pt
\hrule width\hsize
\noindent
\vrule height 13pt depth 5pt
\vtop{\hsize=\tablesize
\hbox to\hsize{\kern 12em 
$6\,\alpha_i-2\,\omega_i=4\,\alpha_i+1\leqslant\xi_i$
\hss\vrule height 13pt depth 5pt\ \,$\phantom{\surd}$\ }
}
\vrule height 13pt depth 5pt\vskip -1pt
\noindent
\vrule height 13pt depth 5pt
\vtop{\hsize=\tablesize
\hbox to\hsize{\kern 12em 
$\xi_i=4\,\alpha_i+1\leqslant 6\,\alpha_i-2\,\omega_i$
\hss\vrule height 13pt depth 5pt\ \,$\surd$\ }
}
\vrule height 13pt depth 5pt\vskip -1pt
\noindent
\vrule height 13pt depth 5pt
\vtop{\hsize=\tablesize
\hbox to\hsize{\kern 12em 
$\xi_i=6\,\alpha_i-2\,\omega_i\leqslant 4\,\alpha_i+1$
\hss\vrule height 13pt depth 5pt\ \,$\phantom{\surd}$\ }
}
\vrule height 13pt depth 5pt
\hrule width\hsize
}
\vskip 1ex plus 3pt
     In the following cases the multiplicity $\omega_i$ coincides
with the multiplicity $\alpha_i$ and, according to \mythetag{3.3}, 
we have $\theta_i=\max(\alpha_i,\beta_i,\omega_i)=\alpha_i$. 
\vskip 1ex plus 3pt
\noindent\mytable{5.23}\par
\vbox{
\hrule width\hsize
\noindent
\vrule height 13pt depth 5pt
\vtop{\hsize=\tablesize
\hbox to\hsize{\ \ $\alpha_i=\omega_i>\beta_i=0$, $p_i\neq 2,3$\ \ 
\vrule height 13pt depth 5pt\ \ $\mu_i>0\ \Rightarrow
\ \mu_i=\alpha_i$ and $\eta_i=0$
\hss}
}
\vrule height 13pt depth 5pt
\hrule width\hsize
\noindent
\vrule height 13pt depth 5pt
\vtop{\hsize=\tablesize
\hbox to\hsize{\ \ \ \ \ $m_C=2\,\alpha_i+2\,\mu_i$,\kern 4.3em
$m_D=\xi_i\geqslant 4\,\alpha_i+2\,\eta_i$,\kern 4.3em 
$m_E=\varkappa_i\geqslant 4\,\alpha_i$
\hss}
}
\vrule height 13pt depth 5pt
\hrule width\hsize
\noindent
\vrule height 13pt depth 5pt
\vtop{\hsize=\tablesize
\hbox to\hsize{\kern 13em 
$2\,\alpha_i+2\,\mu_i=\varkappa_i\leqslant\xi_i$
\hss\vrule height 13pt depth 5pt\ \,$\surd$\ }
}
\vrule height 13pt depth 5pt\vskip -1pt
\noindent
\vrule height 13pt depth 5pt
\vtop{\hsize=\tablesize
\hbox to\hsize{\kern 13em 
$\xi_i=\varkappa_i\leqslant 2\,\alpha_i+2\,\mu_i$
\hss\vrule height 13pt depth 5pt\ \,$\surd$\ }
}
\vrule height 13pt depth 5pt\vskip -1pt
\noindent
\vrule height 13pt depth 5pt
\vtop{\hsize=\tablesize
\hbox to\hsize{\kern 13em 
$\xi_i=2\,\alpha_i+2\,\mu_i\leqslant\varkappa_i$
\hss\vrule height 13pt depth 5pt\ \,$\surd$\ }
}
\vrule height 13pt depth 5pt
\hrule width\hsize
}
\vskip 1ex plus 3pt
\noindent\mytable{5.24}\par
\vbox{
\hrule width\hsize
\noindent
\vrule height 13pt depth 5pt
\vtop{\hsize=\tablesize
\hbox to\hsize{\ \ $\alpha_i=\omega_i>\beta_i=0$, $p_i\neq 2,3$\ \ 
\vrule height 13pt depth 5pt\ \ $\mu_i=0$ and $\eta_i=\alpha_i$
\hss}
}
\vrule height 13pt depth 5pt
\hrule width\hsize
\noindent
\vrule height 13pt depth 5pt
\vtop{\hsize=\tablesize
\hbox to\hsize{\ \ \ \ \ $m_C=\zeta_i\geqslant 2\,\alpha_i$,
\kern 5.8em
$m_D=\xi_i\geqslant 6\,\alpha_i$,\kern 5.8em 
$m_E=\varkappa_i\geqslant 4\,\alpha_i$
\hss}
}
\vrule height 13pt depth 5pt
\hrule width\hsize
\noindent
\vrule height 13pt depth 5pt
\vtop{\hsize=\tablesize
\hbox to\hsize{\kern 15em 
$\zeta_i=\varkappa_i\leqslant\xi_i$
\hss\vrule height 13pt depth 5pt\ \,$\surd$\ }
}
\vrule height 13pt depth 5pt\vskip -1pt
\noindent
\vrule height 13pt depth 5pt
\vtop{\hsize=\tablesize
\hbox to\hsize{\kern 15em 
$\xi_i=\varkappa_i\leqslant\zeta_i$
\hss\vrule height 13pt depth 5pt\ \,$\surd$\ }
}
\vrule height 13pt depth 5pt\vskip -1pt
\noindent
\vrule height 13pt depth 5pt
\vtop{\hsize=\tablesize
\hbox to\hsize{\kern 15em 
$\xi_i=\zeta_i\leqslant\varkappa_i$
\hss\vrule height 13pt depth 5pt\ \,$\surd$\ }
}
\vrule height 13pt depth 5pt
\hrule width\hsize
}\pagebreak
\noindent\mytable{5.25}\par
\vbox{
\hrule width\hsize
\noindent
\vrule height 13pt depth 5pt
\vtop{\hsize=\tablesize
\hbox to\hsize{\ \ $\alpha_i=\omega_i>\beta_i=0$, $p_i=2$\ \ 
\vrule height 13pt depth 5pt\ \ $\mu_i>0\ \Rightarrow
\ \mu_i=\alpha_i$ and $\eta_i=0$
\hss}
}
\vrule height 13pt depth 5pt
\hrule width\hsize
\noindent
\vrule height 13pt depth 5pt
\vtop{\hsize=\tablesize
\hbox to\hsize{\ \ \ \ \ $m_C=2\,\alpha_i+2\,\mu_i$,\kern 3.4em
$m_D=\xi_i\geqslant 4\,\alpha_i+2\,\eta_i$,\kern 3.4em 
$m_E=\varkappa_i\geqslant 4\,\alpha_i+1$
\hss}
}
\vrule height 13pt depth 5pt
\hrule width\hsize
\noindent
\vrule height 13pt depth 5pt
\vtop{\hsize=\tablesize
\hbox to\hsize{\kern 13em 
$2\,\alpha_i+2\,\mu_i=\varkappa_i\leqslant\xi_i$
\hss\vrule height 13pt depth 5pt\ \,$\phantom{\surd}$\ }
}
\vrule height 13pt depth 5pt\vskip -1pt
\noindent
\vrule height 13pt depth 5pt
\vtop{\hsize=\tablesize
\hbox to\hsize{\kern 13em 
$\xi_i=\varkappa_i\leqslant 2\,\alpha_i+2\,\mu_i$
\hss\vrule height 13pt depth 5pt\ \,$\phantom{\surd}$\ }
}
\vrule height 13pt depth 5pt\vskip -1pt
\noindent
\vrule height 13pt depth 5pt
\vtop{\hsize=\tablesize
\hbox to\hsize{\kern 13em 
$\xi_i=2\,\alpha_i+2\,\mu_i\leqslant\varkappa_i$
\hss\vrule height 13pt depth 5pt\ \,$\surd$\ }
}
\vrule height 13pt depth 5pt
\hrule width\hsize
}
\vskip 1ex plus 3pt
\noindent\mytable{5.26}\par
\vbox{
\hrule width\hsize
\noindent
\vrule height 13pt depth 5pt
\vtop{\hsize=\tablesize
\hbox to\hsize{\ \ $\alpha_i=\omega_i>\beta_i=0$, $p_i=2$\ \ 
\vrule height 13pt depth 5pt\ \ $\mu_i=0$ and $\eta_i=\alpha_i$
\hss}
}
\vrule height 13pt depth 5pt
\hrule width\hsize
\noindent
\vrule height 13pt depth 5pt
\vtop{\hsize=\tablesize
\hbox to\hsize{\ \ \ \ \ $m_C=\zeta_i\geqslant 2\,\alpha_i$,
\kern 5em
$m_D=\xi_i\geqslant 6\,\alpha_i$,\kern 5em 
$m_E=\varkappa_i\geqslant 4\,\alpha_i+1$
\hss}
}
\vrule height 13pt depth 5pt
\hrule width\hsize
\noindent
\vrule height 13pt depth 5pt
\vtop{\hsize=\tablesize
\hbox to\hsize{\kern 15em 
$\zeta_i=\varkappa_i\leqslant\xi_i$
\hss\vrule height 13pt depth 5pt\ \,$\surd$\ }
}
\vrule height 13pt depth 5pt\vskip -1pt
\noindent
\vrule height 13pt depth 5pt
\vtop{\hsize=\tablesize
\hbox to\hsize{\kern 15em 
$\xi_i=\varkappa_i\leqslant\zeta_i$
\hss\vrule height 13pt depth 5pt\ \,$\surd$\ }
}
\vrule height 13pt depth 5pt\vskip -1pt
\noindent
\vrule height 13pt depth 5pt
\vtop{\hsize=\tablesize
\hbox to\hsize{\kern 15em 
$\xi_i=\zeta_i\leqslant\varkappa_i$
\hss\vrule height 13pt depth 5pt\ \,$\surd$\ }
}
\vrule height 13pt depth 5pt
\hrule width\hsize
}
\vskip 1ex plus 3pt
\noindent\mytable{5.27}\par
\vbox{
\hrule width\hsize
\noindent
\vrule height 13pt depth 5pt
\vtop{\hsize=\tablesize
\hbox to\hsize{\ \ $\alpha_i=\omega_i>\beta_i=0$, $p_i=3$\ \ 
\vrule height 13pt depth 5pt\ \ $\mu_i>0\ \Rightarrow
\ \mu_i=\alpha_i$ and $\eta_i=0$
\hss}
}
\vrule height 13pt depth 5pt
\hrule width\hsize
\noindent
\vrule height 13pt depth 5pt
\vtop{\hsize=\tablesize
\hbox to\hsize{\ \ \ \ \ $m_C=2\,\alpha_i+2\,\mu_i$,\kern 5.4em
$m_D=\xi_i\geqslant 4\,\alpha_i+2\,\eta_i$,\kern 5.4em 
$m_E=4\,\alpha_i$
\hss}
}
\vrule height 13pt depth 5pt
\hrule width\hsize
\noindent
\vrule height 13pt depth 5pt
\vtop{\hsize=\tablesize
\hbox to\hsize{\kern 13em 
$2\,\alpha_i+2\,\mu_i=4\,\alpha_i\leqslant\xi_i$
\hss\vrule height 13pt depth 5pt\ \,$\surd$\ }
}
\vrule height 13pt depth 5pt\vskip -1pt
\noindent
\vrule height 13pt depth 5pt
\vtop{\hsize=\tablesize
\hbox to\hsize{\kern 13em 
$\xi_i=4\,\alpha_i\leqslant 2\,\alpha_i+2\,\mu_i$
\hss\vrule height 13pt depth 5pt\ \,$\surd$\ }
}
\vrule height 13pt depth 5pt\vskip -1pt
\noindent
\vrule height 13pt depth 5pt
\vtop{\hsize=\tablesize
\hbox to\hsize{\kern 13em 
$\xi_i=2\,\alpha_i+2\,\mu_i\leqslant 4\,\alpha_i$
\hss\vrule height 13pt depth 5pt\ \,$\surd$\ }
}
\vrule height 13pt depth 5pt
\hrule width\hsize
}
\vskip 1ex plus 3pt
\noindent\mytable{5.28}\par
\vbox{
\hrule width\hsize
\noindent
\vrule height 13pt depth 5pt
\vtop{\hsize=\tablesize
\hbox to\hsize{\ \ $\alpha_i=\omega_i>\beta_i=0$, $p_i=3$\ \ 
\vrule height 13pt depth 5pt\ \ $\mu_i=0$ and $\eta_i=\alpha_i$
\hss}
}
\vrule height 13pt depth 5pt
\hrule width\hsize
\noindent
\vrule height 13pt depth 5pt
\vtop{\hsize=\tablesize
\hbox to\hsize{\ \ \ \ \ $m_C=\zeta_i\geqslant 2\,\alpha_i$,
\kern 7em
$m_D=\xi_i\geqslant 6\,\alpha_i$,\kern 7em 
$m_E=4\,\alpha_i$
\hss}
}
\vrule height 13pt depth 5pt
\hrule width\hsize
\noindent
\vrule height 13pt depth 5pt
\vtop{\hsize=\tablesize
\hbox to\hsize{\kern 15em 
$\zeta_i=4\,\alpha_i\leqslant\xi_i$
\hss\vrule height 13pt depth 5pt\ \,$\surd$\ }
}
\vrule height 13pt depth 5pt\vskip -1pt
\noindent
\vrule height 13pt depth 5pt
\vtop{\hsize=\tablesize
\hbox to\hsize{\kern 15em 
$\xi_i=4\,\alpha_i\leqslant\zeta_i$
\hss\vrule height 13pt depth 5pt\ \,$\phantom{\surd}$\ }
}
\vrule height 13pt depth 5pt\vskip -1pt
\noindent
\vrule height 13pt depth 5pt
\vtop{\hsize=\tablesize
\hbox to\hsize{\kern 15em 
$\xi_i=\zeta_i\leqslant 4\,\alpha_i$
\hss\vrule height 13pt depth 5pt\ \,$\phantom{\surd}$\ }
}
\vrule height 13pt depth 5pt
\hrule width\hsize
}
\vskip 1ex plus 3pt
     In the following cases the multiplicity $\omega_i$ is greater than
the multiplicity $\alpha_i$. Therefore, according to \mythetag{3.3}, we
have $\theta_i=\max(\alpha_i,\beta_i,\omega_i)=\omega_i$.
\vskip 1ex plus 3pt
\noindent\mytable{5.29}\par
\vbox{
\hrule width\hsize
\noindent
\vrule height 13pt depth 5pt
\vtop{\hsize=\tablesize
\hbox to\hsize{\ \ $\omega_i>\alpha_i>\beta_i=0$, $p_i\neq 2$\ \ 
\vrule height 13pt depth 5pt\ \ $\mu_i>0$ and $\eta_i>0\ \Rightarrow
\ \mu_i=\alpha_i$ and $\eta_i=\omega_i-\alpha_i$
\hss}
}
\vrule height 13pt depth 5pt
\hrule width\hsize
\noindent
\vrule height 13pt depth 5pt
\vtop{\hsize=\tablesize
\hbox to\hsize{\ \ \ \ \ $m_C=4\,\omega_i-2\,\alpha_i+2\,\mu_i$,\kern 3.6em
$m_D=2\,\omega_i+2\,\alpha_i+2\,\eta_i$,\kern 3.6em 
$m_E=4\,\omega_i$
\hss}
}
\vrule height 13pt depth 5pt
\hrule width\hsize
\noindent
\vrule height 13pt depth 5pt
\vtop{\hsize=\tablesize
\hbox to\hsize{\kern 8em 
$4\,\omega_i-2\,\alpha_i+2\,\mu_i=4\,\omega_i\leqslant 
2\,\omega_i+2\,\alpha_i+2\,\eta_i$
\hss\vrule height 13pt depth 5pt\ \,$\surd$\ }
}
\vrule height 13pt depth 5pt\vskip -1pt
\noindent
\vrule height 13pt depth 5pt
\vtop{\hsize=\tablesize
\hbox to\hsize{\kern 8em 
$2\,\omega_i+2\,\alpha_i+2\,\eta_i=4\,\omega_i\leqslant 4\,\omega_i
-2\,\alpha_i+2\,\mu_i$
\hss\vrule height 13pt depth 5pt\ \,$\surd$\ }
}
\vrule height 13pt depth 5pt\vskip -1pt
\noindent
\vrule height 13pt depth 5pt
\vtop{\hsize=\tablesize
\hbox to\hsize{\kern 8em 
$2\,\omega_i+2\,\alpha_i+2\,\eta_i=4\,\omega_i-2\,\alpha_i+2\,\mu_i
\leqslant 4\,\omega_i$
\hss\vrule height 13pt depth 5pt\ \,$\surd$\ }
}
\vrule height 13pt depth 5pt
\hrule width\hsize
}\pagebreak
\noindent\mytable{5.30}\par
\vbox{
\hrule width\hsize
\noindent
\vrule height 13pt depth 5pt
\vtop{\hsize=\tablesize
\hbox to\hsize{\ \ $\omega_i>\alpha_i>\beta_i=0$, $p_i\neq 2$\ \ 
\vrule height 13pt depth 5pt\ \ $\mu_i=0$ and $\eta_i=\omega_i$
\hss}
}
\vrule height 13pt depth 5pt
\hrule width\hsize
\noindent
\vrule height 13pt depth 5pt
\vtop{\hsize=\tablesize
\hbox to\hsize{\ \ \ \ \ $m_C=\zeta_i\geqslant 4\,\omega_i-2\,\alpha_i$,
\kern 5.25em
$m_D=4\,\omega_i+2\,\alpha_i$,\kern 5.25em
$m_E=4\,\omega_i$
\hss}
}
\vrule height 13pt depth 5pt
\hrule width\hsize
\noindent
\vrule height 13pt depth 5pt
\vtop{\hsize=\tablesize
\hbox to\hsize{\kern 13em 
$\zeta_i=4\,\omega_i\leqslant 4\,\omega_i+2\,\alpha_i$
\hss\vrule height 13pt depth 5pt\ \,$\surd$\ }
}
\vrule height 13pt depth 5pt\vskip -1pt
\noindent
\vrule height 13pt depth 5pt
\vtop{\hsize=\tablesize
\hbox to\hsize{\kern 13em 
$4\,\omega_i+2\,\alpha_i=4\,\omega_i\leqslant\zeta_i$
\hss\vrule height 13pt depth 5pt\ \,$\phantom{\surd}$\ }
}
\vrule height 13pt depth 5pt\vskip -1pt
\noindent
\vrule height 13pt depth 5pt
\vtop{\hsize=\tablesize
\hbox to\hsize{\kern 13em 
$4\,\omega_i+2\,\alpha_i=\zeta_i\leqslant 4\,\omega_i$
\hss\vrule height 13pt depth 5pt\ \,$\phantom{\surd}$\ }
}
\vrule height 13pt depth 5pt
\hrule width\hsize
}
\vskip 1ex plus 3pt
\noindent\mytable{5.31}\par
\vbox{
\hrule width\hsize
\noindent
\vrule height 13pt depth 5pt
\vtop{\hsize=\tablesize
\hbox to\hsize{\ \ $\omega_i>\alpha_i>\beta_i=0$, $p_i\neq 2$\ \ 
\vrule height 13pt depth 5pt\ \ $\mu_i=\omega_i$ and $\eta_i=0$
\hss}
}
\vrule height 13pt depth 5pt
\hrule width\hsize
\noindent
\vrule height 13pt depth 5pt
\vtop{\hsize=\tablesize
\hbox to\hsize{\ \ \ \ \ $m_C=6\,\omega_i-2\,\alpha_i$,
\kern 5.25em
$m_D=\xi_i\geqslant 2\,\omega_i+2\,\alpha_i$,\kern 5.25em
$m_E=4\,\omega_i$
\hss}
}
\vrule height 13pt depth 5pt
\hrule width\hsize
\noindent
\vrule height 13pt depth 5pt
\vtop{\hsize=\tablesize
\hbox to\hsize{\kern 13em 
$6\,\omega_i-2\,\alpha_i=4\,\omega_i\leqslant\xi_i$
\hss\vrule height 13pt depth 5pt\ \,$\phantom{\surd}$\ }
}
\vrule height 13pt depth 5pt\vskip -1pt
\noindent
\vrule height 13pt depth 5pt
\vtop{\hsize=\tablesize
\hbox to\hsize{\kern 13em 
$\xi_i=4\,\omega_i\leqslant 6\,\omega_i-2\,\alpha_i$
\hss\vrule height 13pt depth 5pt\ \,$\surd$\ }
}
\vrule height 13pt depth 5pt\vskip -1pt
\noindent
\vrule height 13pt depth 5pt
\vtop{\hsize=\tablesize
\hbox to\hsize{\kern 13em 
$\xi_i=6\,\omega_i-2\,\alpha_i\leqslant 4\,\omega_i$
\hss\vrule height 13pt depth 5pt\ \,$\phantom{\surd}$\ }
}
\vrule height 13pt depth 5pt
\hrule width\hsize
}
\vskip 1ex plus 3pt
\noindent\mytable{5.32}\par
\vbox{
\hrule width\hsize
\noindent
\vrule height 13pt depth 5pt
\vtop{\hsize=\tablesize
\hbox to\hsize{\ \ $\omega_i>\alpha_i>\beta_i=0$, $p_i=2$\ \ 
\vrule height 13pt depth 5pt\ \ $\mu_i>0$ and $\eta_i>0\ \Rightarrow
\ \mu_i=\alpha_i$ and $\eta_i=\omega_i-\alpha_i$
\hss}
}
\vrule height 13pt depth 5pt
\hrule width\hsize
\noindent
\vrule height 13pt depth 5pt
\vtop{\hsize=\tablesize
\hbox to\hsize{\ \ \ \ \ $m_C=4\,\omega_i-2\,\alpha_i+2\,\mu_i$,\kern 2.8em
$m_D=2\,\omega_i+2\,\alpha_i+2\,\eta_i$,\kern 2.8em 
$m_E=4\,\omega_i+1$
\hss}
}
\vrule height 13pt depth 5pt
\hrule width\hsize
\noindent
\vrule height 13pt depth 5pt
\vtop{\hsize=\tablesize
\hbox to\hsize{\kern 8em 
$4\,\omega_i-2\,\alpha_i+2\,\mu_i=4\,\omega_i+1\leqslant 
2\,\omega_i+2\,\alpha_i+2\,\eta_i$
\hss\vrule height 13pt depth 5pt\ \,$\phantom{\surd}$\ }
}
\vrule height 13pt depth 5pt\vskip -1pt
\noindent
\vrule height 13pt depth 5pt
\vtop{\hsize=\tablesize
\hbox to\hsize{\kern 8em 
$2\,\omega_i+2\,\alpha_i+2\,\eta_i=4\,\omega_i+1\leqslant 4\,\omega_i
-2\,\alpha_i+2\,\mu_i$
\hss\vrule height 13pt depth 5pt\ \,$\phantom{\surd}$\ }
}
\vrule height 13pt depth 5pt\vskip -1pt
\noindent
\vrule height 13pt depth 5pt
\vtop{\hsize=\tablesize
\hbox to\hsize{\kern 8em 
$2\,\omega_i+2\,\alpha_i+2\,\eta_i=4\,\omega_i-2\,\alpha_i+2\,\mu_i
\leqslant 4\,\omega_i+1$
\hss\vrule height 13pt depth 5pt\ \,$\surd$\ }
}
\vrule height 13pt depth 5pt
\hrule width\hsize
}
\vskip 1ex plus 3pt
\noindent\mytable{5.33}\par
\vbox{
\hrule width\hsize
\noindent
\vrule height 13pt depth 5pt
\vtop{\hsize=\tablesize
\hbox to\hsize{\ \ $\omega_i>\alpha_i>\beta_i=0$, $p_i=2$\ \ 
\vrule height 13pt depth 5pt\ \ $\mu_i=0$ and $\eta_i=\omega_i$
\hss}
}
\vrule height 13pt depth 5pt
\hrule width\hsize
\noindent
\vrule height 13pt depth 5pt
\vtop{\hsize=\tablesize
\hbox to\hsize{\ \ \ \ \ $m_C=\zeta_i\geqslant 4\,\omega_i-2\,\alpha_i$,
\kern 4.5em
$m_D=4\,\omega_i+2\,\alpha_i$,\kern 4.5em
$m_E=4\,\omega_i+1$
\hss}
}
\vrule height 13pt depth 5pt
\hrule width\hsize
\noindent
\vrule height 13pt depth 5pt
\vtop{\hsize=\tablesize
\hbox to\hsize{\kern 12em 
$\zeta_i=4\,\omega_i+1\leqslant 4\,\omega_i+2\,\alpha_i$
\hss\vrule height 13pt depth 5pt\ \,$\surd$\ }
}
\vrule height 13pt depth 5pt\vskip -1pt
\noindent
\vrule height 13pt depth 5pt
\vtop{\hsize=\tablesize
\hbox to\hsize{\kern 12em 
$4\,\omega_i+2\,\alpha_i=4\,\omega_i+1\leqslant\zeta_i$
\hss\vrule height 13pt depth 5pt\ \,$\phantom{\surd}$\ }
}
\vrule height 13pt depth 5pt\vskip -1pt
\noindent
\vrule height 13pt depth 5pt
\vtop{\hsize=\tablesize
\hbox to\hsize{\kern 12em 
$4\,\omega_i+2\,\alpha_i=\zeta_i\leqslant 4\,\omega_i+1$
\hss\vrule height 13pt depth 5pt\ \,$\phantom{\surd}$\ }
}
\vrule height 13pt depth 5pt
\hrule width\hsize
}
\vskip 1ex plus 3pt
\noindent\mytable{5.34}\par
\vbox{
\hrule width\hsize
\noindent
\vrule height 13pt depth 5pt
\vtop{\hsize=\tablesize
\hbox to\hsize{\ \ $\omega_i>\alpha_i>\beta_i=0$, $p_i=2$\ \ 
\vrule height 13pt depth 5pt\ \ $\mu_i=\omega_i$ and $\eta_i=0$
\hss}
}
\vrule height 13pt depth 5pt
\hrule width\hsize
\noindent
\vrule height 13pt depth 5pt
\vtop{\hsize=\tablesize
\hbox to\hsize{\ \ \ \ \ $m_C=6\,\omega_i-2\,\alpha_i$,
\kern 4.3em
$m_D=\xi_i\geqslant 2\,\omega_i+2\,\alpha_i$,\kern 4.3em
$m_E=4\,\omega_i+1$
\hss}
}
\vrule height 13pt depth 5pt
\hrule width\hsize
\noindent
\vrule height 13pt depth 5pt
\vtop{\hsize=\tablesize
\hbox to\hsize{\kern 13em 
$6\,\omega_i-2\,\alpha_i=4\,\omega_i+1\leqslant\xi_i$
\hss\vrule height 13pt depth 5pt\ \,$\phantom{\surd}$\ }
}
\vrule height 13pt depth 5pt\vskip -1pt
\noindent
\vrule height 13pt depth 5pt
\vtop{\hsize=\tablesize
\hbox to\hsize{\kern 13em 
$\xi_i=4\,\omega_i+1\leqslant 6\,\omega_i-2\,\alpha_i$
\hss\vrule height 13pt depth 5pt\ \,$\surd$\ }
}
\vrule height 13pt depth 5pt\vskip -1pt
\noindent
\vrule height 13pt depth 5pt
\vtop{\hsize=\tablesize
\hbox to\hsize{\kern 13em 
$\xi_i=6\,\omega_i-2\,\alpha_i\leqslant 4\,\omega_i+1$
\hss\vrule height 13pt depth 5pt\ \,$\phantom{\surd}$\ }
}
\vrule height 13pt depth 5pt
\hrule width\hsize
}
\vskip 1ex plus 3pt
     The case $p_i=3$ for $\omega_i>\alpha_i>\beta_i=0$ leads to the same 
equalities and inequalities as the other cases in the tables \mythetable{5.29},
\pagebreak\mythetable{5.30}, and \mythetable{5.31}.\par
\noindent\mytable{5.35}\par
\vbox{
\hrule width\hsize
\noindent
\vrule height 13pt depth 5pt
\vtop{\hsize=\tablesize
\hbox to\hsize{\ \ $\omega_i>\alpha_i=\beta_i=0$, $p_i\neq 2$\ \ 
\vrule height 13pt depth 5pt\ \ $\eta_i>0\ \Rightarrow
\ \eta_i=\omega_i$ and $\mu_i=0$
\hss}
}
\vrule height 13pt depth 5pt
\hrule width\hsize
\noindent
\vrule height 13pt depth 5pt
\vtop{\hsize=\tablesize
\hbox to\hsize{\ \ \ \ \ $m_C=\zeta_i\geqslant 
4\,\omega_i+2\,\mu_i$,\kern 4.3em
$m_D=2\,\omega_i+2\,\eta_i$,\kern 4.3em 
$m_E=\varkappa_i\geqslant 4\,\omega_i$
\hss}
}
\vrule height 13pt depth 5pt
\hrule width\hsize
\noindent
\vrule height 13pt depth 5pt
\vtop{\hsize=\tablesize
\hbox to\hsize{\kern 13em 
$\zeta_i=\varkappa_i\leqslant 2\,\omega_i+2\,\eta_i$
\hss\vrule height 13pt depth 5pt\ \,$\surd$\ }
}
\vrule height 13pt depth 5pt\vskip -1pt
\noindent
\vrule height 13pt depth 5pt
\vtop{\hsize=\tablesize
\hbox to\hsize{\kern 13em 
$2\,\omega_i+2\,\eta_i=\varkappa_i\leqslant\zeta_i$
\hss\vrule height 13pt depth 5pt\ \,$\surd$\ }
}
\vrule height 13pt depth 5pt\vskip -1pt
\noindent
\vrule height 13pt depth 5pt
\vtop{\hsize=\tablesize
\hbox to\hsize{\kern 13em 
$2\,\omega_i+2\,\eta_i=\zeta_i\leqslant\varkappa_i$
\hss\vrule height 13pt depth 5pt\ \,$\surd$\ }
}
\vrule height 13pt depth 5pt
\hrule width\hsize
}
\vskip 0.1ex plus 3pt
\noindent\mytable{5.36}\par
\vbox{
\hrule width\hsize
\noindent
\vrule height 13pt depth 5pt
\vtop{\hsize=\tablesize
\hbox to\hsize{\ \ $\omega_i>\alpha_i=\beta_i=0$, $p_i\neq 2$\ \ 
\vrule height 13pt depth 5pt\ \ $\eta_i=0$ and $\mu_i=\omega_i$
\hss}
}
\vrule height 13pt depth 5pt
\hrule width\hsize
\noindent
\vrule height 13pt depth 5pt
\vtop{\hsize=\tablesize
\hbox to\hsize{\ \ \ \ \ $m_C=\zeta_i\geqslant 6\,\omega_i$,
\kern 5.8em
$m_D=\xi_i\geqslant 2\,\omega_i$,\kern 5.8em
$m_E=\varkappa_i\geqslant 4\,\omega_i$
\hss}
}
\vrule height 13pt depth 5pt
\hrule width\hsize
\noindent
\vrule height 13pt depth 5pt
\vtop{\hsize=\tablesize
\hbox to\hsize{\kern 15em 
$\zeta_i=\varkappa_i\leqslant\xi_i$
\hss\vrule height 13pt depth 5pt\ \,$\surd$\ }
}
\vrule height 13pt depth 5pt\vskip -1pt
\noindent
\vrule height 13pt depth 5pt
\vtop{\hsize=\tablesize
\hbox to\hsize{\kern 15em 
$\xi_i=\varkappa_i\leqslant\zeta_i$
\hss\vrule height 13pt depth 5pt\ \,$\surd$\ }
}
\vrule height 13pt depth 5pt\vskip -1pt
\noindent
\vrule height 13pt depth 5pt
\vtop{\hsize=\tablesize
\hbox to\hsize{\kern 15em 
$\xi_i=\zeta_i\leqslant\varkappa_i$
\hss\vrule height 13pt depth 5pt\ \,$\surd$\ }
}
\vrule height 13pt depth 5pt
\hrule width\hsize
}
\vskip 0.1ex plus 3pt
\noindent\mytable{5.37}\par
\vbox{
\hrule width\hsize
\noindent
\vrule height 13pt depth 5pt
\vtop{\hsize=\tablesize
\hbox to\hsize{\ \ $\omega_i>\alpha_i=\beta_i=0$, $p_i=2$\ \ 
\vrule height 13pt depth 5pt\ \ $\eta_i>0\ \Rightarrow
\ \eta_i=\omega_i$ and $\mu_i=0$
\hss}
}
\vrule height 13pt depth 5pt
\hrule width\hsize
\noindent
\vrule height 13pt depth 5pt
\vtop{\hsize=\tablesize
\hbox to\hsize{\ \ \ \ \ $m_C=\zeta_i\geqslant 
4\,\omega_i+2\,\mu_i$,\kern 3.4em
$m_D=2\,\omega_i+2\,\eta_i$,\kern 3.4em 
$m_E=\varkappa_i\geqslant 4\,\omega_i+1$
\hss}
}
\vrule height 13pt depth 5pt
\hrule width\hsize
\noindent
\vrule height 13pt depth 5pt
\vtop{\hsize=\tablesize
\hbox to\hsize{\kern 13em 
$\zeta_i=\varkappa_i\leqslant 2\,\omega_i+2\,\eta_i$
\hss\vrule height 13pt depth 5pt\ \,$\phantom{\surd}$\ }
}
\vrule height 13pt depth 5pt\vskip -1pt
\noindent
\vrule height 13pt depth 5pt
\vtop{\hsize=\tablesize
\hbox to\hsize{\kern 13em 
$2\,\omega_i+2\,\eta_i=\varkappa_i\leqslant\zeta_i$
\hss\vrule height 13pt depth 5pt\ \,$\phantom{\surd}$\ }
}
\vrule height 13pt depth 5pt\vskip -1pt
\noindent
\vrule height 13pt depth 5pt
\vtop{\hsize=\tablesize
\hbox to\hsize{\kern 13em 
$2\,\omega_i+2\,\eta_i=\zeta_i\leqslant\varkappa_i$
\hss\vrule height 13pt depth 5pt\ \,$\surd$\ }
}
\vrule height 13pt depth 5pt
\hrule width\hsize
}
\vskip 0.1ex plus 3pt
\noindent\mytable{5.38}\par
\vbox{
\hrule width\hsize
\noindent
\vrule height 13pt depth 5pt
\vtop{\hsize=\tablesize
\hbox to\hsize{\ \ $\omega_i>\alpha_i=\beta_i=0$, $p_i=2$\ \ 
\vrule height 13pt depth 5pt\ \ $\eta_i=0$ and $\mu_i=\omega_i$
\hss}
}
\vrule height 13pt depth 5pt
\hrule width\hsize
\noindent
\vrule height 13pt depth 5pt
\vtop{\hsize=\tablesize
\hbox to\hsize{\ \ \ \ \ $m_C=\zeta_i\geqslant 6\,\omega_i$,
\kern 5em
$m_D=\xi_i\geqslant 2\,\omega_i$,\kern 5em
$m_E=\varkappa_i\geqslant 4\,\omega_i+1$
\hss}
}
\vrule height 13pt depth 5pt
\hrule width\hsize
\noindent
\vrule height 13pt depth 5pt
\vtop{\hsize=\tablesize
\hbox to\hsize{\kern 15em 
$\zeta_i=\varkappa_i\leqslant\xi_i$
\hss\vrule height 13pt depth 5pt\ \,$\surd$\ }
}
\vrule height 13pt depth 5pt\vskip -1pt
\noindent
\vrule height 13pt depth 5pt
\vtop{\hsize=\tablesize
\hbox to\hsize{\kern 15em 
$\xi_i=\varkappa_i\leqslant\zeta_i$
\hss\vrule height 13pt depth 5pt\ \,$\surd$\ }
}
\vrule height 13pt depth 5pt\vskip -1pt
\noindent
\vrule height 13pt depth 5pt
\vtop{\hsize=\tablesize
\hbox to\hsize{\kern 15em 
$\xi_i=\zeta_i\leqslant\varkappa_i$
\hss\vrule height 13pt depth 5pt\ \,$\surd$\ }
}
\vrule height 13pt depth 5pt
\hrule width\hsize
}
\vskip 0.9ex plus 3pt
    Note that in all of the above cases $\alpha_i\geqslant\beta_i$. The rest
of the cases are those where $\alpha_i<\beta_i$. Let's recall that the 
polynomial $P_{abu}(t)$ in \mythetag{1.1} is invariant with respect to exchanging
parameters $a$ and $b$ (see \mythetag{2.1}). The same is true for 
$P_{\sigma(abu)}(t)$ in \mythetag{2.4}. Therefore we can produce the rest of
the cases from those already considered. 
\vskip 1ex
\noindent\mytable{5.39} (symmetric to the table \mythetable{5.1})\par
\vbox{
\hrule width\hsize
\noindent
\vrule height 13pt depth 5pt
\vtop{\hsize=\tablesize
\hbox to\hsize{\ \ $\beta_i>\alpha_i>\omega_i=0$, $p_i\neq 2$\ \ 
\vrule height 13pt depth 5pt\ \ $\mu_i>0$ and $\eta_i>0\ \Rightarrow
\ \mu_i=\alpha_i$ and $\eta_i=\beta_i-\alpha_i$
\hss}
}
\vrule height 13pt depth 5pt
\hrule width\hsize
\noindent
\vrule height 13pt depth 5pt
\vtop{\hsize=\tablesize
\hbox to\hsize{\ \ \ \ \ $m_C=4\,\beta_i-2\,\alpha_i+2\,\mu_i$,\kern 3.6em
$m_D=2\,\beta_i+2\,\alpha_i+2\,\eta_i$,\kern 3.6em
$m_E=4\,\beta_i$
\hss}
}
\vrule height 13pt depth 5pt
\hrule width\hsize
\noindent
\vrule height 13pt depth 5pt
\vtop{\hsize=\tablesize
\hbox to\hsize{\kern 8em 
$4\,\beta_i-2\,\alpha_i+2\,\mu_i=4\,\beta_i\leqslant 
2\,\beta_i+2\,\alpha_i+2\,\eta_i$
\hss\vrule height 13pt depth 5pt\ \,$\surd$\ }
}
\vrule height 13pt depth 5pt\vskip -1pt
\noindent
\vrule height 13pt depth 5pt
\vtop{\hsize=\tablesize
\hbox to\hsize{\kern 8em 
$2\,\beta_i+2\,\alpha_i+2\,\eta_i=4\,\beta_i\leqslant 
4\,\beta_i-2\,\alpha_i+2\,\mu_i$
\hss\vrule height 13pt depth 5pt\ \,$\surd$\ }
}
\vrule height 13pt depth 5pt\vskip -1pt
\noindent
\vrule height 13pt depth 5pt
\vtop{\hsize=\tablesize
\hbox to\hsize{\kern 8em 
$2\,\beta_i+2\,\alpha_i+2\,\eta_i=4\,\beta_i-2\,\alpha_i+2\,\mu_i
\leqslant 4\,\beta_i$
\hss\vrule height 13pt depth 5pt\ \,$\surd$\ }
}
\vrule height 13pt depth 5pt
\hrule width\hsize
}
\vskip 1ex
\noindent\mytable{5.40} (symmetric to the table \mythetable{5.2})\par
\vbox{
\hrule width\hsize
\noindent
\vrule height 13pt depth 5pt
\vtop{\hsize=\tablesize
\hbox to\hsize{\ \ $\beta_i>\alpha_i>\omega_i=0$, $p_i\neq 2$\ \ 
\vrule height 13pt depth 5pt\ \ $\mu_i=0$ and $\eta_i=\beta_i$
\hss}
}
\vrule height 13pt depth 5pt
\hrule width\hsize
\noindent
\vrule height 13pt depth 5pt
\vtop{\hsize=\tablesize
\hbox to\hsize{\ \ \ \ \ $m_C=\zeta_i\geqslant 4\,\beta_i-2\,\alpha_i$, 
\kern 5.2em
$m_D=4\,\beta_i+2\,\alpha_i$,\kern 5.2em
$m_E=4\,\beta_i$
\hss}
}
\vrule height 13pt depth 5pt
\hrule width\hsize
\noindent
\vrule height 13pt depth 5pt
\vtop{\hsize=\tablesize
\hbox to\hsize{\kern 12em 
$\zeta_i=4\,\beta_i\leqslant 4\,\beta_i+2\,\alpha_i$
\hss\vrule height 13pt depth 5pt\ \,$\surd$\ }
}
\vrule height 13pt depth 5pt\vskip -1pt
\noindent
\vrule height 13pt depth 5pt
\vtop{\hsize=\tablesize
\hbox to\hsize{\kern 12em 
$4\,\beta_i+2\,\alpha_i=4\,\beta_i\leqslant\zeta_i$
\hss\vrule height 13pt depth 5pt\ \,$\phantom{\surd}$\ }
}
\vrule height 13pt depth 5pt\vskip -1pt
\noindent
\vrule height 13pt depth 5pt
\vtop{\hsize=\tablesize
\hbox to\hsize{\kern 12em 
$4\,\beta_i+2\,\alpha_i=\zeta_i\leqslant 4\,\beta_i$
\hss\vrule height 13pt depth 5pt\ \,$\phantom{\surd}$\ }
}
\vrule height 13pt depth 5pt
\hrule width\hsize
}
\vskip 1ex plus 3pt
\noindent\mytable{5.41} (symmetric to the table \mythetable{5.3})\par
\vbox{
\hrule width\hsize
\noindent
\vrule height 13pt depth 5pt
\vtop{\hsize=\tablesize
\hbox to\hsize{\ \ $\beta_i>\alpha_i>\omega_i=0$, $p_i\neq 2$\ \ 
\vrule height 13pt depth 5pt\ \ $\mu_i=\beta_i$ and $\eta_i=0$
\hss}
}
\vrule height 13pt depth 5pt
\hrule width\hsize
\noindent
\vrule height 13pt depth 5pt
\vtop{\hsize=\tablesize
\hbox to\hsize{\ \ \ \ \ $m_C=6\,\beta_i-2\,\alpha_i$, 
\kern 5.2em
$m_D=\xi_i\geqslant 2\,\beta_i+2\,\alpha_i$,\kern 5.2em
$m_E=4\,\beta_i$
\hss}
}
\vrule height 13pt depth 5pt
\hrule width\hsize
\noindent
\vrule height 13pt depth 5pt
\vtop{\hsize=\tablesize
\hbox to\hsize{\kern 12em 
$6\,\beta_i-2\,\alpha_i=4\,\beta_i\leqslant\xi_i$
\hss\vrule height 13pt depth 5pt\ \,$\phantom{\surd}$\ }
}
\vrule height 13pt depth 5pt\vskip -1pt
\noindent
\vrule height 13pt depth 5pt
\vtop{\hsize=\tablesize
\hbox to\hsize{\kern 12em 
$\xi_i=4\,\beta_i\leqslant 6\,\beta_i-2\,\alpha_i$
\hss\vrule height 13pt depth 5pt\ \,$\surd$\ }
}
\vrule height 13pt depth 5pt\vskip -1pt
\noindent
\vrule height 13pt depth 5pt
\vtop{\hsize=\tablesize
\hbox to\hsize{\kern 12em 
$\xi_i=6\,\beta_i-2\,\alpha_i\leqslant 4\,\beta_i$
\hss\vrule height 13pt depth 5pt\ \,$\phantom{\surd}$\ }
}
\vrule height 13pt depth 5pt
\hrule width\hsize
}
\vskip 1ex plus 3pt
\noindent\mytable{5.42} (symmetric to the table \mythetable{5.4})\par
\vbox{
\hrule width\hsize
\noindent
\vrule height 13pt depth 5pt
\vtop{\hsize=\tablesize
\hbox to\hsize{\ \ $\beta_i>\alpha_i>\omega_i=0$, $p_i=2$\ \ 
\vrule height 13pt depth 5pt\ \ $\mu_i>0$ and $\eta_i>0\ \Rightarrow
\ \mu_i=\alpha_i$ and $\eta_i=\beta_i-\alpha_i$
\hss}
}
\vrule height 13pt depth 5pt
\hrule width\hsize
\noindent
\vrule height 13pt depth 5pt
\vtop{\hsize=\tablesize
\hbox to\hsize{\ \ \ \ \ $m_C=4\,\beta_i-2\,\alpha_i+2\,\mu_i$,\kern 3em
$m_D=2\,\beta_i+2\,\alpha_i+2\,\eta_i$,\kern 3em
$m_E=4\,\beta_i+1$
\hss}
}
\vrule height 13pt depth 5pt
\hrule width\hsize
\noindent
\vrule height 13pt depth 5pt
\vtop{\hsize=\tablesize
\hbox to\hsize{\kern 7em 
$4\,\beta_i-2\,\alpha_i+2\,\mu_i=4\,\beta_i+1\leqslant 
2\,\beta_i+2\,\alpha_i+2\,\eta_i$
\hss\vrule height 13pt depth 5pt\ \,$\phantom{\surd}$\ }
}
\vrule height 13pt depth 5pt\vskip -1pt
\noindent
\vrule height 13pt depth 5pt
\vtop{\hsize=\tablesize
\hbox to\hsize{\kern 7em 
$2\,\beta_i+2\,\alpha_i+2\,\eta_i=4\,\beta_i+1\leqslant 
4\,\beta_i-2\,\alpha_i+2\,\mu_i$
\hss\vrule height 13pt depth 5pt\ \,$\phantom{\surd}$\ }
}
\vrule height 13pt depth 5pt\vskip -1pt
\noindent
\vrule height 13pt depth 5pt
\vtop{\hsize=\tablesize
\hbox to\hsize{\kern 7em 
$2\,\beta_i+2\,\alpha_i+2\,\eta_i=4\,\beta_i-2\,\alpha_i+2\,\mu_i
\leqslant 4\,\beta_i+1$
\hss\vrule height 13pt depth 5pt\ \,$\surd$\ }
}
\vrule height 13pt depth 5pt
\hrule width\hsize
}
\vskip 1ex plus 3pt
\noindent\mytable{5.43} (symmetric to the table \mythetable{5.5})\par
\vbox{
\hrule width\hsize
\noindent
\vrule height 13pt depth 5pt
\vtop{\hsize=\tablesize
\hbox to\hsize{\ \ $\beta_i>\alpha_i>\omega_i=0$, $p_i=2$\ \ 
\vrule height 13pt depth 5pt\ \ $\mu_i=0$ and $\eta_i=\beta_i$
\hss}
}
\vrule height 13pt depth 5pt
\hrule width\hsize
\noindent
\vrule height 13pt depth 5pt
\vtop{\hsize=\tablesize
\hbox to\hsize{\ \ \ \ \ $m_C=\zeta_i\geqslant 4\,\beta_i-2\,\alpha_i$, 
\kern 4.5em
$m_D=4\,\beta_i+2\,\alpha_i$,\kern 4.5em
$m_E=4\,\beta_i+1$
\hss}
}
\vrule height 13pt depth 5pt
\hrule width\hsize
\noindent
\vrule height 13pt depth 5pt
\vtop{\hsize=\tablesize
\hbox to\hsize{\kern 12em 
$\zeta_i=4\,\beta_i+1\leqslant 4\,\beta_i+2\,\alpha_i$
\hss\vrule height 13pt depth 5pt\ \,$\surd$\ }
}
\vrule height 13pt depth 5pt\vskip -1pt
\noindent
\vrule height 13pt depth 5pt
\vtop{\hsize=\tablesize
\hbox to\hsize{\kern 12em 
$4\,\beta_i+2\,\alpha_i=4\,\beta_i+1\leqslant\zeta_i$
\hss\vrule height 13pt depth 5pt\ \,$\phantom{\surd}$\ }
}
\vrule height 13pt depth 5pt\vskip -1pt
\noindent
\vrule height 13pt depth 5pt
\vtop{\hsize=\tablesize
\hbox to\hsize{\kern 12em 
$4\,\beta_i+2\,\alpha_i=\zeta_i\leqslant 4\,\beta_i+1$
\hss\vrule height 13pt depth 5pt\ \,$\phantom{\surd}$\ }
}
\vrule height 13pt depth 5pt
\hrule width\hsize
}
\vskip 1ex plus 3pt
\noindent\mytable{5.44} (symmetric to the table \mythetable{5.6})\par
\vbox{
\hrule width\hsize
\noindent
\vrule height 13pt depth 5pt
\vtop{\hsize=\tablesize
\hbox to\hsize{\ \ $\beta_i>\alpha_i>\omega_i=0$, $p_i=2$\ \ 
\vrule height 13pt depth 5pt\ \ $\mu_i=\beta_i$ and $\eta_i=0$
\hss}
}
\vrule height 13pt depth 5pt
\hrule width\hsize
\noindent
\vrule height 13pt depth 5pt
\vtop{\hsize=\tablesize
\hbox to\hsize{\ \ \ \ \ $m_C=6\,\beta_i-2\,\alpha_i$, 
\kern 4.5em
$m_D=\xi_i\geqslant 2\,\beta_i+2\,\alpha_i$,\kern 4.5em
$m_E=4\,\beta_i+1$
\hss}
}
\vrule height 13pt depth 5pt	
\hrule width\hsize
\noindent
\vrule height 13pt depth 5pt
\vtop{\hsize=\tablesize
\hbox to\hsize{\kern 12em 
$6\,\beta_i-2\,\alpha_i=4\,\beta_i+1\leqslant\xi_i$
\hss\vrule height 13pt depth 5pt\ \,$\phantom{\surd}$\ }
}
\vrule height 13pt depth 5pt\vskip -1pt
\noindent
\vrule height 13pt depth 5pt
\vtop{\hsize=\tablesize
\hbox to\hsize{\kern 12em 
$\xi_i=4\,\beta_i+1\leqslant 6\,\beta_i-2\,\alpha_i$
\hss\vrule height 13pt depth 5pt\ \,$\surd$\ }
}
\vrule height 13pt depth 5pt\vskip -1pt
\noindent
\vrule height 13pt depth 5pt
\vtop{\hsize=\tablesize
\hbox to\hsize{\kern 12em 
$\xi_i=6\,\beta_i-2\,\alpha_i\leqslant 4\,\beta_i+1$
\hss\vrule height 13pt depth 5pt\ \,$\phantom{\surd}$\ }
}
\vrule height 13pt depth 5pt
\hrule width\hsize
}
\vskip 1ex plus 3pt
     In the following cases the multiplicity $\alpha_i$ coincides with the 
multiplicity $\omega_i=0$. I this case, according to \mythetag{3.3}, we
have $\theta_i=\max(\alpha_i,\beta_i,\omega_i)=\beta_i$.\par
\noindent\mytable{5.45} (symmetric to the table \mythetable{5.7})\par
\vbox{
\hrule width\hsize
\noindent
\vrule height 13pt depth 5pt
\vtop{\hsize=\tablesize
\hbox to\hsize{\ \ $\beta_i>\alpha_i=\omega_i=0$, $p_i\neq 2,3$\ \ 
\vrule height 13pt depth 5pt\ \ $\eta_i>0 \Rightarrow \eta_i=\beta_i$
and $\mu_i=0$
\hss}
}
\vrule height 13pt depth 5pt
\hrule width\hsize
\noindent
\vrule height 13pt depth 5pt
\vtop{\hsize=\tablesize
\hbox to\hsize{\ \ \ \ \ $m_C=\zeta_i\geqslant 4\,\beta_i+2\,\mu_i$, 
\kern 4.1em
$m_D=2\,\beta_i+2\,\eta_i$,\kern 4.1em
$m_E=\varkappa_i\geqslant 4\,\beta_i$
\hss}
}
\vrule height 13pt depth 5pt
\hrule width\hsize
\noindent
\vrule height 13pt depth 5pt
\vtop{\hsize=\tablesize
\hbox to\hsize{\kern 12.5em 
$\varkappa_i=\zeta_i\leqslant 2\,\beta_i+2\,\eta_i$
\hss\vrule height 13pt depth 5pt\ \,$\surd$\ }
}
\vrule height 13pt depth 5pt\vskip -1pt
\noindent
\vrule height 13pt depth 5pt
\vtop{\hsize=\tablesize
\hbox to\hsize{\kern 12.5em 
$\varkappa_i=2\,\beta_i+2\,\eta_i\leqslant\zeta_i$
\hss\vrule height 13pt depth 5pt\ \,$\surd$\ }
}
\vrule height 13pt depth 5pt\vskip -1pt
\noindent
\vrule height 13pt depth 5pt
\vtop{\hsize=\tablesize
\hbox to\hsize{\kern 12.5em 
$\zeta_i=2\,\beta_i+2\,\eta_i\leqslant\varkappa_i$
\hss\vrule height 13pt depth 5pt\ \,$\surd$\ }
}
\vrule height 13pt depth 5pt
\hrule width\hsize
}
\vskip 1ex plus 3pt
\noindent\mytable{5.46} (symmetric to the table \mythetable{5.8})\par
\vbox{
\hrule width\hsize
\noindent
\vrule height 13pt depth 5pt
\vtop{\hsize=\tablesize
\hbox to\hsize{\ \ $\beta_i>\alpha_i=\omega_i=0$, $p_i\neq 2,3$\ \ 
\vrule height 13pt depth 5pt\ \ $\eta_i=0$ and $\mu_i=\beta_i$
\hss}
}
\vrule height 13pt depth 5pt
\hrule width\hsize
\noindent
\vrule height 13pt depth 5pt
\vtop{\hsize=\tablesize
\hbox to\hsize{\ \ \ \ \ $m_C=\zeta_i\geqslant 6\,\beta_i$, 
\kern 5.8em
$m_D=\xi_i\geqslant 2\,\beta_i$,\kern 5.8em
$m_E=\varkappa_i\geqslant 4\,\beta_i$
\hss}
}
\vrule height 13pt depth 5pt
\hrule width\hsize
\noindent
\vrule height 13pt depth 5pt
\vtop{\hsize=\tablesize
\hbox to\hsize{\kern 15em 
$\varkappa_i=\zeta_i\leqslant\xi_i$
\hss\vrule height 13pt depth 5pt\ \,$\surd$\ }
}
\vrule height 13pt depth 5pt\vskip -1pt
\noindent
\vrule height 13pt depth 5pt
\vtop{\hsize=\tablesize
\hbox to\hsize{\kern 15em 
$\varkappa_i=\xi_i\leqslant\zeta_i$
\hss\vrule height 13pt depth 5pt\ \,$\surd$\ }
}
\vrule height 13pt depth 5pt\vskip -1pt
\noindent
\vrule height 13pt depth 5pt
\vtop{\hsize=\tablesize
\hbox to\hsize{\kern 15em 
$\zeta_i=\xi_i\leqslant\varkappa_i$
\hss\vrule height 13pt depth 5pt\ \,$\surd$\ }
}
\vrule height 13pt depth 5pt
\hrule width\hsize
}
\vskip 1ex
     The following tables correspond to the special values of the prime
factor $p_i$, i\.\,e\. to $p_i=2$ and to $p_i=3$.\par 
\vskip 1ex plus 3pt
\noindent\mytable{5.47} (symmetric to the table \mythetable{5.9})\par
\vbox{
\hrule width\hsize
\noindent
\vrule height 13pt depth 5pt
\vtop{\hsize=\tablesize
\hbox to\hsize{\ \ $\beta_i>\alpha_i=\omega_i=0$, $p_i=2$\ \ 
\vrule height 13pt depth 5pt\ \ $\eta_i>0 \Rightarrow \eta_i=\beta_i$
and $\mu_i=0$
\hss}
}
\vrule height 13pt depth 5pt
\hrule width\hsize
\noindent
\vrule height 13pt depth 5pt
\vtop{\hsize=\tablesize
\hbox to\hsize{\ \ \ \ \ $m_C=\zeta_i\geqslant 4\,\beta_i+2\,\mu_i$, 
\kern 3.3em
$m_D=2\,\beta_i+2\,\eta_i$,\kern 3.3em
$m_E=\varkappa_i\geqslant 4\,\beta_i+1$
\hss}
}
\vrule height 13pt depth 5pt
\hrule width\hsize
\noindent
\vrule height 13pt depth 5pt
\vtop{\hsize=\tablesize
\hbox to\hsize{\kern 13em 
$\varkappa_i=\zeta_i\leqslant 2\,\beta_i+2\,\eta_i$
\hss\vrule height 13pt depth 5pt\ \,$\phantom{\surd}$\ }
}
\vrule height 13pt depth 5pt\vskip -1pt
\noindent
\vrule height 13pt depth 5pt
\vtop{\hsize=\tablesize
\hbox to\hsize{\kern 13em 
$\varkappa_i=2\,\beta_i+2\,\eta_i\leqslant\zeta_i$
\hss\vrule height 13pt depth 5pt\ \,$\phantom{\surd}$\ }
}
\vrule height 13pt depth 5pt\vskip -1pt
\noindent
\vrule height 13pt depth 5pt
\vtop{\hsize=\tablesize
\hbox to\hsize{\kern 13em 
$\zeta_i=2\,\beta_i+2\,\eta_i\leqslant\varkappa_i$
\hss\vrule height 13pt depth 5pt\ \,$\surd$\ }
}
\vrule height 13pt depth 5pt
\hrule width\hsize
}
\vskip 1ex plus 3pt
\noindent\mytable{5.48} (symmetric to the table \mythetable{5.10})\par
\vbox{
\hrule width\hsize
\noindent
\vrule height 13pt depth 5pt
\vtop{\hsize=\tablesize
\hbox to\hsize{\ \ $\beta_i>\alpha_i=\omega_i=0$, $p_i=2$\ \ 
\vrule height 13pt depth 5pt\ \ $\eta_i=0$ and $\mu_i=\beta_i$
\hss}
}
\vrule height 13pt depth 5pt
\hrule width\hsize
\noindent
\vrule height 13pt depth 5pt
\vtop{\hsize=\tablesize
\hbox to\hsize{\ \ \ \ \ $m_C=\zeta_i\geqslant 6\,\beta_i$, 
\kern 5em
$m_D=\xi_i\geqslant 2\,\beta_i$,\kern 5em
$m_E=\varkappa_i\geqslant 4\,\beta_i+1$
\hss}
}
\vrule height 13pt depth 5pt
\hrule width\hsize
\noindent
\vrule height 13pt depth 5pt
\vtop{\hsize=\tablesize
\hbox to\hsize{\kern 15em 
$\varkappa_i=\zeta_i\leqslant\xi_i$
\hss\vrule height 13pt depth 5pt\ \,$\surd$\ }
}
\vrule height 13pt depth 5pt\vskip -1pt
\noindent
\vrule height 13pt depth 5pt
\vtop{\hsize=\tablesize
\hbox to\hsize{\kern 15em 
$\varkappa_i=\xi_i\leqslant\zeta_i$
\hss\vrule height 13pt depth 5pt\ \,$\surd$\ }
}
\vrule height 13pt depth 5pt\vskip -1pt
\noindent
\vrule height 13pt depth 5pt
\vtop{\hsize=\tablesize
\hbox to\hsize{\kern 15em 
$\zeta_i=\xi_i\leqslant\varkappa_i$
\hss\vrule height 13pt depth 5pt\ \,$\surd$\ }
}
\vrule height 13pt depth 5pt
\hrule width\hsize
}
\vskip 1ex plus 3pt
\noindent\mytable{5.49} (symmetric to the table \mythetable{5.11})\par
\vbox{
\hrule width\hsize
\noindent
\vrule height 13pt depth 5pt
\vtop{\hsize=\tablesize
\hbox to\hsize{\ \ $\beta_i>\alpha_i=\omega_i=0$, $p_i=3$\ \ 
\vrule height 13pt depth 5pt\ \ $\eta_i>0 \Rightarrow \eta_i=\beta_i$
and $\mu_i=0$
\hss}
}
\vrule height 13pt depth 5pt
\hrule width\hsize
\noindent
\vrule height 13pt depth 5pt
\vtop{\hsize=\tablesize
\hbox to\hsize{\ \ \ \ \ $m_C=\zeta_i\geqslant 4\,\beta_i+2\,\mu_i$, 
\kern 5.3em
$m_D=2\,\beta_i+2\,\eta_i$,\kern 5.3em
$m_E=4\,\beta_i$
\hss}
}
\vrule height 13pt depth 5pt
\hrule width\hsize
\noindent
\vrule height 13pt depth 5pt
\vtop{\hsize=\tablesize
\hbox to\hsize{\kern 13em 
$4\,\beta_i=\zeta_i\leqslant 2\,\beta_i+2\,\eta_i$
\hss\vrule height 13pt depth 5pt\ \,$\surd$\ }
}
\vrule height 13pt depth 5pt\vskip -1pt
\noindent
\vrule height 13pt depth 5pt
\vtop{\hsize=\tablesize
\hbox to\hsize{\kern 13em 
$4\,\beta_i=2\,\beta_i+2\,\eta_i\leqslant\zeta_i$
\hss\vrule height 13pt depth 5pt\ \,$\surd$\ }
}
\vrule height 13pt depth 5pt\vskip -1pt
\noindent
\vrule height 13pt depth 5pt
\vtop{\hsize=\tablesize
\hbox to\hsize{\kern 13em 
$\zeta_i=2\,\beta_i+2\,\eta_i\leqslant4\,\beta_i$
\hss\vrule height 13pt depth 5pt\ \,$\surd$\ }
}
\vrule height 13pt depth 5pt
\hrule width\hsize
}\pagebreak
\noindent\mytable{5.50} (symmetric to the table \mythetable{5.12})\par
\vbox{
\hrule width\hsize
\noindent
\vrule height 13pt depth 5pt
\vtop{\hsize=\tablesize
\hbox to\hsize{\ \ $\beta_i>\alpha_i=\omega_i=0$, $p_i=3$\ \ 
\vrule height 13pt depth 5pt\ \ $\eta_i=0$ and $\mu_i=\beta_i$
\hss}
}
\vrule height 13pt depth 5pt
\hrule width\hsize
\noindent
\vrule height 13pt depth 5pt
\vtop{\hsize=\tablesize
\hbox to\hsize{\ \ \ \ \ $m_C=\zeta_i\geqslant 6\,\beta_i$, 
\kern 7em
$m_D=\xi_i\geqslant 2\,\beta_i$,\kern 7em
$m_E=4\,\beta_i$
\hss}
}
\vrule height 13pt depth 5pt
\hrule width\hsize
\noindent
\vrule height 13pt depth 5pt
\vtop{\hsize=\tablesize
\hbox to\hsize{\kern 15em 
$4\,\beta_i=\zeta_i\leqslant\xi_i$
\hss\vrule height 13pt depth 5pt\ \,$\phantom{\surd}$\ }
}
\vrule height 13pt depth 5pt\vskip -1pt
\noindent
\vrule height 13pt depth 5pt
\vtop{\hsize=\tablesize
\hbox to\hsize{\kern 15em 
$4\,\beta_i=\xi_i\leqslant\zeta_i$
\hss\vrule height 13pt depth 5pt\ \,$\surd$\ }
}
\vrule height 13pt depth 5pt\vskip -1pt
\noindent
\vrule height 13pt depth 5pt
\vtop{\hsize=\tablesize
\hbox to\hsize{\kern 15em 
$\zeta_i=\xi_i\leqslant4\,\beta_i$
\hss\vrule height 13pt depth 5pt\ \,$\phantom{\surd}$\ }
}
\vrule height 13pt depth 5pt
\hrule width\hsize
}
\vskip 1ex plus 3pt
     In the following cases the multiplicity $\omega_i$ is greater than
the multiplicity $\alpha_i$, but less than $\beta_i$. In this case, according 
to \mythetag{3.3}, we have $\theta_i=\max(\alpha_i,\beta_i,\omega_i)=\beta_i$.
\vskip 1ex plus 3pt
\noindent\mytable{5.51} (symmetric to the table \mythetable{5.17})\par
\vbox{
\hrule width\hsize
\noindent
\vrule height 13pt depth 5pt
\vtop{\hsize=\tablesize
\hbox to\hsize{\ \ $\beta_i>\omega_i>\alpha_i=0$, $p_i\neq 2,3$\ \ 
\vrule height 13pt depth 5pt\ \ $\mu_i>0$ and $\eta_i>0\ \Rightarrow
\ \mu_i=\omega_i$ and $\eta_i=\beta_i-\omega_i$
\hss}
}
\vrule height 13pt depth 5pt
\hrule width\hsize
\noindent
\vrule height 13pt depth 5pt
\vtop{\hsize=\tablesize
\hbox to\hsize{\ \ \ \ \ $m_C=4\,\beta_i-2\,\omega_i+2\,\mu_i$,\kern 3.6em
$m_D=2\,\beta_i+2\,\omega_i+2\,\eta_i$,\kern 3.6em 
$m_E=4\,\beta_i$
\hss}
}
\vrule height 13pt depth 5pt
\hrule width\hsize
\noindent
\vrule height 13pt depth 5pt
\vtop{\hsize=\tablesize
\hbox to\hsize{\kern 8em 
$4\,\beta_i-2\,\omega_i+2\,\mu_i=4\,\beta_i\leqslant 2\,\beta_i
+2\,\omega_i+2\,\eta_i$
\hss\vrule height 13pt depth 5pt\ \,$\surd$\ }
}
\vrule height 13pt depth 5pt\vskip -1pt
\noindent
\vrule height 13pt depth 5pt
\vtop{\hsize=\tablesize
\hbox to\hsize{\kern 8em 
$2\,\beta_i+2\,\omega_i+2\,\eta_i=4\,\beta_i\leqslant 4\,\beta_i
-2\,\omega_i+2\,\mu_i$
\hss\vrule height 13pt depth 5pt\ \,$\surd$\ }
}
\vrule height 13pt depth 5pt\vskip -1pt
\noindent
\vrule height 13pt depth 5pt
\vtop{\hsize=\tablesize
\hbox to\hsize{\kern 8em 
$2\,\beta_i+2\,\omega_i+2\,\eta_i=4\,\beta_i-2\,\omega_i+2\,\mu_i
\leqslant 4\,\beta_i$
\hss\vrule height 13pt depth 5pt\ \,$\surd$\ }
}
\vrule height 13pt depth 5pt
\hrule width\hsize
}
\vskip 1ex plus 3pt
\noindent\mytable{5.52} (symmetric to the table \mythetable{5.18})\par
\vbox{
\hrule width\hsize
\noindent
\vrule height 13pt depth 5pt
\vtop{\hsize=\tablesize
\hbox to\hsize{\ \ $\beta_i>\omega_i>\alpha_i=0$, $p_i\neq 2,3$\ \ 
\vrule height 13pt depth 5pt\ \ $\mu_i=0$ and $\eta_i=\beta_i$
\hss}
}
\vrule height 13pt depth 5pt
\hrule width\hsize
\noindent
\vrule height 13pt depth 5pt
\vtop{\hsize=\tablesize
\hbox to\hsize{\ \ \ \ \ $m_C=\zeta_i\geqslant 4\,\beta_i-2\,\omega_i$,
\kern 5.25em
$m_D=4\,\beta_i+2\,\omega_i$,\kern 5.25em
$m_E=4\,\beta_i$
\hss}
}
\vrule height 13pt depth 5pt
\hrule width\hsize
\noindent
\vrule height 13pt depth 5pt
\vtop{\hsize=\tablesize
\hbox to\hsize{\kern 13em 
$\zeta_i=4\,\beta_i\leqslant 4\,\beta_i+2\,\omega_i$
\hss\vrule height 13pt depth 5pt\ \,$\surd$\ }
}
\vrule height 13pt depth 5pt\vskip -1pt
\noindent
\vrule height 13pt depth 5pt
\vtop{\hsize=\tablesize
\hbox to\hsize{\kern 13em 
$4\,\beta_i+2\,\omega_i=4\,\beta_i\leqslant\zeta_i$
\hss\vrule height 13pt depth 5pt\ \,$\phantom{\surd}$\ }
}
\vrule height 13pt depth 5pt\vskip -1pt
\noindent
\vrule height 13pt depth 5pt
\vtop{\hsize=\tablesize
\hbox to\hsize{\kern 13em 
$4\,\beta_i+2\,\omega_i=\zeta_i\leqslant 4\,\beta_i$
\hss\vrule height 13pt depth 5pt\ \,$\phantom{\surd}$\ }
}
\vrule height 13pt depth 5pt
\hrule width\hsize
}
\vskip 1ex plus 3pt
\noindent\mytable{5.53} (symmetric to the table \mythetable{5.19})\par
\vbox{
\hrule width\hsize
\noindent
\vrule height 13pt depth 5pt
\vtop{\hsize=\tablesize
\hbox to\hsize{\ \ $\beta_i>\omega_i>\alpha_i=0$, $p_i\neq 2,3$\ \ 
\vrule height 13pt depth 5pt\ \ $\mu_i=\beta_i$ and $\eta_i=0$
\hss}
}
\vrule height 13pt depth 5pt
\hrule width\hsize
\noindent
\vrule height 13pt depth 5pt
\vtop{\hsize=\tablesize
\hbox to\hsize{\ \ \ \ \ $m_C=6\,\beta_i-2\,\omega_i$,
\kern 5.25em
$m_D=\xi_i\geqslant 2\,\beta_i+2\,\omega_i$,\kern 5.25em
$m_E=4\,\beta_i$
\hss}
}
\vrule height 13pt depth 5pt
\hrule width\hsize
\noindent
\vrule height 13pt depth 5pt
\vtop{\hsize=\tablesize
\hbox to\hsize{\kern 13em 
$6\,\beta_i-2\,\omega_i=4\,\beta_i\leqslant\xi_i$
\hss\vrule height 13pt depth 5pt\ \,$\phantom{\surd}$\ }
}
\vrule height 13pt depth 5pt\vskip -1pt
\noindent
\vrule height 13pt depth 5pt
\vtop{\hsize=\tablesize
\hbox to\hsize{\kern 13em 
$\xi_i=4\,\beta_i\leqslant 6\,\beta_i-2\,\omega_i$
\hss\vrule height 13pt depth 5pt\ \,$\surd$\ }
}
\vrule height 13pt depth 5pt\vskip -1pt
\noindent
\vrule height 13pt depth 5pt
\vtop{\hsize=\tablesize
\hbox to\hsize{\kern 13em 
$\xi_i=6\,\beta_i-2\,\omega_i\leqslant 4\,\beta_i$
\hss\vrule height 13pt depth 5pt\ \,$\phantom{\surd}$\ }
}
\vrule height 13pt depth 5pt
\hrule width\hsize
}
\noindent\mytable{5.54} (symmetric to the table \mythetable{5.20})\par
\vbox{
\hrule width\hsize
\noindent
\vrule height 13pt depth 5pt
\vtop{\hsize=\tablesize
\hbox to\hsize{\ \ $\beta_i>\omega_i>\alpha_i=0$, $p_i=2,3$\ \ 
\vrule height 13pt depth 5pt\ \ $\mu_i>0$ and $\eta_i>0\ \Rightarrow
\ \mu_i=\omega_i$ and $\eta_i=\beta_i-\omega_i$
\hss}
}
\vrule height 13pt depth 5pt
\hrule width\hsize
\noindent
\vrule height 13pt depth 5pt
\vtop{\hsize=\tablesize
\hbox to\hsize{\ \ \ \ \ $m_C=4\,\beta_i-2\,\omega_i+2\,\mu_i$,\kern 2.9em
$m_D=2\,\beta_i+2\,\omega_i+2\,\eta_i$,\kern 2.9em 
$m_E=4\,\beta_i+1$
\hss}
}
\vrule height 13pt depth 5pt
\hrule width\hsize
\noindent
\vrule height 13pt depth 5pt
\vtop{\hsize=\tablesize
\hbox to\hsize{\kern 8em 
$4\,\beta_i-2\,\omega_i+2\,\mu_i=4\,\beta_i+1\leqslant 2\,\beta_i
+2\,\omega_i+2\,\eta_i$
\hss\vrule height 13pt depth 5pt\ \,$\phantom{\surd}$\ }
}
\vrule height 13pt depth 5pt\vskip -1pt
\noindent
\vrule height 13pt depth 5pt
\vtop{\hsize=\tablesize
\hbox to\hsize{\kern 8em 
$2\,\beta_i+2\,\omega_i+2\,\eta_i=4\,\beta_i+1\leqslant 4\,\beta_i
-2\,\omega_i+2\,\mu_i$
\hss\vrule height 13pt depth 5pt\ \,$\phantom{\surd}$\ }
}
\vrule height 13pt depth 5pt\vskip -1pt
\noindent
\vrule height 13pt depth 5pt
\vtop{\hsize=\tablesize
\hbox to\hsize{\kern 8em 
$2\,\beta_i+2\,\omega_i+2\,\eta_i=4\,\beta_i-2\,\omega_i+2\,\mu_i
\leqslant 4\,\beta_i+1$
\hss\vrule height 13pt depth 5pt\ \,$\surd$\ }
}
\vrule height 13pt depth 5pt
\hrule width\hsize
}\pagebreak
\noindent\mytable{5.55} (symmetric to the table \mythetable{5.21})\par
\vbox{
\hrule width\hsize
\noindent
\vrule height 13pt depth 5pt
\vtop{\hsize=\tablesize
\hbox to\hsize{\ \ $\beta_i>\omega_i>\alpha_i=0$, $p_i=2,3$\ \ 
\vrule height 13pt depth 5pt\ \ $\mu_i=0$ and $\eta_i=\beta_i$
\hss}
}
\vrule height 13pt depth 5pt
\hrule width\hsize
\noindent
\vrule height 13pt depth 5pt
\vtop{\hsize=\tablesize
\hbox to\hsize{\ \ \ \ \ $m_C=\zeta_i\geqslant 4\,\beta_i-2\,\omega_i$,
\kern 4.4em
$m_D=4\,\beta_i+2\,\omega_i$,\kern 4.4em
$m_E=4\,\beta_i+1$
\hss}
}
\vrule height 13pt depth 5pt
\hrule width\hsize
\noindent
\vrule height 13pt depth 5pt
\vtop{\hsize=\tablesize
\hbox to\hsize{\kern 12em 
$\zeta_i=4\,\beta_i+1\leqslant 4\,\beta_i+2\,\omega_i$
\hss\vrule height 13pt depth 5pt\ \,$\surd$\ }
}
\vrule height 13pt depth 5pt\vskip -1pt
\noindent
\vrule height 13pt depth 5pt
\vtop{\hsize=\tablesize
\hbox to\hsize{\kern 12em 
$4\,\beta_i+2\,\omega_i=4\,\beta_i+1\leqslant\zeta_i$
\hss\vrule height 13pt depth 5pt\ \,$\phantom{\surd}$\ }
}
\vrule height 13pt depth 5pt\vskip -1pt
\noindent
\vrule height 13pt depth 5pt
\vtop{\hsize=\tablesize
\hbox to\hsize{\kern 12em 
$4\,\beta_i+2\,\omega_i=\zeta_i\leqslant 4\,\beta_i+1$
\hss\vrule height 13pt depth 5pt\ \,$\phantom{\surd}$\ }
}
\vrule height 13pt depth 5pt
\hrule width\hsize
}
\vskip 1ex plus 3pt
\noindent\mytable{5.56} (symmetric to the table \mythetable{5.22})\par
\vbox{
\hrule width\hsize
\noindent
\vrule height 13pt depth 5pt
\vtop{\hsize=\tablesize
\hbox to\hsize{\ \ $\beta_i>\omega_i>\alpha_i=0$, $p_i=2,3$\ \ 
\vrule height 13pt depth 5pt\ \ $\mu_i=\beta_i$ and $\eta_i=0$
\hss}
}
\vrule height 13pt depth 5pt
\hrule width\hsize
\noindent
\vrule height 13pt depth 5pt
\vtop{\hsize=\tablesize
\hbox to\hsize{\ \ \ \ \ $m_C=6\,\beta_i-2\,\omega_i$,
\kern 4.5em
$m_D=\xi_i\geqslant 2\,\beta_i+2\,\omega_i$,\kern 4.5em
$m_E=4\,\beta_i+1$
\hss}
}
\vrule height 13pt depth 5pt
\hrule width\hsize
\noindent
\vrule height 13pt depth 5pt
\vtop{\hsize=\tablesize
\hbox to\hsize{\kern 12em 
$6\,\beta_i-2\,\omega_i=4\,\beta_i+1\leqslant\xi_i$
\hss\vrule height 13pt depth 5pt\ \,$\phantom{\surd}$\ }
}
\vrule height 13pt depth 5pt\vskip -1pt
\noindent
\vrule height 13pt depth 5pt
\vtop{\hsize=\tablesize
\hbox to\hsize{\kern 12em 
$\xi_i=4\,\beta_i+1\leqslant 6\,\beta_i-2\,\omega_i$
\hss\vrule height 13pt depth 5pt\ \,$\surd$\ }
}
\vrule height 13pt depth 5pt\vskip -1pt
\noindent
\vrule height 13pt depth 5pt
\vtop{\hsize=\tablesize
\hbox to\hsize{\kern 12em 
$\xi_i=6\,\beta_i-2\,\omega_i\leqslant 4\,\beta_i+1$
\hss\vrule height 13pt depth 5pt\ \,$\phantom{\surd}$\ }
}
\vrule height 13pt depth 5pt
\hrule width\hsize
}
\vskip 1ex plus 3pt
     In the following cases the multiplicity $\omega_i$ coincides with 
the multiplicity $\beta_i$. Then, according to \mythetag{3.3}, we have
$\theta_i=\max(\alpha_i,\beta_i,\omega_i)=\beta_i$.
\vskip 1ex plus 3pt
\noindent\mytable{5.57} (symmetric to the table \mythetable{5.23})\par
\vbox{
\hrule width\hsize
\noindent
\vrule height 13pt depth 5pt
\vtop{\hsize=\tablesize
\hbox to\hsize{\ \ $\beta_i=\omega_i>\alpha_i=0$, $p_i\neq 2,3$\ \ 
\vrule height 13pt depth 5pt\ \ $\mu_i>0\ \Rightarrow
\ \mu_i=\beta_i$ and $\eta_i=0$
\hss}
}
\vrule height 13pt depth 5pt
\hrule width\hsize
\noindent
\vrule height 13pt depth 5pt
\vtop{\hsize=\tablesize
\hbox to\hsize{\ \ \ \ \ $m_C=2\,\beta_i+2\,\mu_i$,\kern 4.3em
$m_D=\xi_i\geqslant 4\,\beta_i+2\,\eta_i$,\kern 4.3em 
$m_E=\varkappa_i\geqslant 4\,\beta_i$
\hss}
}
\vrule height 13pt depth 5pt
\hrule width\hsize
\noindent
\vrule height 13pt depth 5pt
\vtop{\hsize=\tablesize
\hbox to\hsize{\kern 13em 
$2\,\beta_i+2\,\mu_i=\varkappa_i\leqslant\xi_i$
\hss\vrule height 13pt depth 5pt\ \,$\surd$\ }
}
\vrule height 13pt depth 5pt\vskip -1pt
\noindent
\vrule height 13pt depth 5pt
\vtop{\hsize=\tablesize
\hbox to\hsize{\kern 13em 
$\xi_i=\varkappa_i\leqslant 2\,\beta_i+2\,\mu_i$
\hss\vrule height 13pt depth 5pt\ \,$\surd$\ }
}
\vrule height 13pt depth 5pt\vskip -1pt
\noindent
\vrule height 13pt depth 5pt
\vtop{\hsize=\tablesize
\hbox to\hsize{\kern 13em 
$\xi_i=2\,\beta_i+2\,\mu_i\leqslant\varkappa_i$
\hss\vrule height 13pt depth 5pt\ \,$\surd$\ }
}
\vrule height 13pt depth 5pt
\hrule width\hsize
}
\vskip 1ex plus 3pt
\noindent\mytable{5.58} (symmetric to the table \mythetable{5.24})\par
\vbox{
\hrule width\hsize
\noindent
\vrule height 13pt depth 5pt
\vtop{\hsize=\tablesize
\hbox to\hsize{\ \ $\beta_i=\omega_i>\alpha_i=0$, $p_i\neq 2,3$\ \ 
\vrule height 13pt depth 5pt\ \ $\mu_i=0$ and $\eta_i=\beta_i$
\hss}
}
\vrule height 13pt depth 5pt
\hrule width\hsize
\noindent
\vrule height 13pt depth 5pt
\vtop{\hsize=\tablesize
\hbox to\hsize{\ \ \ \ \ $m_C=\zeta_i\geqslant 2\,\beta_i$,
\kern 5.8em
$m_D=\xi_i\geqslant 6\,\beta_i$,\kern 5.8em 
$m_E=\varkappa_i\geqslant 4\,\beta_i$
\hss}
}
\vrule height 13pt depth 5pt
\hrule width\hsize
\noindent
\vrule height 13pt depth 5pt
\vtop{\hsize=\tablesize
\hbox to\hsize{\kern 15em 
$\zeta_i=\varkappa_i\leqslant\xi_i$
\hss\vrule height 13pt depth 5pt\ \,$\surd$\ }
}
\vrule height 13pt depth 5pt\vskip -1pt
\noindent
\vrule height 13pt depth 5pt
\vtop{\hsize=\tablesize
\hbox to\hsize{\kern 15em 
$\xi_i=\varkappa_i\leqslant\zeta_i$
\hss\vrule height 13pt depth 5pt\ \,$\surd$\ }
}
\vrule height 13pt depth 5pt\vskip -1pt
\noindent
\vrule height 13pt depth 5pt
\vtop{\hsize=\tablesize
\hbox to\hsize{\kern 15em 
$\xi_i=\zeta_i\leqslant\varkappa_i$
\hss\vrule height 13pt depth 5pt\ \,$\surd$\ }
}
\vrule height 13pt depth 5pt
\hrule width\hsize
}
\vskip 1ex plus 3pt
\noindent\mytable{5.59} (symmetric to the table \mythetable{5.25})\par
\vbox{
\hrule width\hsize
\noindent
\vrule height 13pt depth 5pt
\vtop{\hsize=\tablesize
\hbox to\hsize{\ \ $\beta_i=\omega_i>\alpha_i=0$, $p_i=2$\ \ 
\vrule height 13pt depth 5pt\ \ $\mu_i>0\ \Rightarrow
\ \mu_i=\beta_i$ and $\eta_i=0$
\hss}
}
\vrule height 13pt depth 5pt
\hrule width\hsize
\noindent
\vrule height 13pt depth 5pt
\vtop{\hsize=\tablesize
\hbox to\hsize{\ \ \ \ \ $m_C=2\,\beta_i+2\,\mu_i$,\kern 3.4em
$m_D=\xi_i\geqslant 4\,\beta_i+2\,\eta_i$,\kern 3.4em 
$m_E=\varkappa_i\geqslant 4\,\beta_i+1$
\hss}
}
\vrule height 13pt depth 5pt
\hrule width\hsize
\noindent
\vrule height 13pt depth 5pt
\vtop{\hsize=\tablesize
\hbox to\hsize{\kern 13em 
$2\,\beta_i+2\,\mu_i=\varkappa_i\leqslant\xi_i$
\hss\vrule height 13pt depth 5pt\ \,$\phantom{\surd}$\ }
}
\vrule height 13pt depth 5pt\vskip -1pt
\noindent
\vrule height 13pt depth 5pt
\vtop{\hsize=\tablesize
\hbox to\hsize{\kern 13em 
$\xi_i=\varkappa_i\leqslant 2\,\beta_i+2\,\mu_i$
\hss\vrule height 13pt depth 5pt\ \,$\phantom{\surd}$\ }
}
\vrule height 13pt depth 5pt\vskip -1pt
\noindent
\vrule height 13pt depth 5pt
\vtop{\hsize=\tablesize
\hbox to\hsize{\kern 13em 
$\xi_i=2\,\beta_i+2\,\mu_i\leqslant\varkappa_i$
\hss\vrule height 13pt depth 5pt\ \,$\surd$\ }
}
\vrule height 13pt depth 5pt
\hrule width\hsize
}\pagebreak
\noindent\mytable{5.60} (symmetric to the table \mythetable{5.26})\par
\vbox{
\hrule width\hsize
\noindent
\vrule height 13pt depth 5pt
\vtop{\hsize=\tablesize
\hbox to\hsize{\ \ $\beta_i=\omega_i>\alpha_i=0$, $p_i=2$\ \ 
\vrule height 13pt depth 5pt\ \ $\mu_i=0$ and $\eta_i=\beta_i$
\hss}
}
\vrule height 13pt depth 5pt
\hrule width\hsize
\noindent
\vrule height 13pt depth 5pt
\vtop{\hsize=\tablesize
\hbox to\hsize{\ \ \ \ \ $m_C=\zeta_i\geqslant 2\,\beta_i$,
\kern 5em
$m_D=\xi_i\geqslant 6\,\beta_i$,\kern 5em 
$m_E=\varkappa_i\geqslant 4\,\beta_i+1$
\hss}
}
\vrule height 13pt depth 5pt
\hrule width\hsize
\noindent
\vrule height 13pt depth 5pt
\vtop{\hsize=\tablesize
\hbox to\hsize{\kern 15em 
$\zeta_i=\varkappa_i\leqslant\xi_i$
\hss\vrule height 13pt depth 5pt\ \,$\surd$\ }
}
\vrule height 13pt depth 5pt\vskip -1pt
\noindent
\vrule height 13pt depth 5pt
\vtop{\hsize=\tablesize
\hbox to\hsize{\kern 15em 
$\xi_i=\varkappa_i\leqslant\zeta_i$
\hss\vrule height 13pt depth 5pt\ \,$\surd$\ }
}
\vrule height 13pt depth 5pt\vskip -1pt
\noindent
\vrule height 13pt depth 5pt
\vtop{\hsize=\tablesize
\hbox to\hsize{\kern 15em 
$\xi_i=\zeta_i\leqslant\varkappa_i$
\hss\vrule height 13pt depth 5pt\ \,$\surd$\ }
}
\vrule height 13pt depth 5pt
\hrule width\hsize
}
\vskip 1ex plus 3pt
\noindent\mytable{5.61} (symmetric to the table \mythetable{5.27})\par
\vbox{
\hrule width\hsize
\noindent
\vrule height 13pt depth 5pt
\vtop{\hsize=\tablesize
\hbox to\hsize{\ \ $\beta_i=\omega_i>\alpha_i=0$, $p_i=3$\ \ 
\vrule height 13pt depth 5pt\ \ $\mu_i>0\ \Rightarrow
\ \mu_i=\beta_i$ and $\eta_i=0$
\hss}
}
\vrule height 13pt depth 5pt
\hrule width\hsize
\noindent
\vrule height 13pt depth 5pt
\vtop{\hsize=\tablesize
\hbox to\hsize{\ \ \ \ \ $m_C=2\,\beta_i+2\,\mu_i$,\kern 5.4em
$m_D=\xi_i\geqslant 4\,\beta_i+2\,\eta_i$,\kern 5.4em 
$m_E=4\,\beta_i$
\hss}
}
\vrule height 13pt depth 5pt
\hrule width\hsize
\noindent
\vrule height 13pt depth 5pt
\vtop{\hsize=\tablesize
\hbox to\hsize{\kern 13em 
$2\,\beta_i+2\,\mu_i=4\,\beta_i\leqslant\xi_i$
\hss\vrule height 13pt depth 5pt\ \,$\surd$\ }
}
\vrule height 13pt depth 5pt\vskip -1pt
\noindent
\vrule height 13pt depth 5pt
\vtop{\hsize=\tablesize
\hbox to\hsize{\kern 13em 
$\xi_i=4\,\beta_i\leqslant 2\,\beta_i+2\,\mu_i$
\hss\vrule height 13pt depth 5pt\ \,$\surd$\ }
}
\vrule height 13pt depth 5pt\vskip -1pt
\noindent
\vrule height 13pt depth 5pt
\vtop{\hsize=\tablesize
\hbox to\hsize{\kern 13em 
$\xi_i=2\,\beta_i+2\,\mu_i\leqslant 4\,\beta_i$
\hss\vrule height 13pt depth 5pt\ \,$\surd$\ }
}
\vrule height 13pt depth 5pt
\hrule width\hsize
}
\vskip 1ex plus 3pt
\noindent\mytable{5.62} (symmetric to the table \mythetable{5.28})\par
\vbox{
\hrule width\hsize
\noindent
\vrule height 13pt depth 5pt
\vtop{\hsize=\tablesize
\hbox to\hsize{\ \ $\beta_i=\omega_i>\alpha_i=0$, $p_i=3$\ \ 
\vrule height 13pt depth 5pt\ \ $\mu_i=0$ and $\eta_i=\beta_i$
\hss}
}
\vrule height 13pt depth 5pt
\hrule width\hsize
\noindent
\vrule height 13pt depth 5pt
\vtop{\hsize=\tablesize
\hbox to\hsize{\ \ \ \ \ $m_C=\zeta_i\geqslant 2\,\beta_i$,
\kern 7em
$m_D=\xi_i\geqslant 6\,\beta_i$,\kern 7em 
$m_E=4\,\beta_i$
\hss}
}
\vrule height 13pt depth 5pt
\hrule width\hsize
\noindent
\vrule height 13pt depth 5pt
\vtop{\hsize=\tablesize
\hbox to\hsize{\kern 15em 
$\zeta_i=4\,\beta_i\leqslant\xi_i$
\hss\vrule height 13pt depth 5pt\ \,$\surd$\ }
}
\vrule height 13pt depth 5pt\vskip -1pt
\noindent
\vrule height 13pt depth 5pt
\vtop{\hsize=\tablesize
\hbox to\hsize{\kern 15em 
$\xi_i=4\,\beta_i\leqslant\zeta_i$
\hss\vrule height 13pt depth 5pt\ \,$\phantom{\surd}$\ }
}
\vrule height 13pt depth 5pt\vskip -1pt
\noindent
\vrule height 13pt depth 5pt
\vtop{\hsize=\tablesize
\hbox to\hsize{\kern 15em 
$\xi_i=\zeta_i\leqslant 4\,\beta_i$
\hss\vrule height 13pt depth 5pt\ \,$\phantom{\surd}$\ }
}
\vrule height 13pt depth 5pt
\hrule width\hsize
}
\vskip 1ex plus 3pt
     In the following cases the multiplicity $\omega_i$ is greater than
the multiplicity $\beta_i$. Therefore, according to \mythetag{3.3}, we
have $\theta_i=\max(\alpha_i,\beta_i,\omega_i)=\omega_i$.
\vskip 1ex plus 3pt
\noindent\mytable{5.63} (symmetric to the table \mythetable{5.29})\par
\vbox{
\hrule width\hsize
\noindent
\vrule height 13pt depth 5pt
\vtop{\hsize=\tablesize
\hbox to\hsize{\ \ $\omega_i>\beta_i>\alpha_i=0$, $p_i\neq 2$\ \ 
\vrule height 13pt depth 5pt\ \ $\mu_i>0$ and $\eta_i>0\ \Rightarrow
\ \mu_i=\beta_i$ and $\eta_i=\omega_i-\beta_i$
\hss}
}
\vrule height 13pt depth 5pt
\hrule width\hsize
\noindent
\vrule height 13pt depth 5pt
\vtop{\hsize=\tablesize
\hbox to\hsize{\ \ \ \ \ $m_C=4\,\omega_i-2\,\beta_i+2\,\mu_i$,\kern 3.6em
$m_D=2\,\omega_i+2\,\beta_i+2\,\eta_i$,\kern 3.6em 
$m_E=4\,\omega_i$
\hss}
}
\vrule height 13pt depth 5pt
\hrule width\hsize
\noindent
\vrule height 13pt depth 5pt
\vtop{\hsize=\tablesize
\hbox to\hsize{\kern 8em 
$4\,\omega_i-2\,\beta_i+2\,\mu_i=4\,\omega_i\leqslant 
2\,\omega_i+2\,\beta_i+2\,\eta_i$
\hss\vrule height 13pt depth 5pt\ \,$\surd$\ }
}
\vrule height 13pt depth 5pt\vskip -1pt
\noindent
\vrule height 13pt depth 5pt
\vtop{\hsize=\tablesize
\hbox to\hsize{\kern 8em 
$2\,\omega_i+2\,\beta_i+2\,\eta_i=4\,\omega_i\leqslant 4\,\omega_i
-2\,\beta_i+2\,\mu_i$
\hss\vrule height 13pt depth 5pt\ \,$\surd$\ }
}
\vrule height 13pt depth 5pt\vskip -1pt
\noindent
\vrule height 13pt depth 5pt
\vtop{\hsize=\tablesize
\hbox to\hsize{\kern 8em 
$2\,\omega_i+2\,\beta_i+2\,\eta_i=4\,\omega_i-2\,\beta_i+2\,\mu_i
\leqslant 4\,\omega_i$
\hss\vrule height 13pt depth 5pt\ \,$\surd$\ }
}
\vrule height 13pt depth 5pt
\hrule width\hsize
}
\vskip 1ex plus 3pt
\noindent\mytable{5.64} (symmetric to the table \mythetable{5.30})\par
\vbox{
\hrule width\hsize
\noindent
\vrule height 13pt depth 5pt
\vtop{\hsize=\tablesize
\hbox to\hsize{\ \ $\omega_i>\beta_i>\alpha_i=0$, $p_i\neq 2$\ \ 
\vrule height 13pt depth 5pt\ \ $\mu_i=0$ and $\eta_i=\omega_i$
\hss}
}
\vrule height 13pt depth 5pt
\hrule width\hsize
\noindent
\vrule height 13pt depth 5pt
\vtop{\hsize=\tablesize
\hbox to\hsize{\ \ \ \ \ $m_C=\zeta_i\geqslant 4\,\omega_i-2\,\beta_i$,
\kern 5.25em
$m_D=4\,\omega_i+2\,\beta_i$,\kern 5.25em
$m_E=4\,\omega_i$
\hss}
}
\vrule height 13pt depth 5pt
\hrule width\hsize
\noindent
\vrule height 13pt depth 5pt
\vtop{\hsize=\tablesize
\hbox to\hsize{\kern 13em 
$\zeta_i=4\,\omega_i\leqslant 4\,\omega_i+2\,\beta_i$
\hss\vrule height 13pt depth 5pt\ \,$\surd$\ }
}
\vrule height 13pt depth 5pt\vskip -1pt
\noindent
\vrule height 13pt depth 5pt
\vtop{\hsize=\tablesize
\hbox to\hsize{\kern 13em 
$4\,\omega_i+2\,\beta_i=4\,\omega_i\leqslant\zeta_i$
\hss\vrule height 13pt depth 5pt\ \,$\phantom{\surd}$\ }
}
\vrule height 13pt depth 5pt\vskip -1pt
\noindent
\vrule height 13pt depth 5pt
\vtop{\hsize=\tablesize
\hbox to\hsize{\kern 13em 
$4\,\omega_i+2\,\beta_i=\zeta_i\leqslant 4\,\omega_i$
\hss\vrule height 13pt depth 5pt\ \,$\phantom{\surd}$\ }
}
\vrule height 13pt depth 5pt
\hrule width\hsize
}\pagebreak
\noindent\mytable{5.65} (symmetric to the table \mythetable{5.31})\par
\vbox{
\hrule width\hsize
\noindent
\vrule height 13pt depth 5pt
\vtop{\hsize=\tablesize
\hbox to\hsize{\ \ $\omega_i>\beta_i>\alpha_i=0$, $p_i\neq 2$\ \ 
\vrule height 13pt depth 5pt\ \ $\mu_i=\omega_i$ and $\eta_i=0$
\hss}
}
\vrule height 13pt depth 5pt
\hrule width\hsize
\noindent
\vrule height 13pt depth 5pt
\vtop{\hsize=\tablesize
\hbox to\hsize{\ \ \ \ \ $m_C=6\,\omega_i-2\,\beta_i$,
\kern 5.25em
$m_D=\xi_i\geqslant 2\,\omega_i+2\,\beta_i$,\kern 5.25em
$m_E=4\,\omega_i$
\hss}
}
\vrule height 13pt depth 5pt
\hrule width\hsize
\noindent
\vrule height 13pt depth 5pt
\vtop{\hsize=\tablesize
\hbox to\hsize{\kern 13em 
$6\,\omega_i-2\,\beta_i=4\,\omega_i\leqslant\xi_i$
\hss\vrule height 13pt depth 5pt\ \,$\phantom{\surd}$\ }
}
\vrule height 13pt depth 5pt\vskip -1pt
\noindent
\vrule height 13pt depth 5pt
\vtop{\hsize=\tablesize
\hbox to\hsize{\kern 13em 
$\xi_i=4\,\omega_i\leqslant 6\,\omega_i-2\,\beta_i$
\hss\vrule height 13pt depth 5pt\ \,$\surd$\ }
}
\vrule height 13pt depth 5pt\vskip -1pt
\noindent
\vrule height 13pt depth 5pt
\vtop{\hsize=\tablesize
\hbox to\hsize{\kern 13em 
$\xi_i=6\,\omega_i-2\,\beta_i\leqslant 4\,\omega_i$
\hss\vrule height 13pt depth 5pt\ \,$\phantom{\surd}$\ }
}
\vrule height 13pt depth 5pt
\hrule width\hsize
}
\vskip 1ex plus 3pt
\noindent\mytable{5.66} (symmetric to the table \mythetable{5.32})\par
\vbox{
\hrule width\hsize
\noindent
\vrule height 13pt depth 5pt
\vtop{\hsize=\tablesize
\hbox to\hsize{\ \ $\omega_i>\beta_i>\alpha_i=0$, $p_i=2$\ \ 
\vrule height 13pt depth 5pt\ \ $\mu_i>0$ and $\eta_i>0\ \Rightarrow
\ \mu_i=\beta_i$ and $\eta_i=\omega_i-\beta_i$
\hss}
}
\vrule height 13pt depth 5pt
\hrule width\hsize
\noindent
\vrule height 13pt depth 5pt
\vtop{\hsize=\tablesize
\hbox to\hsize{\ \ \ \ \ $m_C=4\,\omega_i-2\,\beta_i+2\,\mu_i$,\kern 2.8em
$m_D=2\,\omega_i+2\,\beta_i+2\,\eta_i$,\kern 2.8em 
$m_E=4\,\omega_i+1$
\hss}
}
\vrule height 13pt depth 5pt
\hrule width\hsize
\noindent
\vrule height 13pt depth 5pt
\vtop{\hsize=\tablesize
\hbox to\hsize{\kern 8em 
$4\,\omega_i-2\,\beta_i+2\,\mu_i=4\,\omega_i+1\leqslant 
2\,\omega_i+2\,\beta_i+2\,\eta_i$
\hss\vrule height 13pt depth 5pt\ \,$\phantom{\surd}$\ }
}
\vrule height 13pt depth 5pt\vskip -1pt
\noindent
\vrule height 13pt depth 5pt
\vtop{\hsize=\tablesize
\hbox to\hsize{\kern 8em 
$2\,\omega_i+2\,\beta_i+2\,\eta_i=4\,\omega_i+1\leqslant 4\,\omega_i
-2\,\beta_i+2\,\mu_i$
\hss\vrule height 13pt depth 5pt\ \,$\phantom{\surd}$\ }
}
\vrule height 13pt depth 5pt\vskip -1pt
\noindent
\vrule height 13pt depth 5pt
\vtop{\hsize=\tablesize
\hbox to\hsize{\kern 8em 
$2\,\omega_i+2\,\beta_i+2\,\eta_i=4\,\omega_i-2\,\beta_i+2\,\mu_i
\leqslant 4\,\omega_i+1$
\hss\vrule height 13pt depth 5pt\ \,$\surd$\ }
}
\vrule height 13pt depth 5pt
\hrule width\hsize
}
\vskip 1ex plus 3pt
\noindent\mytable{5.67} (symmetric to the table \mythetable{5.33})\par
\vbox{
\hrule width\hsize
\noindent
\vrule height 13pt depth 5pt
\vtop{\hsize=\tablesize
\hbox to\hsize{\ \ $\omega_i>\beta_i>\alpha_i=0$, $p_i=2$\ \ 
\vrule height 13pt depth 5pt\ \ $\mu_i=0$ and $\eta_i=\omega_i$
\hss}
}
\vrule height 13pt depth 5pt
\hrule width\hsize
\noindent
\vrule height 13pt depth 5pt
\vtop{\hsize=\tablesize
\hbox to\hsize{\ \ \ \ \ $m_C=\zeta_i\geqslant 4\,\omega_i-2\,\beta_i$,
\kern 4.5em
$m_D=4\,\omega_i+2\,\beta_i$,\kern 4.5em
$m_E=4\,\omega_i+1$
\hss}
}
\vrule height 13pt depth 5pt
\hrule width\hsize
\noindent
\vrule height 13pt depth 5pt
\vtop{\hsize=\tablesize
\hbox to\hsize{\kern 12em 
$\zeta_i=4\,\omega_i+1\leqslant 4\,\omega_i+2\,\beta_i$
\hss\vrule height 13pt depth 5pt\ \,$\surd$\ }
}
\vrule height 13pt depth 5pt\vskip -1pt
\noindent
\vrule height 13pt depth 5pt
\vtop{\hsize=\tablesize
\hbox to\hsize{\kern 12em 
$4\,\omega_i+2\,\beta_i=4\,\omega_i+1\leqslant\zeta_i$
\hss\vrule height 13pt depth 5pt\ \,$\phantom{\surd}$\ }
}
\vrule height 13pt depth 5pt\vskip -1pt
\noindent
\vrule height 13pt depth 5pt
\vtop{\hsize=\tablesize
\hbox to\hsize{\kern 12em 
$4\,\omega_i+2\,\beta_i=\zeta_i\leqslant 4\,\omega_i+1$
\hss\vrule height 13pt depth 5pt\ \,$\phantom{\surd}$\ }
}
\vrule height 13pt depth 5pt
\hrule width\hsize
}
\vskip 1ex plus 3pt
\noindent\mytable{5.68} (symmetric to the table \mythetable{5.34})\par
\vbox{
\hrule width\hsize
\noindent
\vrule height 13pt depth 5pt
\vtop{\hsize=\tablesize
\hbox to\hsize{\ \ $\omega_i>\beta_i>\alpha_i=0$, $p_i=2$\ \ 
\vrule height 13pt depth 5pt\ \ $\mu_i=\omega_i$ and $\eta_i=0$
\hss}
}
\vrule height 13pt depth 5pt
\hrule width\hsize
\noindent
\vrule height 13pt depth 5pt
\vtop{\hsize=\tablesize
\hbox to\hsize{\ \ \ \ \ $m_C=6\,\omega_i-2\,\beta_i$,
\kern 4.3em
$m_D=\xi_i\geqslant 2\,\omega_i+2\,\beta_i$,\kern 4.3em
$m_E=4\,\omega_i+1$
\hss}
}
\vrule height 13pt depth 5pt
\hrule width\hsize
\noindent
\vrule height 13pt depth 5pt
\vtop{\hsize=\tablesize
\hbox to\hsize{\kern 13em 
$6\,\omega_i-2\,\beta_i=4\,\omega_i+1\leqslant\xi_i$
\hss\vrule height 13pt depth 5pt\ \,$\phantom{\surd}$\ }
}
\vrule height 13pt depth 5pt\vskip -1pt
\noindent
\vrule height 13pt depth 5pt
\vtop{\hsize=\tablesize
\hbox to\hsize{\kern 13em 
$\xi_i=4\,\omega_i+1\leqslant 6\,\omega_i-2\,\beta_i$
\hss\vrule height 13pt depth 5pt\ \,$\surd$\ }
}
\vrule height 13pt depth 5pt\vskip -1pt
\noindent
\vrule height 13pt depth 5pt
\vtop{\hsize=\tablesize
\hbox to\hsize{\kern 13em 
$\xi_i=6\,\omega_i-2\,\beta_i\leqslant 4\,\omega_i+1$
\hss\vrule height 13pt depth 5pt\ \,$\phantom{\surd}$\ }
}
\vrule height 13pt depth 5pt
\hrule width\hsize
}
\vskip 1ex plus 3pt
     Thus, totally we have 68 cases placed into 68 tables. They describe
completely the structure of the prime factors $p_1,\,\ldots,\,p_n$ in 
\mythetag{5.1}.\par
\head
6. The structural theorem. 
\endhead
     Analyzing the whole variety of data in the tables~\mythetable{5.1} through 
\mythetable{5.68}, we subdivide the set of prime factors $S=\{p_1,\,\ldots,\,p_n\}$ 
from \mythetag{5.1} into a disjoint union of several sets. The first three of such 
sets are given by the formulas
$$
\gather
\hskip -2em
S_1=\{p_i\in S:\ \alpha_i>\beta_i>\omega_i=0,\ \mu_i=\beta_i,
\ \eta_i=\alpha_i-\beta_i\},\\
\hskip -2em
S_2=\{p_i\in S:\ \alpha_i>\beta_i>\omega_i=0,\ \mu_i=0,
\ \eta_i=\alpha_i\},
\mytag{6.1}\\
\hskip -2em
S_3=\{p_i\in S:\ \alpha_i>\beta_i>\omega_i=0,\ \mu_i=\alpha_i,
\ \eta_i=0\}.
\endgather
$$
The formulas \mythetag{6.1} correspond to the tables~\mythetable{5.1} through 
\mythetable{5.6}. Using the formulas \mythetag{6.1}, we define the following 
integer numbers:
$$
\xalignat 2
&\hskip -2em
b_1=\prod_{p_i\in S_1}p_i^{\,\beta_i},
&&\tilde b_1=\prod_{p_i\in S_1}p_i^{\,\alpha_i-\beta_i},\\
&\hskip -2em
b_2=\prod_{p_i\in S_2}p_i^{\,\beta_i},
&&\tilde b_2=\prod_{p_i\in S_2}p_i^{\,\alpha_i-\beta_i},
\mytag{6.2}\\
&\hskip -2em
b_3=\prod_{p_i\in S_3}p_i^{\,\beta_i},
&&\tilde b_3=\prod_{p_i\in S_3}p_i^{\,\alpha_i-\beta_i}.
\endxalignat
$$
\par
     The next three sets $S_4$, $S_5$, and $S_6$ are defined by the following formulas:
$$
\gather
\hskip -2em
S_4=\{p_i\in S:\ \omega_i>\beta_i>\alpha_i=0,\ \mu_i=\beta_i,
\ \eta_i=\omega_i-\beta_i\},\\
\hskip -2em
S_5=\{p_i\in S:\ \omega_i>\beta_i>\alpha_i=0,\ \mu_i=0,
\ \eta_i=\omega_i\},
\mytag{6.3}\\
\hskip -2em
S_6=\{p_i\in S:\ \omega_i>\beta_i>\alpha_i=0,\ \mu_i=\omega_i,
\ \eta_i=0\}.
\endgather
$$
The formulas \mythetag{6.3} correspond to the tables~\mythetable{5.63} through
\mythetable{5.68}.  Using the formulas \mythetag{6.3}, we define the following 
integer numbers:
$$
\xalignat 2
&\hskip -2em
b_4=\prod_{p_i\in S_4}p_i^{\,\beta_i},
&&\tilde b_4=\prod_{p_i\in S_4}p_i^{\,\omega_i-\beta_i},\\
&\hskip -2em
b_5=\prod_{p_i\in S_5}p_i^{\,\beta_i},
&&\tilde b_5=\prod_{p_i\in S_5}p_i^{\,\omega_i-\beta_i},
\mytag{6.4}\\
&\hskip -2em
b_6=\prod_{p_i\in S_6}p_i^{\,\beta_i},
&&\tilde b_6=\prod_{p_i\in S_6}p_i^{\,\omega_i-\beta_i}.
\endxalignat
$$
Some of the sets $S_1$, $S_2$, $S_3$, $S_4$, $S_5$, $S_6$ can be empty. Therefore 
we interpret the formulas \mythetag{6.2} and \mythetag{6.4} so that $b_k=\tilde 
b_k=1$ if the corresponding set $S_k$ is empty. Moreover, $b_k\neq 1$ implies 
$\tilde b_k\neq 1$ and vice versa. If $b_k\cdot\tilde b_k\neq 1$, then the prime 
factors of the number $b_k$ coincide with the prime factors of the number 
$\tilde b_k$.\par
     The next three sets $S_7$, $S_8$, and $S_9$ are defined by the 
following formulas:
$$
\gather
\hskip -2em
S_7=\{p_i\in S:\ \beta_i>\omega_i>\alpha_i=0,\ \mu_i=\omega_i,
\ \eta_i=\beta_i-\omega_i\},\\
\hskip -2em
S_8=\{p_i\in S:\ \beta_i>\omega_i>\alpha_i=0,\ \mu_i=0,
\ \eta_i=\beta_i\},
\mytag{6.5}\\
\hskip -2em
S_9=\{p_i\in S:\ \beta_i>\omega_i>\alpha_i=0,\ \mu_i=\beta_i,
\ \eta_i=0\}.
\endgather
$$
The formulas \mythetag{6.5} correspond to the tables~\mythetable{5.51} through
\mythetable{5.56}. Using the formulas \mythetag{6.5}, we define the following 
integer numbers:
$$
\xalignat 2
&\hskip -2em
u_1=\prod_{p_i\in S_7}p_i^{\,\omega_i},
&&\tilde u_1=\prod_{p_i\in S_7}p_i^{\,\beta_i-\omega_i},\\
&\hskip -2em
u_2=\prod_{p_i\in S_8}p_i^{\,\omega_i},
&&\tilde u_2=\prod_{p_i\in S_8}p_i^{\,\beta_i-\omega_i},
\mytag{6.6}\\
&\hskip -2em
u_3=\prod_{p_i\in S_9}p_i^{\,\omega_i},
&&\tilde u_3=\prod_{p_i\in S_9}p_i^{\,\beta_i-\omega_i}.
\endxalignat
$$
\par
     The next three sets $S_{10}$, $S_{11}$, and $S_{12}$ are defined by the 
following formulas:
$$
\gather
\hskip -2em
S_{10}=\{p_i\in S:\ \alpha_i>\omega_i>\beta_i=0,\ \mu_i=\omega_i,
\ \eta_i=\alpha_i-\omega_i\},\\
\hskip -2em
S_{11}=\{p_i\in S:\ \alpha_i>\omega_i>\beta_i=0,\ \mu_i=0,
\ \eta_i=\alpha_i\},
\mytag{6.7}\\
\hskip -2em
S_{12}=\{p_i\in S:\ \alpha_i>\omega_i>\beta_i=0,\ \mu_i=\alpha_i,
\ \eta_i=0\}.
\endgather
$$
The formulas \mythetag{6.7} correspond to the tables~\mythetable{5.17} through
\mythetable{5.22}. Using the formulas \mythetag{6.7}, we define the following 
integer numbers:
$$
\xalignat 2
&\hskip -2em
u_4=\prod_{p_i\in S_{10}}p_i^{\,\omega_i},
&&\tilde u_4=\prod_{p_i\in S_{10}}p_i^{\,\alpha_i-\omega_i},\\
&\hskip -2em
u_5=\prod_{p_i\in S_{11}}p_i^{\,\omega_i},
&&\tilde u_5=\prod_{p_i\in S_{11}}p_i^{\,\alpha_i-\omega_i},
\mytag{6.8}\\
&\hskip -2em
u_6=\prod_{p_i\in S_{12}}p_i^{\,\omega_i},
&&\tilde u_6=\prod_{p_i\in S_{12}}p_i^{\,\alpha_i-\omega_i}.
\endxalignat
$$
Some of the sets $S_7$, $S_8$, $S_9$, $S_{10}$, $S_{11}$, $S_{12}$ can be empty. 
Therefore we interpret the formulas \mythetag{6.6} and \mythetag{6.8} so that 
$u_k=\tilde u_k=1$ if the corresponding set $S_k$ is empty. Moreover, $u_k\neq 1$ 
implies $\tilde u_k\neq 1$ and vice versa. If $u_k\cdot\tilde u_k\neq 1$, then 
the prime factors of the number $u_k$ coincide with the prime factors of the number 
$\tilde u_k$.\par
     The next three sets $S_{13}$, $S_{14}$, and $S_{15}$ are defined by the 
following formulas:
$$
\gather
\hskip -2em
S_{13}=\{p_i\in S:\ \omega_i>\alpha_i>\beta_i=0,\ \mu_i=\alpha_i,
\ \eta_i=\omega_i-\alpha_i\},\\
\hskip -2em
S_{14}=\{p_i\in S:\ \omega_i>\alpha_i>\beta_i=0,\ \mu_i=0,
\ \eta_i=\omega_i\},
\mytag{6.9}\\
\hskip -2em
S_{15}=\{p_i\in S:\ \omega_i>\alpha_i>\beta_i=0,\ \mu_i=\omega_i,
\ \eta_i=0\}.
\endgather
$$
The formulas \mythetag{6.9} correspond to the tables~\mythetable{5.29} through
\mythetable{5.34}. Using the formulas \mythetag{6.9}, we define the following 
integer numbers:
$$
\xalignat 2
&\hskip -2em
a_1=\prod_{p_i\in S_{13}}p_i^{\,\alpha_i},
&&\tilde a_1=\prod_{p_i\in S_{13}}p_i^{\,\omega_i-\alpha_i},\\
&\hskip -2em
a_2=\prod_{p_i\in S_{14}}p_i^{\,\alpha_i},
&&\tilde a_2=\prod_{p_i\in S_{14}}p_i^{\,\omega_i-\alpha_i},
\mytag{6.10}\\
&\hskip -2em
a_3=\prod_{p_i\in S_{15}}p_i^{\,\alpha_i},
&&\tilde a_3=\prod_{p_i\in S_{15}}p_i^{\,\omega_i-\alpha_i}.
\endxalignat
$$
\par
     The next three sets $S_{16}$, $S_{17}$, and $S_{18}$ are defined similarly.
For this purpose we use the following three formulas analogous to \mythetag{6.9}:
$$
\gather
\hskip -2em
S_{16}=\{p_i\in S:\ \beta_i>\alpha_i>\omega_i=0,\ \mu_i=\alpha_i,
\ \eta_i=\beta_i-\alpha_i\},\\
\hskip -2em
S_{17}=\{p_i\in S:\ \beta_i>\alpha_i>\omega_i=0,\ \mu_i=0,
\ \eta_i=\beta_i\},
\mytag{6.11}\\
\hskip -2em
S_{18}=\{p_i\in S:\ \beta_i>\alpha_i>\omega_i=0,\ \mu_i=\beta_i,
\ \eta_i=0\}.
\endgather
$$
The formulas \mythetag{6.11} correspond to the tables~\mythetable{5.39} through
\mythetable{5.44}. We use the formulas \mythetag{6.11} in order to define six
numbers $a_4$, $a_5$, $a_6$, $\tilde a_4$, $\tilde a_5$, $\tilde a_6$ similar
to the numbers $a_1$, $a_2$, $a_3$, $\tilde a_1$, $\tilde a_2$, $\tilde a_3$
in \mythetag{6.10}. For this purpose we write 
$$
\xalignat 2
&\hskip -2em
a_4=\prod_{p_i\in S_{16}}p_i^{\,\alpha_i},
&&\tilde a_4=\prod_{p_i\in S_{16}}p_i^{\,\beta_i-\alpha_i},\\
&\hskip -2em
a_5=\prod_{p_i\in S_{17}}p_i^{\,\alpha_i},
&&\tilde a_5=\prod_{p_i\in S_{17}}p_i^{\,\beta_i-\alpha_i},
\mytag{6.12}\\
&\hskip -2em
a_6=\prod_{p_i\in S_{18}}p_i^{\,\alpha_i},
&&\tilde a_6=\prod_{p_i\in S_{18}}p_i^{\,\beta_i-\alpha_i}.
\endxalignat
$$
Some of the sets $S_{13}$, $S_{14}$, $S_{15}$, $S_{16}$, $S_{17}$, $S_{18}$ can 
be empty. Therefore we interpret the formulas \mythetag{6.10} and \mythetag{6.12} 
so that $a_k=\tilde a_k=1$ if the corresponding set $S_k$ is empty. Moreover, 
$a_k\neq 1$ implies $\tilde a_k\neq 1$ and vice versa. If $a_k\cdot\tilde a_k\neq 1$, 
then the prime factors of the number $a_k$ coincide with the prime factors of the 
number $\tilde a_k$.\par
     The following sets and their associated numbers are defined in a slightly
different manner. The sets $S_{19}$, $S_{20}$, $S_{21}$, $S_{22}$ are given by
the formulas 
$$
\gather
\hskip -2em
S_{19}=\{p_i\in S:\ \alpha_i>\beta_i=\omega_i=0,\ \mu_i=\alpha_i,
\ \eta_i=0\},\\
\vspace{-1.7ex}
\mytag{6.13}\\
\vspace{-1.7ex}
\hskip -2em
S_{20}=\{p_i\in S:\ \alpha_i>\beta_i=\omega_i=0,\ \mu_i=0,
\ \eta_i=\alpha_i\},\\
\vspace{1ex}
\hskip -2em
S_{21}=\{p_i\in S:\ \omega_i=\beta_i>\alpha_i=0,\ \mu_i=\omega_i,
\ \eta_i=0\},\\
\vspace{-1.7ex}
\mytag{6.14}\\
\vspace{-1.7ex}
\hskip -2em
S_{22}=\{p_i\in S:\ \omega_i=\beta_i>\alpha_i=0,\ \mu_i=0,
\ \eta_i=\omega_i\}.
\endgather
$$
The formulas \mythetag{6.13} and \mythetag{6.14} correspond to the 
tables~\mythetable{5.7} through \mythetable{5.12} and \mythetable{5.57} 
through \mythetable{5.62}. Using them, we define the following integer 
numbers:
$$
\xalignat 2
&\hskip -2em
a_7=\prod_{p_i\in S_{19}}p_i^{\,\alpha_i},
&&\tilde a_7=\prod_{p_i\in S_{20}}p_i^{\,\alpha_i},
\mytag{6.15}\\
&\hskip -2em
a_8=\prod_{p_i\in S_{21}}p_i^{\,\omega_i},
&&\tilde a_8=\prod_{p_i\in S_{22}}p_i^{\,\omega_i},
\mytag{6.16}\\
\endxalignat
$$
Some of the sets $S_{19}$, $S_{20}$, $S_{21}$, $S_{22}$ can be empty. 
Therefore we interpret the formulas \mythetag{6.15} and \mythetag{6.16} 
so that $a_k=1$ or $\tilde a_k=1$ if the corresponding set $S_k$ is 
empty. Unlike the numbers \mythetag{6.12}, the numbers \mythetag{6.15} 
and \mythetag{6.16} are not correlated within their pairs.\par
     The sets $S_{23}$, $S_{24}$, $S_{25}$, $S_{26}$ are given by the 
following formulas 
$$
\gather
\hskip -2em
S_{23}=\{p_i\in S:\ \beta_i>\alpha_i=\omega_i=0,\ \mu_i=\beta_i,
\ \eta_i=0\},\\
\vspace{-1.7ex}
\mytag{6.17}\\
\vspace{-1.7ex}
\hskip -2em
S_{24}=\{p_i\in S:\ \beta_i>\alpha_i=\omega_i=0,\ \mu_i=0,
\ \eta_i=\beta_i\},\\
\vspace{1ex}
\hskip -2em
S_{25}=\{p_i\in S:\ \alpha_i=\omega_i>\beta_i=0,\ \mu_i=\alpha_i,
\ \eta_i=0\},\\
\vspace{-1.7ex}
\mytag{6.18}\\
\vspace{-1.7ex}
\hskip -2em
S_{26}=\{p_i\in S:\ \alpha_i=\omega_i>\beta_i=0,\ \mu_i=0,
\ \eta_i=\alpha_i\}.
\endgather
$$
The formulas \mythetag{6.17} and \mythetag{6.18} correspond to the 
tables~\mythetable{5.45} through \mythetable{5.50} and \mythetable{5.23} 
through \mythetable{5.28}. Using them, we define the following integer 
numbers:
$$
\allowdisplaybreaks
\xalignat 2
&\hskip -2em
b_7=\prod_{p_i\in S_{23}}p_i^{\,\beta_i},
&&\tilde b_7=\prod_{p_i\in S_{24}}p_i^{\,\beta_i},
\mytag{6.19}\\
&\hskip -2em
b_8=\prod_{p_i\in S_{25}}p_i^{\,\alpha_i},
&&\tilde b_8=\prod_{p_i\in S_{26}}p_i^{\,\alpha_i},
\mytag{6.20}\\
\endxalignat
$$
Some of the sets $S_{23}$, $S_{24}$, $S_{25}$, $S_{26}$ can be empty. 
Therefore we interpret the formulas \mythetag{6.19} and \mythetag{6.20} 
so that $b_k=1$ or $\tilde b_k=1$ if the corresponding set $S_k$ is 
empty. The numbers \mythetag{6.19} and \mythetag{6.20} are also not 
correlated within their pairs.\par
     The sets $S_{27}$, $S_{28}$, $S_{29}$, $S_{30}$ are given by the 
following formulas: 
$$
\gather
\hskip -2em
S_{27}=\{p_i\in S:\ \omega_i>\alpha_i=\beta_i=0,\ \mu_i=\omega_i,
\ \eta_i=0\},\\
\vspace{-1.7ex}
\mytag{6.21}\\
\vspace{-1.7ex}
\hskip -2em
S_{28}=\{p_i\in S:\ \omega_i>\alpha_i=\beta_i=0,\ \mu_i=0,
\ \eta_i=\omega_i\},\\
\vspace{1ex}
\hskip -2em
S_{29}=\{p_i\in S:\ \alpha_i=\beta_i>\omega_i=0,\ \mu_i=\alpha_i,
\ \eta_i=0\},\\
\vspace{-1.7ex}
\mytag{6.22}\\
\vspace{-1.7ex}
\hskip -2em
S_{30}=\{p_i\in S:\ \alpha_i=\beta_i>\omega_i=0,\ \mu_i=0,
\ \eta_i=\alpha_i\}.
\endgather
$$
The formulas \mythetag{6.21} and \mythetag{6.22} correspond to the 
tables~\mythetable{5.35} through \mythetable{5.38} and \mythetable{5.13} 
through \mythetable{5.16}. Using them, we define the following integer 
numbers:
$$
\xalignat 2
&\hskip -2em
u_7=\prod_{p_i\in S_{27}}p_i^{\,\omega_i},
&&\tilde u_7=\prod_{p_i\in S_{28}}p_i^{\,\omega_i},
\mytag{6.23}\\
&\hskip -2em
u_8=\prod_{p_i\in S_{29}}p_i^{\,\alpha_i},
&&\tilde u_8=\prod_{p_i\in S_{30}}p_i^{\,\alpha_i},
\mytag{6.24}\\
\endxalignat
$$
\par
     Now we can compare the formulas \mythetag{6.1}, \mythetag{6.3}, 
\mythetag{6.5}, \mythetag{6.7}, \mythetag{6.9},  \mythetag{6.11}, 
\mythetag{6.13}, \mythetag{6.14}, \mythetag{6.17}, \mythetag{6.18}, 
\mythetag{6.21}, and \mythetag{6.22} with the tables~\mythetable{5.1} 
through \mythetable{5.68} considered in the previous section. As a 
result we derive that 
$$
\hskip -2em
S=\{p_1,\,\ldots,\,p_n\}=\bigcup^{30}_{i=1}S_i.
\mytag{6.25}
$$
Then we apply \mythetag{6.25} to the formulas \mythetag{6.2}, 
\mythetag{6.4}, \mythetag{6.2}, \mythetag{6.6}, \mythetag{6.8}, 
\mythetag{6.10}, \mythetag{6.12}, \mythetag{6.15}, \mythetag{6.16}, 
\mythetag{6.19}, \mythetag{6.20}, \mythetag{6.23}, and \mythetag{6.24}.
This yields the following formulas:
$$
\allowdisplaybreaks
\gather
\hskip -2em
a=a_7\,\tilde a_7\,b_8\,\tilde b_8\,u_8\,\tilde u_8
\,\prod^6_{i=1}a_i\,\prod^3_{i=1}b_i\,\tilde b_i
\,\prod^6_{i=4}u_i\,\tilde u_i,
\mytag{6.26}\\
\hskip -2em
b=b_7\,\tilde b_7\,u_8\,\tilde u_8\,a_8\,\tilde a_8
\,\prod^6_{i=1}b_i\,\prod^3_{i=1}u_i\,\tilde u_i
\,\prod^6_{i=4}a_i\,\tilde a_i,
\mytag{6.27}\\
\hskip -2em
u=u_7\,\tilde u_7\,a_8\,\tilde a_8\,b_8\,\tilde b_8\,
\,\prod^6_{i=1}u_i\,\prod^3_{i=1}a_i\,\tilde a_i
\,\prod^6_{i=4}b_i\,\tilde b_i,
\mytag{6.28}\\
\hskip -2em
\tilde a=a_8\,\tilde a_8\,b_7\,\tilde b_7\,u_7\,\tilde u_7\,
\,\prod^6_{i=1}\tilde a_i\,\prod^3_{i=1}u_i\,\tilde u_i
\,\prod^6_{i=4}b_i\,\tilde b_i,
\mytag{6.29}\\
\hskip -2em
\tilde b=b_8\,\tilde b_8\,u_7\,\tilde u_7\,a_7\,\tilde a_7
\,\prod^6_{i=1}\tilde b_i\,\prod^3_{i=1}a_i\,\tilde a_i
\,\prod^6_{i=4}u_i\,\tilde u_i,
\mytag{6.30}\\
\hskip -2em
\tilde u=u_8\,\tilde u_8\,a_7\,\tilde a_7\,b_7\,\tilde b_7\,
\,\prod^6_{i=1}\tilde u_i\,\prod^3_{i=1}b_i\,\tilde b_i
\,\prod^6_{i=4}a_i\,\tilde a_i.
\mytag{6.31}\\
\endgather
$$
Along with the formulas \mythetag{6.26}, \mythetag{6.27}, 
\mythetag{6.28}, \mythetag{6.29}, \mythetag{6.30}, 
\mythetag{6.31}, we derive the following formulas for the
parameters $Z$, $A_0$, and $B_0$: 
$$
\gather
\hskip -2em
Z=\prod^8_{i=1}a_i\,\tilde a_i\,\prod^8_{i=1}b_i\,\tilde b_i\,
\prod^8_{i=1}u_i\,\tilde u_i,
\mytag{6.32}\\
\hskip -2em
A_0\kern 4pt=\kern -4pt\prod\Sb i=1,3,4\\i=6,7,8\endSb\kern -6pt 
a_i\,b_i\,u_i\kern -3pt\prod_{i=3,6}\kern -3pt 
\tilde a_i\,\tilde b_i\,\tilde u_i,
\mytag{6.33}\\
\hskip -2em
B_0\kern 4pt=\kern -4pt\prod\Sb i=1,2,4\\i=5,7,8\endSb\kern -6pt 
\tilde a_i\,\tilde b_i\,\tilde u_i\kern -3pt\prod_{i=2,5}\kern -3pt 
a_i\,b_i\,u_i.
\mytag{6.34}\\
\endgather
$$
Thus, the variables $a$, $b$, $u$, $\tilde a$, $\tilde b$, $\tilde u$,
$Z$, $A_0$, $B_0$ are expressed through $48$ new variables 
$a_1,\,\ldots,\,a_8$, $\tilde a_1,\,\ldots,\,\tilde a_8$, $b_1,\,\ldots,
\,b_8$, $\tilde b_1,\,\ldots,\,\tilde b_8$, $u_1,\,\ldots,\,u_8$, 
$\tilde u_1,\,\ldots,\,\tilde u_8$ by means of the formulas \mythetag{6.26},
\mythetag{6.27}, \mythetag{6.28}, \mythetag{6.29}, \mythetag{6.30},
\mythetag{6.31}, \mythetag{6.32}, \mythetag{6.33}, and \mythetag{6.34}.
Most of the new variable are coprime by their definition. Indeed, we have 
the following coprimality conditions to be fulfilled:
$$
\gather
\hskip -4em
\gcd(a_i,a_j)=1,\ \gcd(b_i,b_j)=1,\ \gcd(u_i,u_j)=1\text{\ \ for }i\neq j;\\
\hskip -4em
\gcd(\tilde a_i,\tilde a_j)=1,\ \gcd(\tilde b_i,\tilde b_j)=1,\ \gcd(\tilde u_i,
\tilde u_j)=1\text{\ \ for }i\neq j;\\
\hskip -2em
\gcd(a_i,\tilde a_j)=\gcd(b_i,\tilde b_j)=\gcd(u_i,\tilde u_j)=1\text{\ \ for }
i\neq j\text{\ or }i>6\text{\ or }j>6;\qquad\\
\vspace{-1.5ex}
\mytag{6.35}\\
\vspace{-1.5ex}
\hskip -4em
\gcd(a_i,b_j)=\gcd(a_i,\tilde b_j)=\gcd(a_i,u_j)=\gcd(a_i,\tilde u_j)=1;\\
\hskip -4em
\gcd(\tilde a_i,b_j)=\gcd(\tilde a_i,\tilde b_j)=\gcd(\tilde a_i,u_j)
=\gcd(\tilde a_i,\tilde u_j)=1;\\
\hskip -4em
\gcd(b_i,u_j)=\gcd(b_i,\tilde u_j)=\gcd(\tilde b_i,u_j)=\gcd(\tilde b_i,
\tilde u_j)=1.
\endgather
$$
The only exception from the above coprimality conditions \mythetag{6.35} are the 
numbers within the pairs $(a_i,\tilde a_i)$, $(b_i,\tilde b_i)$, $(u_i,\tilde u_i)$
for $1\leqslant i\leqslant 6$. In such pairs we have 
$$
\gather
\hskip -2em
a_i=1\Longleftrightarrow \tilde a_i=1\text{\ \ for \ }i=1,\,\ldots,\,6;\\
\hskip -2em
b_i=1\Longleftrightarrow \tilde b_i=1\text{\ \ for \ }i=1,\,\ldots,\,6;
\mytag{6.36}\\
\hskip -2em
u_i=1\Longleftrightarrow \tilde u_i=1\text{\ \ for \ }i=1,\,\ldots,\,6;\\
\vspace{1ex}
\hskip -2em
p\mid a_i\Longleftrightarrow p\mid \tilde a_i\text{\ \ for\ }
p\text{\ is prime and\ }i=1,\,\ldots,\,6;\\
\hskip -2em
p\mid b_i\Longleftrightarrow p\mid \tilde b_i\text{\ \ for\ }
p\text{\ is prime and\ }i=1,\,\ldots,\,6;
\mytag{6.37}\\
\hskip -2em
p\mid u_i\Longleftrightarrow p\mid \tilde u_i\text{\ \ for\ }
p\text{\ is prime and\ }i=1,\,\ldots,\,6.
\endgather
$$
The conditions \mythetag{6.36} and \mythetag{6.37} are called the {\it cohesion
conditions}.\par
     The next step now is to substitute the formulas \mythetag{6.26}, \mythetag{6.27}, 
\mythetag{6.28}, \mythetag{6.29}, \mythetag{6.30}, \mythetag{6.31}, \mythetag{6.32}, 
\mythetag{6.33}, and \mythetag{6.34} into the equation \mythetag{4.26}. As a result 
we get a polynomial equation with $27$ terms. These terms have the common factor
$$
\gathered
C=a_1^4\,a_2^2\,a_3^4\,a_4^4\,a_5^2\,a_6^4\,a_7^2\,a_8^4\,\tilde a_1^4\,\tilde a_2^4\,
\tilde a_3^2\,\tilde a_4^4\,\tilde a_5^4\,\tilde a_6^2\,\tilde a_7^4\,\tilde a_8^2\,
b_1^4\,b_2^2\,b_3^4\,b_4^4\,b_5^2\,b_6^4\,b_7^2\,b_8^4\,\cdot\\
\cdot\,\tilde b_1^4\,\tilde b_2^4\,\tilde b_3^2\,\tilde b_4^4\,\tilde b_5^4\,
\tilde b_6^2\,\tilde b_7^4\,\tilde b_8^2\,u_1^4\,u_2^2\,u_3^4\,u_4^4\,u_5^2\,u_6^4\,
u_7^2\,u_8^4\,\tilde u_1^4\,\tilde u_2^4\,\tilde u_3^2\,\tilde u_4^4\,\tilde u_5^4\,
\tilde u_6^2\,\tilde u_7^4\,\tilde u_8^2.
\endgathered\quad
\mytag{6.38}
$$
Even upon splitting out the nonzero common factor \mythetag{6.38} the equation 
\mythetag{4.26} appears to be rather huge. It is written as follows:
$$
\allowdisplaybreaks
\gather
a_1^4\,a_2^8\,a_3^4\,a_5^4\,a_7^4\,\tilde a_1^4\,\tilde a_2^4\,\tilde a_3^4\,
\tilde a_7^4\,\tilde a_8^4\,b_2^4\,b_5^4\,b_8^4\,\tilde b_1^4\,\tilde b_2^4\,
\tilde b_3^4\,\tilde b_4^4\,\tilde b_5^4\,\tilde b_6^4\,\tilde b_8^8\,
u_2^4\,u_4^4\,u_5^8\,u_6^4\,u_7^4\,
\tilde u_4^4\,\tilde u_5^4\,\tilde u_6^4\,\tilde u_7^4\,\tilde u_8^4\,+\\
+\ 6\,a_1^4\,a_2^6\,a_3^4\,a_4^2\,a_5^4\,a_6^2\,a_7^6\,\tilde a_1^2\,\tilde a_2^2\,
\tilde a_3^4\,\tilde a_6^2\,\tilde a_7^4\,\tilde a_8^2\,b_1^2\,b_2^4\,b_3^2\,b_5^2\,
b_7^2\,b_8^4\,\tilde b_1^4\,\tilde b_2^4\,\tilde b_3^6\,\tilde b_4^2\,\tilde b_5^2\,
\tilde b_6^4\,\tilde b_8^6\,u_2^2\,u_4^4\,u_5^6\ \cdot\\
\cdot\ u_6^4\,u_7^4\,u_8^2\,\tilde u_3^2\,\tilde u_4^4\,\tilde u_5^4\,\tilde u_6^6\,
\tilde u_7^2\,\tilde u_8^4-2\,a_1^4\,a_2^6\,a_3^4\,a_5^2\,a_7^4\,a_8^2\,
\tilde a_1^4\,\tilde a_2^4\,\tilde a_3^6\,\tilde a_6^2\,\tilde a_7^2\,\tilde a_8^4\,
b_2^2\,b_4^2\,b_5^4\,b_6^2\,b_7^2\,b_8^4\,\tilde b_1^2\ \cdot\\
\cdot\ \tilde b_2^2\,\tilde b_3^4\,\tilde b_4^4\,\tilde b_5^4\,\tilde b_6^6\,
\tilde b_8^6\,u_1^2\,u_2^4\,u_3^2\,u_4^4\,u_5^6\,u_6^4\,u_7^6\,\tilde u_3^2\,
\tilde u_4^2\,\tilde u_5^2\,\tilde u_6^4\,\tilde u_7^4\,\tilde u_8^2+a_1^4\,a_2^4\,
a_3^4\,a_4^4\,a_5^4\,a_6^4\,a_7^8\,\tilde a_3^4\,\tilde a_6^4\ \cdot\\
\cdot\ \tilde a_7^4\,b_1^4\,b_2^4\,b_3^4\,b_7^4\,b_8^4\,\tilde b_1^4\,\tilde b_2^4\,
\tilde b_3^8\,\tilde b_6^4\,\tilde b_8^4\,u_4^4\,u_5^4\,u_6^4\,u_7^4\,u_8^4\,
\tilde u_3^4\,\tilde u_4^4\,\tilde u_5^4\,\tilde u_6^8\,\tilde u_8^4
+4\,a_1^4\,a_2^4\,a_3^4\,a_4^2\,a_5^2\,a_6^2\,a_7^6\ \cdot\\
\cdot\ a_8^2\,\tilde a_1^2\,\tilde a_2^2\,\tilde a_3^6\,\tilde a_6^4\,\tilde a_7^2\,
\tilde a_8^2\,b_1^2\,b_2^2\,b_3^2\,b_4^2\,b_5^2\,b_6^2\,b_7^4\,b_8^4\,\tilde b_1^2\,
\tilde b_2^2\,\tilde b_3^6\,\tilde b_4^2\,\tilde b_5^2\,\tilde b_6^6\,\tilde b_8^4\,
u_1^2\,u_2^2\,u_3^2\,u_4^4\,u_5^4\,u_6^4\,u_7^6\,u_8^2\ \cdot\\
\cdot\ \tilde u_3^4\,\tilde u_4^2\,\tilde u_5^2\,\tilde u_6^6\,\tilde u_7^2\,
\tilde u_8^2+a_1^4\,a_2^4\,a_3^4\,a_7^4\,a_8^4\,\tilde a_1^4\,\tilde a_2^4\,
\tilde a_3^8\,\tilde a_6^4\,\tilde a_8^4\,b_4^4\,b_5^4\,b_6^4\,b_7^4\,b_8^4\,
\tilde b_3^4\,\tilde b_4^4\,\tilde b_5^4\,\tilde b_6^8\,\tilde b_8^4\,u_1^4\,
u_2^4\,u_3^4\ \cdot\\
\cdot\ u_4^4\,u_5^4\,u_6^4\,u_7^8\,\tilde u_3^4\,\tilde u_6^4\,\tilde u_7^4
-2\,a_1^4\,a_2^2\,a_3^4\,a_4^4\,a_5^2\,a_6^4\,a_7^8\,a_8^2\,\tilde a_3^6\,
\tilde a_6^6\,\tilde a_7^2\,b_1^4\,b_2^2\,b_3^4\,b_4^2\,b_6^2\,b_7^6\,b_8^4\,
\tilde b_1^2\,\tilde b_2^2\,\tilde b_3^8\ \cdot\\
\cdot\ \tilde b_6^6\,\tilde b_8^2\,u_1^2\,u_3^2\,u_4^4\,u_5^2\,u_6^4\,u_7^6\,u_8^4\,
\tilde u_3^6\,\tilde u_4^2\,\tilde u_5^2\,\tilde u_6^8\,\tilde u_8^2+6\,a_1^4\,
a_2^2\,a_3^4\,a_4^2\,a_6^2\,a_7^6\,a_8^4\,\tilde a_1^2\,\tilde a_2^2\,\tilde a_3^8\,
\tilde a_6^6\,\tilde a_8^2\,b_1^2\,b_3^2\ \cdot\\
\cdot\ b_4^4\,b_5^2\,b_6^4\,b_7^6\,b_8^4\,\tilde b_3^6\,\tilde b_4^2\,\tilde b_5^2\,
\tilde b_6^8\,\tilde b_8^2\,u_1^4\,u_2^2\,u_3^4\,u_4^4\,u_5^2\,u_6^4\,u_7^8\,u_8^2\,
\tilde u_3^6\,\tilde u_6^6\,\tilde u_7^2+a_1^4\,a_3^4\,a_4^4\,a_6^4\,a_7^8\,
a_8^4\,\tilde a_3^8\,\tilde a_6^8\ \cdot\\
\cdot\ b_1^4\,b_3^4\,b_4^4\,b_6^4\,b_7^8\,b_8^4\,\tilde b_3^8\,\tilde b_6^8\,u_1^4\,
u_3^4\,u_4^4\,u_6^4\,u_7^8\,u_8^4\,\tilde u_3^8\,\tilde u_6^8-2\,a_1^2\,a_2^8\,
a_3^2\,a_5^6\,a_7^2\,\tilde a_1^4\,\tilde a_2^4\,\tilde a_3^2\,\tilde a_4^2\,
\tilde a_5^2\,\tilde a_7^4\,\tilde a_8^6\ \cdot\\
\cdot\ b_2^6\,b_5^6\,b_8^2\,\tilde b_1^4\,\tilde b_2^4\,\tilde b_3^2\,\tilde b_4^4\,
\tilde b_5^4\,\tilde b_6^2\,\tilde b_7^2\,\tilde b_8^8\,u_2^6\,u_4^2\,u_5^8\,u_6^2\,
u_7^2\,\tilde u_1^2\,\tilde u_2^2\,\tilde u_4^4\,\tilde u_5^4\,\tilde u_6^2\,
\tilde u_7^4\,\tilde u_8^6+4\,a_1^2\,a_2^6\,a_3^2\,a_4^2\,a_5^6\ \cdot\\
\cdot\ a_6^2\,a_7^4\,\tilde a_1^2\,\tilde a_2^2\,\tilde a_3^2\,\tilde a_4^2\,
\tilde a_5^2\,\tilde a_6^2\,\tilde a_7^4\,\tilde a_8^4\,b_1^2\,b_2^6\,b_3^2\,b_5^4\,
b_7^2\,b_8^2\,\tilde b_1^4\,\tilde b_2^4\,\tilde b_3^4\,\tilde b_4^2\,\tilde b_5^2\,
\tilde b_6^2\,\tilde b_7^2\,\tilde b_8^6\,u_2^4\,u_4^2\,u_5^6\,u_6^2\,u_7^2\,
u_8^2\ \cdot\\
\cdot\ \tilde u_1^2\,\tilde u_2^2\,\tilde u_3^2\,\tilde u_4^4\,\tilde u_5^4\,
\tilde u_6^4\,\tilde u_7^2\,\tilde u_8^6-12\,a_1^2\,a_2^6\,a_3^2\,a_5^4\,a_7^2\,
a_8^2\,\tilde a_1^4\,\tilde a_2^4\,\tilde a_3^4\,\tilde a_4^2\,\tilde a_5^2\,
\tilde a_6^2\,\tilde a_7^2\,\tilde a_8^6\,b_2^4\,b_4^2\,b_5^6\,b_6^2\,b_7^2\ \cdot\\
\cdot\ b_8^2\,\tilde b_1^2\,\tilde b_2^2\,\tilde b_3^2\,\tilde b_4^4\,\tilde b_5^4\,
\tilde b_6^4\,\tilde b_7^2\,\tilde b_8^6\,u_1^2\,u_2^6\,u_3^2\,u_4^2\,u_5^6\,u_6^2\,
u_7^4\,\tilde u_1^2\,\tilde u_2^2\,\tilde u_3^2\,\tilde u_4^2\,\tilde u_5^2\,
\tilde u_6^2\,\tilde u_7^4\,\tilde u_8^4-2\,a_1^2\,a_2^4\,a_3^2\,a_4^4\ \cdot\\
\cdot\ a_5^6\,a_6^4\,a_7^6\,
\tilde a_3^2\,\tilde a_4^2\,\tilde a_5^2\,\tilde a_6^4\,\tilde a_7^4\,\tilde a_8^2\,b_1^4\,b_2^6\,b_3^4\,b_5^2\,b_7^4\,b_8^2\,\tilde b_1^4\,\tilde b_2^4\,
\tilde b_3^6\,\tilde b_6^2\,\tilde b_7^2\,\tilde b_8^4\,u_2^2\,u_4^2\,u_5^4\,u_6^2\,
u_7^2\,u_8^4\,\tilde u_1^2\,\tilde u_2^2\,\tilde u_3^4\ \cdot\\
\cdot\ \tilde u_4^4\,\tilde u_5^4\,\tilde u_6^6\,\tilde u_8^6-8\,a_1^2\,a_2^4\,
a_3^2\,a_4^2\,a_5^4\,a_6^2\,a_7^4\,a_8^2\,\tilde a_1^2\,\tilde a_2^2\,\tilde a_3^4\,
\tilde a_4^2\,\tilde a_5^2\,\tilde a_6^4\,\tilde a_7^2\,\tilde a_8^4\,b_1^2\,b_2^4\,
b_3^2\,b_4^2\,b_5^4\,b_6^2\,b_7^4\,b_8^2\ \cdot\\
\cdot\ \tilde b_1^2\,\tilde b_2^2\,\tilde b_3^4\,\tilde b_4^2\,\tilde b_5^2\,
\tilde b_6^4\,\tilde b_7^2\,\tilde b_8^4\,u_1^2\,u_2^4\,u_3^2\,u_4^2\,u_5^4\,u_6^2\,
u_7^4\,u_8^2\,\tilde u_1^2\,\tilde u_2^2\,\tilde u_3^4\,\tilde u_4^2\,\tilde u_5^2\,
\tilde u_6^4\,\tilde u_7^2\,\tilde u_8^4-2\,a_1^2\,a_2^4\,a_3^2\,a_5^2\ \cdot\\
\cdot\ a_7^2\,a_8^4\,\tilde a_1^4\,\tilde a_2^4\,\tilde a_3^6\,\tilde a_4^2\,
\tilde a_5^2\,\tilde a_6^4\,\tilde a_8^6\,b_2^2\,b_4^4\,b_5^6\,b_6^4\,b_7^4\,b_8^2\,
\tilde b_3^2\,\tilde b_4^4\,\tilde b_5^4\,\tilde b_6^6\,\tilde b_7^2\,\tilde b_8^4\,
u_1^4\,u_2^6\,u_3^4\,u_4^2\,u_5^4\,u_6^2\,u_7^6\,\tilde u_1^2\,
\tilde u_2^2\ \cdot\\
\cdot\ \tilde u_3^4\,\tilde u_6^2\,\tilde u_7^4\,\tilde u_8^2-12\,a_1^2\,a_2^2\,
a_3^2\,a_4^4\,a_5^4\,a_6^4\,a_7^6\,a_8^2\,\tilde a_3^4\,\tilde a_4^2\,\tilde a_5^2\,
\tilde a_6^6\,\tilde a_7^2\,\tilde a_8^2\,b_1^4\,b_2^4\,b_3^4\,b_4^2\,b_5^2\,b_6^2\,
b_7^6\,b_8^2\,\tilde b_1^2\,\tilde b_2^2\ \cdot\\
\cdot\ \tilde b_3^6\,\tilde b_6^4\,\tilde b_7^2\,\tilde b_8^2\,u_1^2\,u_2^2\,u_3^2\,
u_4^2\,u_5^2\,u_6^2\,u_7^4\,u_8^4\,\tilde u_1^2\,\tilde u_2^2\,\tilde u_3^6\,
\tilde u_4^2\,\tilde u_5^2\,\tilde u_6^6\,\tilde u_8^4+4\,a_1^2\,a_2^2\,a_3^2\,
a_4^2\,a_5^2\,a_6^2\,a_7^4\,a_8^4\ \cdot\\
\cdot\ \tilde a_1^2\,\tilde a_2^2\,\tilde a_3^6\,\tilde a_4^2\,\tilde a_5^2\,
\tilde a_6^6\,\tilde a_8^4\,b_1^2\,b_2^2\,b_3^2\,b_4^4\,b_5^4\,b_6^4\,b_7^6\,b_8^2\,
\tilde b_3^4\,\tilde b_4^2\,\tilde b_5^2\,\tilde b_6^6\,\tilde b_7^2\,\tilde b_8^2\,
u_1^4\,u_2^4\,u_3^4\,u_4^2\,u_5^2\,u_6^2\,u_7^6\,u_8^2\,\tilde u_1^2\ \cdot\\
\cdot\ \tilde u_2^2\,\tilde u_3^6\,\tilde u_6^4\,\tilde u_7^2\,\tilde u_8^2
-2\,a_1^2\,a_3^2\,a_4^4\,a_5^2\,a_6^4\,a_7^6\,a_8^4\,\tilde a_3^6\,\tilde a_4^2\,
\tilde a_5^2\,\tilde a_6^8\,\tilde a_8^2\,b_1^4\,b_2^2\,b_3^4\,b_4^4\,b_5^2\,b_6^4\,
b_7^8\,b_8^2\,\tilde b_3^6\,\tilde b_6^6\,\tilde b_7^2\ \cdot\\
\cdot\ u_1^4\,u_2^2\,u_3^4\,u_4^2\,u_6^2\,u_7^6\,u_8^4\,\tilde u_1^2\,\tilde u_2^2\,
\tilde u_3^8\,\tilde u_6^6\,\tilde u_8^2+a_2^8\,a_5^8\,\tilde a_1^4\,\tilde a_2^4\,
\tilde a_4^4\,\tilde a_5^4\,\tilde a_7^4\,\tilde a_8^8\,b_2^8\,b_5^8\,\tilde b_1^4\,
\tilde b_2^4\,\tilde b_4^4\,\tilde b_5^4\,\tilde b_7^4\,\tilde b_8^8\ \cdot\\
\cdot\ u_2^8\,u_5^8\,\tilde u_1^4\,\tilde u_2^4\,\tilde u_4^4\,\tilde u_5^4\,
\tilde u_7^4\,\tilde u_8^8+6\,a_2^6\,a_4^2\,a_5^8\,a_6^2\,a_7^2\,\tilde a_1^2\,
\tilde a_2^2\,\tilde a_4^4\,\tilde a_5^4\,\tilde a_6^2\,\tilde a_7^4\,\tilde a_8^6\,
b_1^2\,b_2^8\,b_3^2\,b_5^6\,b_7^2\,\tilde b_1^4\,\tilde b_2^4\,\tilde b_3^2\ \cdot\\
\cdot\ \tilde b_4^2\,\tilde b_5^2\,\tilde b_7^4\,\tilde b_8^6\,u_2^6\,u_5^6\,u_8^2\,
\tilde u_1^4\,\tilde u_2^4\,\tilde u_3^2\,\tilde u_4^4\,\tilde u_5^4\,\tilde u_6^2\,
\tilde u_7^2\,\tilde u_8^8-2\,a_2^6\,a_5^6\,a_8^2\,\tilde a_1^4\,\tilde a_2^4\,
\tilde a_3^2\,\tilde a_4^4\,\tilde a_5^4\,\tilde a_6^2\,\tilde a_7^2\,\tilde a_8^8\,
b_2^6\,b_4^2\ \cdot\\
\cdot\ b_5^8\,b_6^2\,b_7^2\,\tilde b_1^2\,\tilde b_2^2\,\tilde b_4^4\,\tilde b_5^4\,
\tilde b_6^2\,\tilde b_7^4\,\tilde b_8^6\,u_1^2\,u_2^8\,u_3^2\,u_5^6\,u_7^2\,
\tilde u_1^4\,\tilde u_2^4\,\tilde u_3^2\,\tilde u_4^2\,\tilde u_5^2\,\tilde u_7^4\,
\tilde u_8^6+a_2^4\,a_4^4\,a_5^8\,a_6^4\,a_7^4\,\tilde a_4^4\ \cdot\\
\cdot\ \tilde a_5^4\,\tilde a_6^4\,\tilde a_7^4\,\tilde a_8^4\,b_1^4\,b_2^8\,b_3^4\,
b_5^4\,b_7^4\,\tilde b_1^4\,\tilde b_2^4\,\tilde b_3^4\,\tilde b_7^4\,\tilde b_8^4\,
u_2^4\,u_5^4\,u_8^4\,\tilde u_1^4\,\tilde u_2^4\,\tilde u_3^4\,\tilde u_4^4\,
\tilde u_5^4\,\tilde u_6^4\,\tilde u_8^8
+4\,a_2^4\,a_4^2\,a_5^6\,a_6^2\ \cdot\\
\cdot\ a_7^2\,a_8^2\,\tilde a_1^2\,\tilde a_2^2\,\tilde a_3^2\,\tilde a_4^4\,
\tilde a_5^4\,\tilde a_6^4\,\tilde a_7^2\,\tilde a_8^6\,b_1^2\,b_2^6\,b_3^2\,b_4^2\,
b_5^6\,b_6^2\,b_7^4\,\tilde b_1^2\,\tilde b_2^2\,\tilde b_3^2\,\tilde b_4^2\,
\tilde b_5^2\,\tilde b_6^2\,\tilde b_7^4\,\tilde b_8^4\,u_1^2\,u_2^6\,u_3^2\,u_5^4\,
u_7^2\ \cdot\\
\cdot\ u_8^2\,\tilde u_1^4\,\tilde u_2^4\,\tilde u_3^4\,\tilde u_4^2\,\tilde u_5^2\,
\tilde u_6^2\,\tilde u_7^2\,\tilde u_8^6+a_2^4\,a_5^4\,a_8^4\,\tilde a_1^4\,
\tilde a_2^4\,\tilde a_3^4\,\tilde a_4^4\,\tilde a_5^4\,\tilde a_6^4\,\tilde a_8^8\,
b_2^4\,b_4^4\,b_5^8\,b_6^4\,b_7^4\,\tilde b_4^4\,\tilde b_5^4\,\tilde b_6^4\,
\tilde b_7^4\,\tilde b_8^4\,\ \cdot\\
\cdot\ u_1^4\,u_2^8\,u_3^4\,u_5^4\,u_7^4\,\tilde u_1^4\,\tilde u_2^4\,\tilde u_3^4\,
\tilde u_7^4\,\tilde u_8^4-2\,a_2^2\,a_4^4\,a_5^6\,a_6^4\,a_7^4\,a_8^2\,
\tilde a_3^2\,\tilde a_4^4\,\tilde a_5^4\,\tilde a_6^6\,\tilde a_7^2\,\tilde a_8^4\,
b_1^4\,b_2^6\,b_3^4\,b_4^2\,b_5^4\,b_6^2\ \cdot\\
\cdot\ b_7^6\,\tilde b_1^2\,\tilde b_2^2\,\tilde b_3^4\,\tilde b_6^2\,\tilde b_7^4\,
\tilde b_8^2\,u_1^2\,u_2^4\,u_3^2\,u_5^2\,u_7^2\,u_8^4\,\tilde u_1^4\,\tilde u_2^4\,
\tilde u_3^6\,\tilde u_4^2\,\tilde u_5^2\,\tilde u_6^4\,\tilde u_8^6+6\,a_2^2\,
a_4^2\,a_5^4\,a_6^2\,a_7^2\,a_8^4\,
\tilde a_1^2\,\tilde a_2^2\ \cdot\\
\cdot\ \tilde a_3^4\,\tilde a_4^4\,\tilde a_5^4\,\tilde a_6^6\,\tilde a_8^6\,b_1^2\,
b_2^4\,b_3^2\,b_4^4\,b_5^6\,b_6^4\,b_7^6\,\tilde b_3^2\,\tilde b_4^2\,\tilde b_5^2\,
\tilde b_6^4\,\tilde b_7^4\,\tilde b_8^2\,u_1^4\,u_2^6\,u_3^4\,u_5^2\,u_7^4\,u_8^2\,
\tilde u_1^4\,\tilde u_2^4\,\tilde u_3^6\,\tilde u_6^2\,\tilde u_7^2\ \cdot\\
\cdot\ \tilde u_8^4+a_4^4\,a_5^4\,a_6^4\,a_7^4\,a_8^4\,\tilde a_3^4\,\tilde a_4^4\,
\tilde a_5^4\,\tilde a_6^8\,\tilde a_8^4\,b_1^4\,b_2^4\,b_3^4\,b_4^4\,b_5^4\,b_6^4\,
b_7^8\,\tilde b_3^4\,\tilde b_6^4\,\tilde b_7^4\,u_1^4\,u_2^4\,u_3^4\,u_7^4\,u_8^4\,
\tilde u_1^4\ \cdot\\
\cdot\ \tilde u_2^4\,\tilde u_3^8\,\tilde u_6^4\,\tilde u_8^4=0.
\endgather
$$
The above equation is called the {\it structural equation}. It is a Diophantine
equation with respect to $48$ integer variables $a_1,\,\ldots,\,a_8$, $\tilde a_1,
\,\ldots,\,\tilde a_8$, $b_1,\,\ldots,\,b_8$, $\tilde b_1,\,\ldots,\,\tilde b_8$, 
$u_1,\,\ldots,\,u_8$, $\tilde u_1,\,\ldots,\,\tilde u_8$. Now, summarizing the
results of the sections 5 and 6, then applying the lemma~\mythelemma{4.3}, we 
derive the following theorem. 
\mytheorem{6.1} For a given triple of positive coprime integer numbers $a$, $b$, 
and $u$ such that none of the conditions \mythetag{1.2} is satisfied the polynomial 
Diophantine equation \mythetag{1.3} is resolvable if and only if there are $48$ 
positive integer numbers $a_1,\,\ldots,\,a_8$, $\tilde a_1,\,\ldots,\,\tilde a_8$, 
$b_1,\,\ldots,\,b_8$, $\tilde b_1,\,\ldots,\,\tilde b_8$, $u_1,\,\ldots,\,u_8$, 
$\tilde u_1,\,\ldots,\,\tilde u_8$ obeying the structural equation on the pages 
32 and 33, obeying the cohesion conditions \mythetag{6.36} and \mythetag{6.37}, 
obeying the coprimality conditions \mythetag{6.35}, and such that $a$, $b$, and 
$u$ are expressed through them by means of the formulas \mythetag{6.26},
\mythetag{6.27}, \mythetag{6.28}. Under these conditions the equation 
\mythetag{1.3} has at least two solutions given by the formulas 
$$
\xalignat 2
&\hskip -2em
t\kern 4pt=\kern -4pt\prod\Sb i=1,3,4\\i=6,7,8\endSb\kern -6pt 
a_i\,b_i\,u_i\kern -3pt\prod_{i=3,6}\kern -3pt 
\tilde a_i\,\tilde b_i\,\tilde u_i,
&&t\kern 4pt=\,-\kern -4pt\prod\Sb i=1,3,4\\i=6,7,8\endSb\kern -6pt 
a_i\,b_i\,u_i\kern -3pt\prod_{i=3,6}\kern -3pt 
\tilde a_i\,\tilde b_i\,\tilde u_i.\quad
\mytag{6.39}
\endxalignat
$$
\endproclaim
     The theorem~\mythetheorem{6.1} is the required structural theorem for
the solutions of the Diophantine equation \mythetag{1.3}. The formulas 
\mythetag{6.39} in this structural theorem are immediate from the formulas 
\mythetag{6.33} and \mythetag{4.2}.\par
\head
7. Conclusions. 
\endhead
     The structural theorem~\mythetheorem{6.1} is the main result of this
paper. It can be used in computer search for perfect Euler cuboids or maybe 
in proving their non-existence in the case of the third cuboid 
conjecture~\mytheconjecture{1.1}. The theorem~\mythetheorem{6.1} is
analogous to the structural theorem~4.1 from \mycite{3} associated with the
second cuboid conjecture. 

\enddocument
\end